\documentclass{article}

\usepackage{amsmath}
\usepackage{amsfonts}
\usepackage{amssymb}
\usepackage{amsthm}
\usepackage{diagrams}
\usepackage{graphicx}

\newcommand{\fancy}{\mathcal}

\newcommand{\ZZ}{\mathbb{Z}}
\newcommand{\QQ}{\mathbb{Q}}

\newtheorem{lemma}{Lemma}
\newtheorem{theorem}{Theorem}
\newtheorem{corollary}{Corollary}
\newtheorem{conjecture}{Conjecture}
\newtheorem{proposition}{Proposition}
\newtheorem{question}{Question}

\theoremstyle{definition}
\newtheorem{definition}{Definition}
\newtheorem{example}{Example}
\newtheorem{remark}{Remark}

\numberwithin{theorem}{section}
\numberwithin{lemma}{section}
\numberwithin{proposition}{section}
\numberwithin{definition}{section}
\numberwithin{remark}{section}
\numberwithin{corollary}{section}
\numberwithin{example}{section}
\numberwithin{conjecture}{section}

\DeclareMathOperator{\cc}{cc}

\DeclareMathOperator{\ch}{ch}

\DeclareMathOperator{\cyc}{cyc}
\DeclareMathOperator{\cw}{cw}
\DeclareMathOperator{\dinv}{dinv}

\DeclareMathOperator{\ID}{ID}
\DeclareMathOperator{\inv}{inv}
\DeclareMathOperator{\invcode}{invcode}
\DeclareMathOperator{\InvPlot}{InvPlot}
\DeclareMathOperator{\maj}{maj}
\DeclareMathOperator{\majcode}{majcode}

\DeclareMathOperator{\Standardize}{Standardize}
\DeclareMathOperator{\stat}{stat}

\newcommand{\tH}{\widetilde{H}}

\newtheorem*{lemmaDuplicates}{Lemma \ref{duplicates}}
\newtheorem*{GeneralRecursion}{Proposition \ref{recursion}}
\newtheorem*{InvcodeTheorem}{Theorem \ref{Invcode}}
\newtheorem*{OneColumnStandardize}{Proposition \ref{OneColumnStandardize}}
\newtheorem*{MainLemma}{Proposition \ref{ZeroBump}}
\newtheorem*{PullUpLemma}{Lemma \ref{PullUp}}
\newtheorem*{ThreeRowsLemma}{Lemma \ref{threerows}}
\newtheorem*{ThreeRowsBijectiveLemma}{Lemma \ref{threerowsbijective}}
\newtheorem*{ThreeRowsConclusionTheorem}{Theorem \ref{threerowsconclusion}}
\newtheorem*{StandardInvCorollary}{Corollary \ref{standardinv}}
\newtheorem*{RectangleBumpCorollary3}{Proposition \ref{RectBumpCor3}}

\title{A combinatorial approach to the $q,t$-symmetry relation in Macdonald polynomials}
\author{Maria Monks Gillespie \\ \textit{University of California, Berkeley, CA 94720} \\  monks@math.berkeley.edu}

\begin{document}

\renewcommand{\thefootnote}{\fnsymbol{footnote}}
\renewcommand{\thefootnote}{\arabic{footnote}}

\maketitle{}

\begin{abstract}
Using the combinatorial formula for the transformed Macdonald polynomials of Haglund, Haiman, and Loehr, we investigate the combinatorics of the symmetry relation $\tH_\mu(\mathbf{x};q,t)=\tH_{\mu^\ast}(\mathbf{x};t,q)$. We provide a purely combinatorial proof of the relation in the case of Hall-Littlewood polynomials ($q=0$) when $\mu$ is a partition with at most three rows, and for the coefficients of the square-free monomials in $\mathbf{x}$ for all shapes $\mu$.  We also provide a proof for the full relation in the case when $\mu$ is a hook shape, and for all shapes at the specialization $t=1$.  Our work in the Hall-Littlewood case reveals a new recursive structure for the cocharge statistic on words.
\end{abstract}

\section{Introduction}

Let $\Lambda_{q,t}(x)$ denote the ring of symmetric polynomials in the countably many indeterminates $x_1,x_2,\ldots,$ with coefficients in the field $\QQ(q,t)$ of rational functions in two variables.  The \textit{(transformed) Macdonald polynomials} $\tH_{\mu}(x;q,t)\in \Lambda_{q,t}(x)$, indexed by the set of all partitions $\mu$, form an orthogonal basis of $\Lambda_{q,t}(x)$, and have specializations $\tH_{\mu}(x;0,1)=h_\mu$ and $\tH_{\mu}(x;1,1)=e_1^n$, where $h_\lambda$ and $e_\lambda$ are the homogeneous and elementary symmetric functions, respectively. The polynomials $\tH_\mu$ are a transformation of the functions $P_\lambda$ originally defined by Macdonald in \cite{Macdonald}, and have been the subject of much recent attention in combinatorics and algebraic geometry.  (See \cite{HaglundBook}, \cite{HHL}, and  \cite{CDM}, for instance.)

The symmetric functions $\tH_\mu$ may be defined as the unique collection of polynomials that satisfy certain triangularity conditions. To state them, recall that the Schur functions $s_\lambda$ form a basis for $\Lambda$.  Define the \textit{dominance order} to be the partial order $\le$ on partitions given by $\lambda>\mu$ if and only if $\lambda_1+\cdots+\lambda_k\ge \mu_1+\cdots+\mu_k$ for all $k>0$.  Finally, define $\mu^\ast$ to be the conjugate of a given partition $\mu$, formed by reflecting its Young diagram about the diagonal.

\begin{definition}
  The symmetric functions $\tH_\lambda(x;q,t)$ are the unique elements of $\Lambda_{q,t}(x)$ satisfying:
  \begin{itemize}
    \item $\tH_\mu[(1-q)X;q,t]\in\QQ(q,t)\{s_\lambda:\lambda\ge \mu\}$
    \item $\tH_\mu[(1-t)X;q,t]\in\QQ(q,t)\{s_\lambda:\lambda\ge \mu^\ast\}$
    \item $\tH_\mu[1;q,t]=1$
  \end{itemize}
  Above, the notation $(1-q)X$ means that we substitute for the variables $x_1,x_2,\ldots$ the monomials $x_1,-qx_1,x_2,-qx_2,\ldots$. 
\end{definition}

 The Macdonald polynomials $\tH_\mu$ are orthogonal with respect to the inner product on $\Lambda_{q,t}$ defined by $$\langle f,g\rangle_{q,t}=\left\langle f(x), g\left[\frac{1-q}{1-t}X\right]\right\rangle,$$ where the inner product on the right is the classical Hall inner product.  That is, $$\langle\tH_\mu,\tH_{\lambda}\rangle_{q,t}=0$$ whenever $\mu\neq \lambda$.  (See \cite{CDM} for details.)

 Recall the well-known Schur expansion $$h_\mu=\sum K_{\lambda\mu}s_\lambda$$ where the coefficients $K_{\lambda\mu}$  are the \textit{Kostka numbers}, defined combinatorially as the number of semistandard Young tableaux with shape $\lambda$ and content $\mu$.  Since $\tH_\mu(x;0,1)=h_\mu$, it is natural to define a $q,t$-analog of the Kostka numbers by expanding the transformed Macdonald polynomials $\widetilde{H}_\mu(x;q,t)$ in terms of the Schur basis.

\begin{definition}
  The \textit{$q,t$-Kostka polynomials} are the coefficients in the expansion $$\widetilde{H}_{\mu}(x;q,t)=\sum_{\lambda}\widetilde{K}_{\lambda\mu}(q,t)s_\lambda$$
\end{definition}

It was conjectured by Macdonald, and later proven by Haiman \cite{Haiman}, that the $q,t$-Kostka polynomials $\widetilde{K}_{\lambda\mu}(q,t)$ are polynomials in $q$ and $t$ with nonnegative integer coefficients.  This fact is known as the \textit{Macdonald positivity conjecture}.  Haiman's proof involves showing that the polynomial $\widetilde{K}_{\lambda\mu}(q,t)$ is the Hilbert series of a certain bi-graded module arising from the geometry of the Hilbert scheme of $n$ points in the plane, and relies heavily on geometric methods.  The problem of finding a purely combinatorial explanation of their positivity is still open, in the sense that there is no known formula for the coefficients of the form $\widetilde{K}_{\lambda\mu}(q,t)=\sum_T q^{s(T)}t^{r(T)}$, where $T$ ranges over an appropriate set of Young tableaux and $r$ and $s$ are some combinatorial statistics.

However, a different combinatorial formula for the transformed Macdonald polynomials $\widetilde{H}_\mu$ has been found, and appeared in the literature in \cite{HHL} in 2004.  The authors prove that
\begin{equation}\label{CombForm}
\widetilde{H}_\mu(x;q,t)=\sum_{\sigma} q^{\inv(\sigma)}t^{\maj(\sigma)}x^{\sigma},
\end{equation} 
 where the sum ranges over all fillings $\sigma$ of the diagram of $\mu$ with positive integers, and $x^{\sigma}$ is the monomial $x_1^{m_1}x_2^{m_2}\cdots$ where $m_i$ is the number of times the letter $i$ occurs in $\sigma$.  The statistics $\inv$ and $\maj$ are generalizations of the Mahonian statistics $\inv$ and $\maj$ for permutations.  Their precise definitions can be stated as follows.

\begin{definition}
  Given a word $w=w_1\cdots w_n$ where the letters $w_i$ are taken from some totally ordered alphabet $A$, a \textit{descent} of $w$ is an index $i$ for which $w_i>w_{i+1}$.  The \textit{major index} of $w$, denoted $\maj(w)$, is the sum of the descents of $w$.  
\end{definition}

\begin{definition}  
  Given a filling $\sigma$ of a Young diagram of shape $\mu$ drawn in French notation, let $w^{(1)},\ldots,w^{(\mu_1)}$ be the words formed by the successive columns of $\sigma$, read from top to bottom.  Then $$\maj(\sigma)=\sum_s \maj(w^{(s)}).$$
\end{definition}

\begin{example}
  The major index of the filling in Figure \ref{examplefilling} is $7$, since the first column has major index $6$, the second has major index $0$, and the third column, $1$.  
  \begin{figure}[t] 
    \begin{center}
     \includegraphics{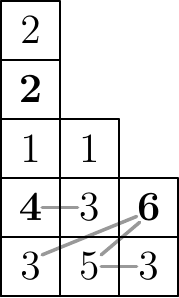}
    \end{center}
    \caption{\label{examplefilling} A filling of a Young diagram.  Descents are shown in boldface, and attacking pairs are connected with a gray line.}
  \end{figure}
\end{example}

\begin{remark}
  The major index restricts to the usual major index on words in the case that the partition is a single column.
\end{remark}

For the statistic $\inv$, we start with the definition provided in \cite{HHL}.  We use the notion of the \textit{arm} of an entry, which is defined to be the number of squares strictly to the right of the entry.  A \textit{descent} is an entry which is strictly greater than the entry just below it.

\begin{definition}
  An \textit{attacking pair} in a filling $\sigma$ of a Young diagram is a pair of entries $u$ and $v$ with $u>v$ satisfying one of the following conditions:
  \begin{enumerate}
  \item $u$ and $v$ are in the same row, with $u$ to the left of $v$, or
  \item $u$ is in the row above $v$ and strictly to its right.
  \end{enumerate}
  \end{definition}
  
 \begin{definition}
  The quantity $\inv(\sigma)$ is defined to be the number of attacking pairs in $\sigma$ minus the sum of the arms of the descents.
\end{definition}

\begin{example}
  In Figure \ref{examplefilling}, there are $4$ attacking pairs, and the arms of the descents have lengths $0$, $2$, and $0$.  Thus $\inv(\sigma)=4-2=2$ in this case.
\end{example}

For our purposes, we will also need the following cleaner definition of the $\inv$ statistic.  This more closely resembles the $\inv$ statistic on a permutation $\pi$, defined to be the number of pairs $i<j$ for which $\pi(i)>\pi(j)$.

\begin{definition}
   Let $\sigma$ be any filling of a Young diagram with letters from a totally ordered alphabet $A$, allowing repeated letters.  A \textit{relative inversion} of a filling $\sigma$ of a Young diagram is a pair of entries $u$ and $v$ in the same row, with $u$ to the left of $v$, such that if $b$ is the entry directly below $u$, one of the following conditions is satisfied:
   \begin{itemize}
     \item $u<v$ and $b$ is between $u$ and $v$ in size, in particular $u\le b<v$.
     \item $u>v$ and $b$ is not between $u$ and $v$ in size, in particular either $b<v<u$ or $v<u\le b$,
   \end{itemize}
   If $u$ and $v$ are on the bottom row, we treat $b$ as any value less than $\min(u,v)$, usually $0$ in the case $A=\ZZ_+$.
\end{definition}

\begin{remark}
  The conditions above for the triple $(u,v,b)$ to form an inversion can also be thought of as saying that the ordering of the sizes of $u,b,v$ orients the triple counterclockwise: either $b<v<u$, $v<u\le b$, or $u\le b<v$.
\end{remark}

\begin{example}
  In Figure \ref{examplefilling}, there are 2 relative inversions: $(5,3)$ in the bottom row, and $(3,6)$ in the second row.
\end{example}

In fact, the number of relative inversions in a filling $\sigma$ is always equal to $\inv(\sigma)$.  In \cite{HHL}, the authors introduce the related notion of an \textit{inversion triple}.  Relative inversions are simply the inversion triples that contribute $1$ to $\inv(\sigma)$.  The description in terms of relative inversions allows us to think of the $\inv$ as being computed row by row (just as $\maj$ is computed column by column).

For completeness, we include here a proof that $\inv(\sigma)$ is equal to the number of relative inversions of $\sigma$.

\begin{proposition}\label{inversions}
  The quantity $\inv(\sigma)$ is equal to the number of relative inversions of $\sigma$.
\end{proposition}

\begin{proof}
  Recall that $\inv(\sigma)$ is defined as the total number of attacking pairs minus the arms of the descents.  Each descent of the form $u>b$ where $b$ is the entry directly below $u$ contributes $-1$ towards $\inv(\sigma)$ for each $v$ to the right of $u$ in the same row.  Call such pairs $(u,v)$ \textit{descent-arm pairs}.  Each attacking pair contributes $+1$ towards $\inv(\sigma)$.  
  
  Define a \textit{good triple} to be a triple of entries $(u,v,b)$ where $u$ is directly above and adjacent to $b$ and $v$ is to the right of $u$ in its row, where we also allow $b$ to be directly below the entire tableau with a value of $0$.  Then each descent-arm pair or attacking pair is a member of a unique good triple, and contributes $-1$ or $+1$, respectively, to $\inv(\sigma)$.  Therefore, $\inv(\sigma)$ is the sum of the contributions of all such pairs in each such triple.
  
  A simple case analysis shows that each good triple contributes a total of $1$ if it is a relative inversion and $0$ otherwise.  Thus $\inv(\sigma)$ is the total number of relative inversions.
\end{proof}

Since this combinatorial formula for $\widetilde{H}_\mu(x;q,t)$ is an expansion in terms of monomials rather than Schur functions, it does not give an immediate answer to the Macdonald positivity conjecture.  Indeed, it perhaps raises more questions than it answers.  For one, there is a well-known $q,t$-symmetry relation for the transformed Macdonald polynomials $\widetilde{H}_\mu(x;q,t)$, namely $$\widetilde{H}_\mu(x;q,t)=\widetilde{H}_{\mu^{\ast}}(x;t,q).$$  This is obvious from the triangularity conditions that define $\tH_{\mu}$, and is also clear from Haiman's geometric interpretation \cite{Haiman}.  When combined with the combinatorial formula, however, we obtain a remarkable generating function identity:
\begin{equation}\label{SymmetryEqn}
\sum_{\sigma:\mu\to \ZZ_+} q^{\inv(\sigma)}t^{\maj(\sigma)}x^{\sigma}=\sum_{\rho:\mu^{\ast}\to \ZZ_+}q^{\maj(\rho)}t^{\inv(\rho)}x^{\rho}.
\end{equation}

  Setting $t=1$ and $\mu=(n)$ and taking the coefficient of $x_1\cdots x_n$ on both sides, this reduces to the well-known equation $$\sum_{w\in S_n}q^{\inv(w)}=\sum_{w\in S_n}q^{\maj(w)},$$ which demonstrates the equidistribution of the Mahonian statistics $\inv$ and $\maj$ on permutations.  There are several known bijective proofs of this identity (see \cite{Carlitz}, \cite{Foata}, \cite{Skandera}).

  In light of this, it is natural to ask if there is an elementary combinatorial proof of (\ref{SymmetryEqn}), in the sense of Conjecture \ref{InvMaj} below.

\begin{definition}
  The \textit{content} of a filling $\sigma$, denoted $|\sigma|$, is the sequence $\alpha=(\alpha_1,\cdots,\alpha_k)$ where $\alpha_i$ is the number of $i$'s used in the filling.  We also define the symbols:
  \begin{itemize}
  \item $\fancy{F}$ - set of all fillings of Young diagrams with positive integers
  \item $\fancy{F}_\mu^\alpha$ - set of fillings of shape $\mu$ and content $\alpha$
  \item $\fancy{F}_\mu^\alpha|_{\inv=a,\maj=b}$ - set of fillings $\sigma\in \fancy{F}_\mu^\alpha$  for which $\inv(\sigma)=a$ and $\maj(\sigma)=b$.
  \end{itemize}
  We also define a \textit{weighted set} to be a set $S$ equipped with a number of statistics $\stat_1,\stat_2,\ldots$, and a \textit{morphism of weighted sets} to be a map that preserves their statistics.  We write $$(S;\stat_1,\stat_2,\ldots)$$ to denote the weighted set if the statistics are not understood.
\end{definition}

\begin{conjecture}\label{InvMaj}
    There is a natural isomorphism of weighted sets $$\varphi:(\fancy{F};\inv,\maj)\to (\fancy{F};\maj,\inv)$$ which interchanges $\inv$ and $\maj$ and sends a partition shape to its conjugate.  That is, for any $a,b,\mu,\alpha$, the map $\varphi$ restricts to a bijection $$\varphi:\fancy{F}_\mu^\alpha|_{\inv=a,\maj=b}\to \fancy{F}_{\mu^\ast}^\alpha|_{\inv=b,\maj=a}.$$
\end{conjecture}

\begin{remark}
  In \cite{HHL}, the authors give a combinatorial proof of the fact that the polynomials $\widetilde{H}_\mu$ are symmetric in the variables $x_i$.  We will make use of this fact repeatedly, rearranging the entries of $\alpha$ as needed.  In other words, to prove Conjecture \ref{InvMaj}, it suffices to find a map $\varphi$ that restricts to bijections $\fancy{F}_\mu^\alpha|_{\inv=a,\maj=b}\to \fancy{F}_{\mu^\ast}^{r(\alpha)}|_{\inv=b,\maj=a}$ where $r$ is some bijective map that rearranges the entries of $\alpha$.
\end{remark}

In this paper, we provide explicit bijections $\varphi$ for several infinite families of values of $a$, $b$, $\alpha$, and $\mu$.  Our bijections naturally extend Carlitz's bijection on permutations, which is defined in section \ref{Carlitz}.  In Section \ref{t=1} we proceed to give a combinatorial proof of the symmetry relation for the specialization $t=1$, and in Section \ref{HookShapes}, we give an explicit bijection $\varphi$ in the case that $\mu$ is a hook shape.

The bulk of our results are developed in Section \ref{q=0}.  Here we investigate the Hall-Littlewood specialization $a=0$, which corresponds to setting $q=0$ in the Macdonald polynomials.  We give a combinatorial proof in this case for all shapes $\mu$ having at most three rows, and also for all shapes $\mu$ when the content $\alpha$ is fixed to be $(1,1,\ldots,1)$.   We also conjecture a strategy for the general problem that draws on the work of Garsia and Procesi on the $S_n$-modules $R_\mu$, which arise as the cohomology rings of the Springer fibers in type A.  \cite{GarsiaProcesi}  

In Section \ref{Applications}, we state some applications of the results on the Hall-Littlewood case to understanding the rings $R_\mu$, in particular regarding the \textit{cocharge} statistic of Lascoux and Schutzenberger (see \cite{GarsiaProcesi} or \cite{CDM}, for instance).  In particular, we demonstrate a new recursive structure exhibited by the cocharge statistic on words.

  Particularly technical proofs of results throughout the paper are deferred to Section \ref{Proofs}.

 \section{The Carlitz bijection}\label{Carlitz}

  Our approach to the symmetry problem is motivated by Carlitz's bijection $(S_n;\inv)\to (S_n;\maj)$, an alternative to the better-known Foata bijection that demonstrates the equidistribution of $\inv$ and $\maj$ on permutations.  A full proof of this bijection can be found in Carlitz's original paper \cite{Carlitz}, or in a somewhat cleaner form in \cite{Skandera}.  For the reader's convenience we will define it here.  
  
  The bijection makes use of certain \textit{codes}:
  
  \begin{definition}
    A \textit{Carlitz code} of length $n$ is a word $w=w_1\cdots w_n$ consisting of nonnegative integers such that $w_{n-i}<i$ for all $i$.  Let $C_n$ denote the set of all Carlitz codes of length $n$, equipped with the combinatorial statistic $\Sigma$ taking a word to the sum of its entries.
  \end{definition}
  
  Notice that the number of Carlitz codes of length $n$ is equal to $n!$.  This allows us to make use of the combinatorial object $(C_n;\sigma)$ of Carlitz codes as an intermediate object connecting $(S_n;\inv)$ to $(S_n;\maj)$.   In particular, the Carlitz bijection is the composite 
  \begin{diagram}
    (S_n;\inv) &\rTo^{\invcode}& (C_n;\Sigma) &\rTo^{\majcode^{-1}}& (S_n;\maj)
  \end{diagram}
   of two simple isomorphisms of weighted sets, defined as follows.
  
  \begin{definition}
  The \textit{inversion code} of a permutation $\pi$, denoted $\invcode(\pi)$, is the sequence $c_1,c_2,\ldots,c_n$ where $c_i$ is the number of inversions of the form $(j,i)$ in $\pi$, i.e. where $i<j$ and $i$ is to the right of $j$.
  \end{definition}
  
  \begin{example}
    We have $\invcode(4132)=1210$, because the $1$ is the smaller entry of one inversion $(4,1)$, the $2$ is the smaller entry of the two inversions $(3,2)$ and $(4,2)$, the $3$ is the smaller entry of the inversion $(4,3)$, and the $4$ is not the smaller entry of any inversion.
  \end{example}
  
  Clearly $\invcode$ is a map $S_n\to C_n$, and it is not hard to see that it is bijective: given a Carlitz code $c_1,\ldots,c_n$, we can reconstruct the permutation $\pi$ it came from as follows.  First write down the number $n$, corresponding to $c_n=0$.  Then, $c_{n-1}$ is either $0$ or $1$, and respectively determines whether to write down $n-1$ to the left or to the right of the $n$.  The entry $c_{n-2}$ then determines where to insert $n-2$ in the sequence, and so on until we have reconstructed $\pi$.  It is also clear that $\invcode$ is an isomorphism of weighted sets, sending $\inv(\pi)$ to $\Sigma(\{c_i\})$.
 
  \begin{definition}
   The map $\majcode:S_n\to C_n$ is defined as follows.  Given $\pi\in S_n$ written as a permutation in word form, remove the $n$ from $\pi$ and set $c_1$ to be the amount the major index decreases as a result.  Then remove the $n-1$ and set $c_2$ to be the amount the major index decreases by, and so on until we have formed a sequence $c_1,c_2,\ldots,c_n$.  Then we define $\majcode(\pi)=c_1,c_2,\ldots,c_n$.
  \end{definition}
  
  \begin{example}
    Let $\pi=3241$.  Its major index is $1+3=4$.  Removing the $4$ results in the permutation $321$, which has major index $3$, so the major index has decreased by $1$ and we set $c_1=1$.  Removing the $3$ results in $21$, which decreased the major index by $2$.  Hence $c_2=2$.  Removing the $2$ decreases the major index by $c_3=1$, and removing the $1$ decreases it by $c_4=0$, so $\majcode(\pi)=1210$.
  \end{example}
 
   As in the case of $\invcode$ above, it is not hard to construct an inverse for $\majcode$, making it an isomorphism of weighted sets $(S_n;\maj)\to (C_n,\Sigma)$.
   
   \begin{definition}
     The \textit{Carlitz bijection} is the isomorphism $$\majcode^{-1}\circ \invcode:(S_n;\inv)\to (S_n;\maj).$$
   \end{definition}
   
   \begin{example}
     We have $\majcode^{-1}\circ \invcode(4132)=\majcode^{-1}(1210)=3241$ by the examples above.
   \end{example}

  \subsection{Carlitz bijection on words}\label{CarlitzWords}
  
  Notice that the Carlitz bijection gives rise to a bijection $\phi$ satisfying Conjecture \ref{InvMaj} for one-column shapes $\mu=(1,1,\ldots,1)$ having content $\alpha=(1,1,\ldots,1)$.  Indeed, $\inv(\sigma)=0$ for any filling $\sigma$ of a one-column shape $\mu$, and $\maj(\rho)=0$ for any filling $\rho$ of its one-row conjugate $\mu^\ast$.  Since $\maj(\sigma)$ and $\inv(\rho)$ in this case are the same as $\maj$ and $\inv$ of their reading words, this determines a bijection for distinct entries ($\alpha=(1,1,\ldots,1)$.)
  
  We now generalize the Carlitz bijection to words, i.e. fillings with \textit{any} content $\alpha$ for one-column shapes $\mu$.
  
  \begin{definition}\label{Aweaklyincreasing}
   Let $A=(a_1^{\alpha_1}, a_2^{\alpha_2},\ldots, a_k^{\alpha_k})$ be any finite multiset of size $n$, with an ordering ``$<$'' such that $a_1<a_2<\cdots<a_k$, and let $\mu$ be a partition of $n$.  We say that a word $c$ of length $n$ is \textit{$A$-weakly increasing} if every subword of the form $$c_{\alpha_1+\cdots +\alpha_i},c_{\alpha_1+\cdots \alpha_i+1},\ldots c_{\alpha_1+\cdots+\alpha_i+\alpha_{i+1}-1}$$ is weakly increasing.   
  \end{definition}
  
  For instance, if $A=\{1,1,2,3,3,3,4,4\}$, ordered by magnitude, then the word $23711213$ is $A$-weakly increasing, since the subwords $23$, $7$, $112$, and $13$, corresponding to each letter of $A$, are weakly increasing.
  
  We also will make use of Macdonald symmetry in the variables $x_i$ by defining a weight-preserving bijection on alphabets.
  
  \begin{definition}\label{reverse}
   The \textit{reverse} of the content $\alpha=(\alpha_1,\ldots,\alpha_M)$ is the tuple $$r(\alpha)=(\alpha_M,\alpha_{M-1},\ldots,\alpha_1).$$  In terms of alphabets, let $A$ be a finite multiset of positive integers with maximum element $M$.  The \textit{content} of $A$ is $\alpha$ if $\alpha_i$ is the multiplicity of $i$ in $A$. The \textit{complement} of $A$, denoted $\overline{A}$, is the multiset consisting of the elements $M+1-a$ for all $a\in A$.  Notice that the content of $\overline{A}$ is $r(\alpha)$. 
  \end{definition}
  
  For instance, the complement of the multiset $$\{1,2,2,2,2,3,4,4\}$$ is $\{4,3,3,3,3,2,1,1\},$ and correspondingly, $r(1,4,1,2)=(2,1,4,1)$.
  
  We generalize Carlitz's codes as follows.
  
  \begin{definition}
    Let $C_{(1^n),A}$ denote the subset of $C_n$ consisting of all Carlitz codes of length $n$ which are $A$-weakly increasing.  This subset inherits the $\Sigma$ statistic from $C_n$.
  \end{definition}
  
  We now can define bijections $$\invcode:(\mathcal{F}_{(1^n)}^\alpha;\inv)\to (C_{(1^n),A};\Sigma)$$ and $$\majcode:(\mathcal{F}_{(n)}^{r(\alpha)};\maj)\to (C_{(1^n),A};\Sigma).$$
  
  \begin{definition}
    Let $w$ be a word consisting of the letters in the ordered alphabet $A=a_1\le \cdots \le a_n$ (corresponding to a filling of a horizontal shape), with ties among the letters broken in the order they appear in $w$.  The \textit{inversion code} of $w$ is the code $\invcode(w)=c_1\cdots c_n$ where $c_i$ is the number of inversions having $a_i$ as the smaller entry of the inversion.
  \end{definition}
  
  For example, the inversion code of the filling
  \begin{center}
  \includegraphics{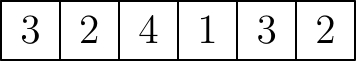}
  \end{center}
   is $312010$, since the $1$ is the smaller entry of $3$ inversions, the first $2$ is the smaller entry of $1$ inversion, the second $2$ is the smaller entry of $2$ inversions, and so on.
  
  \begin{proposition}
    The map $\invcode$ is an isomorphism of weighted sets $$\invcode:\mathcal{F}_{(1^n)}^\alpha\to C_{(1^n),A}.$$
  \end{proposition}
  
  The above proposition will be implied by Proposition \ref{Invcode}, and so we omit the proof.
  
  To define the map $\majcode$, we first require a standardization rule for fillings of columns.
  
  \begin{definition}
    Let $\sigma$ be any filling of a column of height $n$ with positive integers.     We define the \textbf{standardization labeling} on repeated entries as follows.  
    \begin{enumerate}
     \item Let $i$ be a letter that occurs $k$ times in $\sigma$.  Remove any entries larger than $i$ to form a smaller column $\sigma'$.
     \item Find the bottommost $i$ that is either on the very bottom of $\sigma'$ or has entries $a$ and $b$ above and below it with $a>b$.  Assign this $i$ a label of $k$ and remove it.  Repeat this process, labeling the next $i$ by $k-1$ and so on, until there are no $i$'s left that satisfy this condition.
     \item Finally, remove and label any remaining $i$'s in order from top to bottom, decreasing the label by one each time.
   \end{enumerate}       
      
  We define $\Standardize(\sigma)$ is the unique column filling using labels $1,2,\ldots,n$ that respects the ordering of the entries of $\sigma$ and breaks ties according to the standardization labeling.
    \end{definition}
    
  \begin{proposition}\label{OneColumnStandardize}
     For any column filling $\sigma$ with alphabet $A$, let $\rho=\Standardize(\sigma)$.  Then $\rho$ and $\sigma$ have the same major index, and $\majcode(\rho)$ is $A$-weakly increasing.
  \end{proposition}
    
    We defer the proof to Section \ref{Proofs}.  The key step is the following technical lemma.  Define a \textit{consecutive block} of $n$'s in a filling to be a maximal consecutive run of entries in a column which are all filled with the letter $n$.
    
    \begin{lemma}\label{duplicates}
      Given a filling of a one-column shape $\mu=(1^r)$ and largest entry $n$, there is a unique way of ordering the $n$'s in the filling, say $n_1,\ldots,n_{\alpha_n}$, such that the following two conditions are satisfied.
      \begin{enumerate}
         \item Any consecutive block of $n$'s in the column appears in the sequence in order from bottom to top, and
         \item If we remove $n_1,\ldots,n_{\alpha_n}$ in that order, and let $d_i$ be the amount that the major index of the column decreases at the $i$th step, then the sequence $d_1,d_2,\ldots,d_{\alpha_n}$ is weakly increasing.
      \end{enumerate}
    \end{lemma}
    
    We now can define the map $\majcode$ on words, that is, for one-column fillings.
    
    \begin{definition}
      Let $\sigma$ be any filling of a column shape $\mu=(1^r)$.  We define $\majcode(\sigma)=\majcode(\Standardize(\sigma))$, where $\majcode$ of a standard filling is defined to be the $\majcode$ of its reading word (which is a permutation). 
    \end{definition}
    
     \begin{example}
       Let $\sigma$ be the one-column filling whose reading word is $6434666251664$, the standardization labeling on the $6$'s is shown by the subscripts:
          $$6_2\,4\, 3\, 4\, 6_3\, 6_4\, 6_5\, 2\, 5\, 1\, 6_1\, 6_6\, 4$$
          
       Since this one-column shape has size $13$, the filling $\Standardize(\sigma)$ will have the $6$'s relabeled as the numbers from $8$ to $13$ according to the subscripts above: $$\mathbf{9}\,4\,3\,4\,\mathbf{10}\,\mathbf{11}\,\mathbf{12}\,2\,5\,1\,\mathbf{8}\,\mathbf{13}\,4$$
          
       We then remove the $13,12,\ldots,8$ in order.  This results in a sequence of difference values $1,3,3,3,5,7$, which is weakly increasing.
          
       We are left with a column with reading word $4342514$, in which there is only one $5$, so $\Standardize$ changes that to a $7$.  We remove this to obtain a difference of $1$ in the major index.  We are left with $434214$, in which the $4$'s are standardized as follows: $$4_1 3 4_2 2 1 4_3\to 435216.$$  Removing these in order from $6$ down to $1$ decreases the major index by $0,2,3,2,1,0$, respectively.  Therefore, $$\majcode(\sigma)=1,3,3,3,5,7,1,0,2,3,2,1,0.$$  Note that this sequence is $\{6,6,6,6,6,6,5,4,4,4,3,2,1\}$-weakly increasing.
     \end{example}

    \begin{proposition}
      The map $\majcode$ is a weighted set isomorphism $\fancy{F}_{(1^n)}^{r(\alpha)}\to C_{(1^n),A}$ for any alphabet $A$ with content $\alpha$, and any one-column partition shape $(1^n)$.
    \end{proposition}
    
    \begin{proof}
      Carlitz's work shows that $\majcode$ is an isomorphism in the case that $\alpha=(1,1,\ldots,1)$, i.e. $A$ has one of each letter from $1$ to $n$.  In the case of repeated entries, we note that $\majcode$ is still injective.  Indeed, given a code corresponding to a filling, there is a unique place to insert the next number at each step - by applying the $\Standardize$ map, using Carlitz's bijection, and then un-standardizing in the unique way so that the order of entries is preserved and the resulting alphabet is $A$.
      
      Now, notice that by our definition of $\majcode$ and Lemma \ref{duplicates}, the codes we get are $A$-weakly increasing.  We claim that they are also Carlitz codes:  at the $i$th step, there are $n-i+1$ letters remaining, and the difference $d_i$ is either the position of the letter we're removing plus the number of descents strictly below it, or the number of descents weakly below it.  Therefore, the maximum value of $d_i$ is $n-i+1$, and so $d_1d_2\cdots d_n$ is a Carlitz code and is $A$-weakly increasing.  It follows that $\majcode$ is an injective morphism of weighted sets $\fancy{F}_{(1^n)}^{r(\alpha)}|_{\inv=0}\to C_{(1^n),A}$.
      
      Finally, notice that the two sets have the same cardinality: each has cardinality $\binom{n}{\alpha}$ where $\alpha$ is the content of the alphabet $A$.  It follows that $\majcode$ is bijective, as desired.
    \end{proof}

 \section{Specialization at $t=1$}\label{t=1}
 
 In this section, we give a combinatorial proof of the specialization of Conjecture \ref{InvMaj} at $t=1$, namely $\widetilde{H}_\mu(x;q,1)=\widetilde{H}_{\mu^\ast}(x;1,q)$. 
 
  By the combinatorial formula in \cite{HHL}, it suffices to prove that, for any content $\alpha$,  
 \begin{equation}\label{onlyq}
 \mathop{\sum_{\sigma:\mu\to \ZZ_+}}_{|\sigma|=\alpha} q^{\maj(\sigma)}=\mathop{\sum_{\rho:\mu^\ast\to \ZZ_+}}_{|\rho|=\alpha}q^{\inv(\rho)}.
 \end{equation}
 
  To prove this, we build on the Carlitz bijection, defined in Section \ref{Carlitz}.  Let $$f=\invcode^{-1}\circ \majcode$$ be the Carlitz bijection on permutations of a given ordered alphabet with $n$ distinct entries.  We first prove Equation \ref{onlyq} in the case that $\alpha=(1,1,\ldots,1)$.
 
  \begin{definition}\label{cyclic}
    We say that a sequence of numbers $a_1,\ldots, a_n$ are in \textit{cyclic order} if there exists an index $i\in [n]$ for which $$a_{i+1}\le a_{i+2}\le \cdots \le a_{n}\le a_{1}\le a_2 \le \cdots \le a_i.$$
  \end{definition}
  
 \begin{proposition}\label{distinct}
   For any fixed partition $\lambda$, we have $\sum_{\sigma} q^{\maj(\sigma)}=\sum_{\rho} q^{\inv(\rho)}$ where the first sum ranges over all fillings $\sigma:\lambda\to \ZZ_+$ of $\lambda$ with distinct entries, and the second ranges over all fillings $\rho:\lambda^{\ast}\to \ZZ_+$ of the conjugate partition $\lambda^\ast$ with distinct entries.
 \end{proposition}
 
 \begin{proof}
   We extend the bijection $f$ as follows.
   
   Given a filling $\sigma$ of $\lambda$, let $v^{(1)},v^{(2)},\ldots,v^{(k)}$ be the words formed by reading each of the columns of $\lambda$ from top to bottom.  Let $w^{(i)}=f(v^{(i)})$ for each $i$, so that $\maj(v^{(i)})=\inv(w^{(i)})$.  Notice that $\maj(\lambda)=\sum_{i=1}^{k} \maj(v^{(i)})$.  We aim to construct a filling $\rho$ of $\lambda^{\ast}$ such that $\inv(\rho)=\sum_{i=1}^k \inv(w^{(i)})$.
   
   Let the bottom row of $\rho$ be $w^{(1)}$.  To construct the second row, let $t_1=w^{(1)}_1$ be the corner letter.  Let $x_1,x_2,\ldots,x_{r}$ be the unique ordering of the letters of $w^{(2)}$ for which the sequence $t_1,x_1,x_2,\ldots,x_{r}$ is in cyclic order.  Notice that if $x_i$ is placed in the square above $t$, it would be part of exactly $i$ relative inversions to the right of it, since $x_1,\ldots,x_{i-1}$ would form inversions with it and the others would not.
   
   Now, in $w^{(2)}$, let $i_k$ be the number of inversions whose left element is the $k$th letter of $w^{(2)}$.  Then write $x_{i_1}$ in the square above $t_1$ in order to preserve the number of inversions the first letter is a part of.  Then for the square above $t_2=w^{(1)}_2$, similarly order the remaining $x$'s besides $x_{i_1}$ in cyclic order after $t_2$, and write down in this square the unique such $x_{i_2}$ for which it is the left element of exactly $i_2$ inversions in its row.  Continue this process for each $k\le r$ to form the second row of the tableau.
   
   Continue this process on each subsequent row, using the words $w^{(3)}, w^{(4)},\ldots$, to form a tableau $\rho$.  We define $f(\sigma)=\rho$, and it is easy to see that this construction process is reversible (strip off the top row and rearrange according to inversion numbers, then strip off the second, and so on.)  Thus we have extended the Carlitz bijection to tableaux of content $\alpha=(1,1,\ldots,1)$, proving the result in this case.
 \end{proof}
 
 \begin{remark}
   This proof did not depend on any aspect of the bijection other than the fact that it preserves the statistics.  Thus $f$ can be replaced by, say, the Foata bijection \cite{Foata} and the entire proof is still valid.
 \end{remark}
 
 Using this proposition, we prove two technical lemmata about the $q$-series involved.  Define $\inv_w(R)$ to be the number of relative inversions in a row $R$ given a filling $w$ of the row directly beneath it.
 
 \begin{lemma}\label{technical1}
   Let $R$ be the $i+1$st row in a partition diagram $\lambda$ for some $i\ge 1$. Let $w=w_1,\ldots,w_{\lambda_{i}}$ be a fixed filling of the $i$th row, underneath $R$.  Let $a_1,\ldots,a_{\lambda_{i+1}}$ be any $\lambda_{i+1}$ distinct positive integers.  Then $$\sum q^{\inv_w(R)}=(\lambda_{i})_q!$$ where the sum ranges over all fillings of the row $R$ with the integers $a_1,\ldots,a_{\lambda_{i+1}}$ in some order.
 \end{lemma}
 
 \begin{proof}
   We know that $$\sum_{r\in S_{\lambda_{i+1}}\cdot(a)} q^{\inv(r)}=(n)_q!.$$  We use a similar process to that in Proposition \ref{distinct} to construct a bijection $\phi$ from the set of permutations $r$ of $a_1,\ldots,a_{\lambda_{i+1}}$ to itself such that $\inv_w(\phi(r))=\inv(r)$.
   
   Namely, let $r=r_1,\ldots,r_{\lambda_{i+1}}$ be a permutation of $a_{1},\ldots,a_{\lambda_{i+1}}$ and let $i_k$ be the number of inversions that $r_k$ is a part of in $r$ for each $k$.  Let $x_0,\ldots,x_{\lambda_{i+1}}$ be the ordering of the letters of $r$ for which $w_1, x_0,\ldots,x_{\lambda_{i+1}}$ is in cyclic order.   Let the first letter of $\phi(r)$ be $x_{i_1}$, remove $x_{i_1}$ from the sequence, and repeat the process to form the entire row from the letters of $r$.   Let $\phi(r)$ be this row. 
   
   The map $\phi$ can be reversed by using the the all-$0$'s word for $w$ and using the same process as above to recover $r$ from $\phi(r)$.  Thus $\phi$ is bijective.  Moreover $\inv_w(\phi(r))=\inv(r)$ by construction.  This completes the proof.
 \end{proof}
 
 \begin{lemma}\label{technical2}
   Let $r$ be the $(i+1)$st row in a partition diagram $\lambda$ for some $i\ge 1$. Let $w=w_1,\ldots,w_{\lambda_{i}}$ be a fixed filling of the row directly underneath $r$.  Let $a_1,\ldots,a_{\lambda_{i+1}}$ be positive integers, with multiplicities $m_1,\ldots,m_k$.  Then $$\sum q^{\inv_w(r)}=\binom{\lambda_{i+1}}{m_1,\ldots,m_k}_q=\frac{(\lambda_{i+1})_q!}{(m_1)_q!\cdots (m_k)_q!}$$ where the sum ranges over all distinct fillings of the row $r$ with the integers $a_1,\ldots,a_{\lambda_{i+1}}$ in some order.
 \end{lemma}
 
 \begin{proof}
   Multiplying both sides of the relation by $(m_1)_q!\cdots(m_k)_q!$, we wish to show that $$(m_1)_q!\cdots(m_k)_q!\sum q^{\inv_w(r)}=(\lambda_{i+1})_q!.$$
   This follows immediately by interpreting $(\lambda_{i+1})_q!$ and each $(m_i)_q!$ as in Lemma \ref{technical1}, and assigning all possible orderings to the repeated elements and counting the total number of relative inversions in each case.
 \end{proof}
 
 We are now ready to prove Equation \ref{onlyq}.
 
 \begin{theorem}
   We have $$\mathop{\sum_{\sigma:\mu\to \ZZ_+}}_{|\sigma|=\alpha} q^{\maj(\sigma)}=\mathop{\sum_{\rho:\mu^\ast\to \ZZ_+}}_{|\rho|=\alpha}q^{\inv(\rho)}.$$
 \end{theorem}
 
 \begin{proof}
   We break down each sum according to the contents of the columns of $\mu$ and the rows of $\mu^\ast$, respectively.  For a given multiset of contents of the columns, where the entries in the $i$th column have multiplicities $m^{(i)}_1,\ldots,m^{(i)}_{k_i}$, we have that $$\sum_{\sigma} q^{\maj(\sigma)}=\prod_i\binom{\mu_i'}{m^{(i)}_1,\ldots,m^{(i)}_{k_i}}_q,$$ where the sum ranges over all fillings $\sigma$ with the given column entries.  By Lemma \ref{technical2}, we have that the corresponding sum over fillings $\rho$ with the given contents in the rows of $\mu^\ast$ is the same: $$\sum_{\rho}q^{\inv(\rho)}=\prod_i\binom{\mu_i'}{m^{(i)}_1,\ldots,m^{(i)}_{k_i}}_q.$$
   
   Summing over all possible choices of the entries from $\alpha$ for each column of $\mu$, the result follows.
 \end{proof}

 \section{Hook Shapes} \label{HookShapes}
     
     We now demonstrate a bijective proof of Conjecture \ref{InvMaj} in the case that $\mu$ is a \textit{hook shape}, that is, $\mu=(m,1,1,1,\ldots,1)$ for some $m$.  There is a known combinatorial formula for the $q,t$-Kostka poloynomials in the case of hook shapes $\mu$ given by Stembridge \cite{Stembridge}, but it does not involve the $\inv$ and $\maj$ statistics.  
     
     The symmetry of $\inv$ and $\maj$ was demonstrated for fillings of hook shapes having \textit{distinct} entries in \cite{DMM}, and makes use of the Foata bijection.  In this section, we instead use the Carlitz bijection to prove the result, which will hold for arbitrary fillings by the results in Section \ref{CarlitzWords}.
     
     \begin{lemma}\label{leftzero} We have the following two facts about one-column and one-row shapes respectively.
       \begin{itemize}
       \item Given a filling $\sigma$ of a one-column shape, suppose $A=a_1\ge \cdots \ge a_n$ is the alphabet of its entries written in the standardization order as in Proposition \ref{OneColumnStandardize}, from greatest to least.  Then if $a_i$ is the bottommost entry in $\sigma$, then the first $0$ in $\majcode(\sigma)$ is in position $i$ from the left.
       \item Given a filling $\rho$ of a one-row shape, suppose $A=a_1\le \cdots \le a_n$ is the alphabet of its entries written in order with ties broken in reading order.  Then if $a_i$ is the leftmost entry in $\sigma$, then the first $0$ in $\invcode(\sigma)$ is in position $i$ from the left.
       \end{itemize}
     \end{lemma}
     
     \begin{proof}
       For the filling $\sigma$ of a one-column shape, recall that we define $\majcode$ by removing the entries one at a time from greatest to least in standardization order.  The only time the difference in major index is $0$ is when the entry is on the bottom, and so the first time this occurs is when we remove the bottommost entry $a_i$ from the filling (i.e. at the $i$th step).
       
       For the filling $\rho$ of a one-row shape, note that the leftmost entry $a_i$ always has an inversion code number of $0$.  Moreover, if any entry $b$ to its right also has an inversion code number of $0$, then $b\ge a_i$ for otherwise it would be the smaller entry of an inversion (with $a_i$ itself).  It follows that $a_i$ is the smallest entry whose inversion code number is $0$.
     \end{proof}
     
     We now define a map from fillings of hook shapes to \textit{pairs} of partial codes that we call \textit{hook codes}.  
     \begin{definition}
     Let $\sigma$ be a filling of a hook shape $\mu$.  We define the \textit{hook codes} of $\sigma$ to be the pair of codes consisting of the $\invcode$ of its bottom row and the $\majcode$ of its leftmost column, along with the data of which entries occur in the row and which occur in the column.
     \end{definition} 
     
     Notice that, by the standardization orderings on the row and column of $\mu$ as defined in Section \ref{CarlitzWords}, if the corner square in $\mu$ is one of the repeated letters $a$ of the filling, then it is considered the largest $a$ in its column and the smallest $a$ in its row.  
     
     Thus we can define a standardization ordering on fillings of hook shapes: we order the letters is smallest to largest, with the following tie-breaking rules.  
     \begin{itemize}
     \item If two copies of the letter $a$ appear in the left column, the tie is broken as in Section \ref{CarlitzWords}.  
     \item If they appear in the bottom row, then the leftmost $a$ comes first.  
     \item If one appears in the column and the other in the row, the one in the column comes first.
     \end{itemize}
     
    This enables us to represent hook codes visually, as shown in the following example.
     
     \begin{example}\label{hookexample1}
     Consider the filling $\sigma$ of a hook shape shown below.  The $2$ in the corner is considered to be greater than the $2$ above it and less than the $2$ to its right.  To represent the hook code of $\sigma$, we write the entries of the filling in the standardization ordering, and write the $\invcode$ and (the reverse of) $\majcode$ of the bottom row and left column respectively underneath the corresponding letters.  
     
     \begin{center}
       \includegraphics{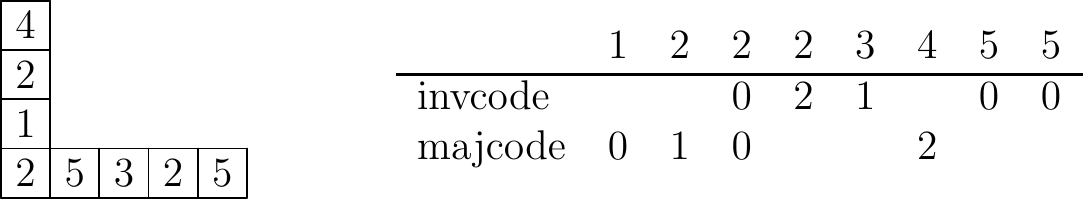}
     \end{center}
     
     Notice that the $\majcode$ is written \textit{backwards}, because the entries are in increasing order.
     
     \end{example}
     
     We now characterize the pairs of codes that correspond to fillings of hook shapes.
     
     \begin{lemma}\label{HookCodes}
       Let $\mu$ be a hook shape of height $h$ and width $l$ with $h+l=n$, and let $A=\{a_1\le \cdots \le a_n\}$ be an ordered multiset.  A pair of partial codes $(X,Y)$ of lengths $l$ and $h$ respectively is a hook code of some filling $\sigma$ of $\mu$ if and only if the four conditions below are satisfied.
       \begin{enumerate}
       \item The leftmost $0$ of $X$ matches the rightmost $0$ of $Y$.
       \item The two codes do not overlap in any other position, and every position is part of at least one of the two codes.
       \item The code $X$ is an element of $C_l$ and is $A$-weakly increasing, where we restrict $A$ to the $l$ letters corresponding to the positions of the entries of $X$.
       \item The code $Y$, \textit{when read backwards}, is an element of $C_h$ and is $A$-weakly increasing, where we restrict $A$ to the $h$ letters corresponding to the positions of the entries of $Y$.
       \end{enumerate}
     
     \end{lemma}
     
     \begin{proof}
     First we show that the hook code of any filling $\sigma$ of $\mu$ satisfies the four conditions.  Condition 1 follows immediately from Lemma \ref{leftzero}, because the major index code is written in reverse order.  Condition 2 is clear since every entry is in either the row or the column and only the corner square is in both.  Conditions 3 and 4 follow immediately from the definition of hook codes.
     
     Now, suppose we have a pair of codes satisfying conditions 1-4.  Then there is a unique way to form a row and a column of entries based on their elements, since they are both valid Carlitz codes and are $A$-weakly increasing by conditions 3 and 4.  Because of condition 1 and Lemma \ref{leftzero}, the leftmost entry of the row is the same as the bottommost entry of the column, and so we can put them together to form a filling $\sigma$ of a hook shape.  Because of condition 2, the hook shape $\mu$ has the appropriate size and shape, and we are done.
     \end{proof}

     Using Lemma \ref{HookCodes}, we can now define our bijection.
     
     \begin{definition}
       For any hook shape $\mu$ and content $\alpha$, let $\phi:\fancy{F_\mu^\alpha}\to\fancy{F}_{\mu^\ast}^{r(\alpha)}$ be the map defined by interchanging the pair of hook codes of a given filling and writing them backwards, and also reversing its alphabet.
     \end{definition}
     
     \begin{example}
       Starting with the tableau in Example \ref{hookexample1}, if we reverse the alphabet, interchange $\invcode$ and $\majcode$, and write the codes in backwards order, then we obtain the filling and pair of codes below.  It follows that the filling in Example \ref{hookexample1} maps to the filling below under $\phi$.
       \begin{center}
         \includegraphics{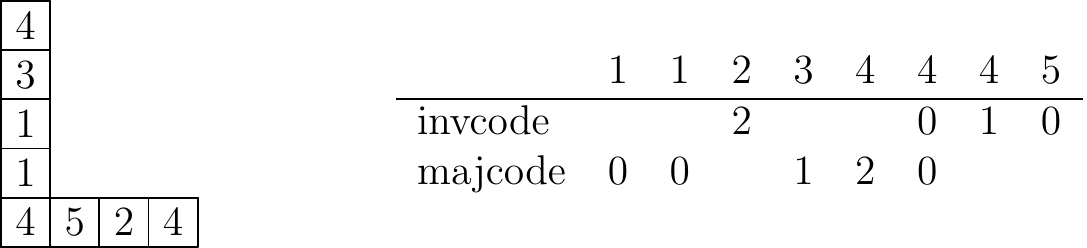}
       \end{center}
     \end{example}
     
     \begin{theorem}
       We have that $$\maj(\phi(\sigma))=\inv(\sigma)$$ and $$\inv(\phi(\sigma))=\maj(\sigma)$$ for any filling $\sigma$ of a given hook shape $\mu$.  Moreover, $\phi$ is a bijection from $\fancy{F}_\mu^\alpha$ to $\fancy{F}_{\mu^\ast}^{r(\alpha)}$ for any content $\alpha$.
     \end{theorem}
     
     \begin{proof}
     
       Clearly $\phi$ interchanges $\inv$ and $\maj$, since it interchanges the $\invcode$ and $\majcode$ of the filling.  To show it is a well-defined map into fillings of the conjugate shape, note that reversing and interchanging the codes and reversing the alphabet results in a pair of codes that satisfy conditions 1-4 of Lemma \ref{HookCodes}.
       
       Finally, $\phi$ is a bijection - in fact, it is an involution - because the operations of reversing the alphabet, interchanging the pair of codes, and writing the codes in the reverse order are all involutions.  
     \end{proof}
     
     \begin{corollary}
       The map $\phi$ above satisfies the conditions of the Conjecture \ref{InvMaj}, proving combinatorially that $$\widetilde{H}_\mu(x;q,t)=\widetilde{H}_{\mu^\ast}(x;t,q)$$ when $\mu$ is a hook shape.
     \end{corollary}

\section{Hall-Littlewood Specialization at $q=0$} \label{q=0}

We now turn to the specialization in which one of the statistics is zero.  In particular, setting $q=0$, the symmetry relation becomes 
$$\widetilde{H}_{\mu}(x;0,t)=\widetilde{H}_{\mu^\ast}(x;t,0),$$ which is a symmetry relation between the transformed \textit{Hall-Littlewood polynomials} $\widetilde{H}_\mu(x;t):=\widetilde{H}_\mu(x;0,t)$.  In this case the symmetry relation becomes
\begin{equation}\label{HallLittlewoodSymmetry}
\mathop{\sum_{\sigma:\mu\to \ZZ_+}}_{\inv(\sigma)=0} t^{\maj(\sigma)}x^{\sigma}=\mathop{\sum_{\rho:\mu^{\ast}\to \ZZ_+}}_{\maj(\rho)=0} t^{\inv(\rho)}x^{\rho}.
\end{equation} Combinatorially, we are trying to find natural morphisms $$\varphi:\fancy{F}_\mu^\alpha|_{\inv=0}\,\to\, \fancy{F}_{\mu^\ast}^{r(\alpha)}|_{\maj=0}$$ of weighted sets, where $\fancy{F}_\mu^\alpha|_{\inv=0}$ is equipped with the $\maj$ statistic, and $\fancy{F}_{\mu^\ast}^\alpha|_{\maj=0}$ is equipped with the $\inv$ statistic.
For the bijection $r(\alpha)$, we will use the \textit{reverse} map of Definition \ref{reverse}.

\subsection{Generalized Carlitz Codes}
  
  In the context of Hall-Littlewood symmetry, we can think of the Carlitz bijection as a solution to the case in which $\mu=(1^n)$ is a straight shape with one column, filled with distinct entries.  Thus, we wish to generalize the notion of a Carlitz code to fillings of arbitrary shapes having $\inv$ or $\maj$ equal to $0$, using arbitrary alphabets.
  
  Our generalization is motivated by the monomial basis of the Garsia-Procesi modules in \cite{GarsiaProcesi}, which are closely connected to the cocharge ($\maj$) statistic.  We define a generalized Carlitz code as follows.
  
  \begin{definition}
  A word having letters in $\{0,1,2,\ldots\}$ is \textit{Yamanouchi} if every suffix contains at least as many $i$'s as $i+1$'s for all $i\ge 0$.

    A word $w$ has \textit{content} $\alpha=(\alpha_1,\ldots,\alpha_k)$ if exactly $\alpha_i$ of the entries of $w$ are equal to $i-1$ for each $i$.  We also sometimes say it has content $A$ where $A$ is the multiset of letters of $w$.

    Finally, a word $w=w_1\cdots w_n$ is $\mu$-\textit{sub-Yamanouchi}, or \textit{$\mu$-Carlitz}, if there exists a Yamanouchi word $v=v_1\cdots v_n$ of content $\mu$ such that $w_i<v_i$ for all $i$.
  \end{definition}
  
  \begin{example}
    The sub-Yamanouchi words for shape $\mu=(1,1,1,\ldots,1)$ are precisely the classical Carlitz codes.
  \end{example}
  
  We will see that the $\mu$-sub-Yamanouchi words are the correct analog of Carlitz codes in the case that our Young diagram fillings have distinct entries.  However, in general we require the following more precise definition.
  
  \begin{definition}
    We define $C_{\mu,A}$ to be the collection of all $\mu$-sub-Yamanouchi codes which are $A$-weakly increasing (see Definition \ref{Aweaklyincreasing}).  We call such codes \textit{generalized Carlitz codes}, and we equip this collection with the statistic $\Sigma:C_{\mu,A}\to \ZZ$ by $\Sigma(c)=\sum c_i$, forming a weighted set $(C_{\mu,A};\Sigma)$.
  \end{definition}
  
  We now introduce the concept of the \textit{monomial} of a code.  The next three definitions are compatible with the notation in \cite{GarsiaProcesi}.
    
    \begin{definition}
    Fix variables $x_1,x_2,\ldots$.   For any finite code $c$ of length $n$, define its \textit{monomial} to be $$x^c=x_n^{c_1}x_{n-1}^{c_2}\cdots x_1^{c_n}.$$  Also let $\fancy{C}_A(\mu)$ be the set of all monomials $x^c$ of $\mu$-sub-Yamanouchi words $c$ that are $A$-weakly increasing.
    \end{definition}
    
    In \cite{GarsiaProcesi}, the authors define a similar set of monomials $\fancy{B}(\mu)$, which are the generators of the modules $R_\mu$ that arise naturally in the study of the Hall-Littlewood polynomials.  We will see that in the case $A=\{1,2,\ldots,n\}$, we have $\fancy{C}_A(\mu)=\fancy{B}(\mu)$, by showing that the sets $\fancy{C}_A(\mu)$ satisfy a generalized version of the recursion in \cite{GarsiaProcesi}.  To state this recursion we require two more definitions, which follow the notation in \cite{GarsiaProcesi}.
   
    \begin{definition}
      Given a partition $\mu$, define $\mu^{(i)}$ to be the partition formed by removing the corner square from the column $a_i$ containing the last square in the $i$th row $\mu_i$.
    \end{definition}
   
    \begin{definition}
    Given a set of monomials $\fancy{C}$ and a monomial $m$, we write $m\cdot \fancy{C}$ to denote the set of all monomials of the form $m\cdot x$ where $x\in \fancy{C}$.
    \end{definition}
    
    The following recursion defines the sets $\fancy{B}(\mu)$.  
    
    \begin{definition}
    The sets $\fancy{B}(\mu)$ are defined by $\fancy{B}((1))=\{1\}$ and the recursion
     $$\fancy{B}(\mu)=\bigsqcup_{i-1}^{\mu_1^\ast}x_n^{i-1}\cdot \fancy{B}(\mu^{(i)}).$$  We refer to these sets as the \textit{Garsia-Procesi module bases}.
    \end{definition}

    We require one new definition in order to state our general recursion in the next proposition.
      
    \begin{definition}
      Let $A=\{a_1,a_2,\ldots,a_n\}$ with $a_1\le a_2\le \cdots \le a_n$ be a multiset of positive integers, and let $\lambda$ be a partition of $n-1$.  We define $\fancy{C}^{(t)}_A(\lambda)$ to be the set of all monomials $x^d$ of $\lambda$-sub-Yamanouchi words $d_1,\ldots,d_{n-1}$ that are $A\setminus \{a_1\}$-weakly increasing and if $a_1=a_2$ then $d_1\ge t$.
    \end{definition}
    
    \begin{proposition}[General Recursion]\label{recursion}
      For any partition $\mu$ of $n$ and any multiset of positive integers $A=\{a_1,a_2,\ldots,a_n\}$ with $a_1\le a_2\le \cdots \le a_n$, we have $$\fancy{C}_{A}(\mu)=\bigsqcup_{i=1}^{\mu_1^\ast}x_n^{i-1}\cdot \fancy{C}^{(i-1)}_{A}(\mu^{(i)}).$$
    \end{proposition}
    
    We defer the proof of this recursion to section \ref{Proofs}.
    
    Notice that in the case $A=\{1,2,\ldots,n\}$, since there are no repeated entries, Proposition \ref{recursion} reduces to  $$\fancy{C}(\mu)=\bigsqcup_{i-1}^{\mu_1^\ast}x_n^{i-1}\cdot \fancy{C}(\mu^{(i)}).$$  Since this is the same as the recursion given for the sets $\fancy{B}(\mu)$ described in the previous section, and $\fancy{C}_{\{1\}}((1))=\{x_1\}=\fancy{B}((1))$, we have the following corollary.
      
    \begin{corollary}
      If $A=\{1,2,\ldots,n\}$, we have  $\fancy{C}_A(\mu)=\fancy{B}(\mu)$.
    \end{corollary}
        
      As noted in \cite{GarsiaProcesi}, we can now also enumerate the sets $\fancy{C}_A(\mu)$ in the case $A=\{1,2,\ldots,n\}$.  For, in this case the simplified recursion gives $$|\fancy{C}_A(\mu)|=\sum_i|\fancy{C}_A(\mu^{(i)})|$$ with $|\fancy{C}_{\{1\}}((1))|=1$.  But the multinomial coefficients $\binom{n}{\mu}$ satisfy $\binom{1}{1}=1$ and the same recursion: $$\binom{n}{\mu}=\sum_i \binom{n}{\mu^{i}}.$$
        
     \begin{corollary}\label{standardinv}
       If $A=\{1,2,\ldots,n\}$, we have $$|\fancy{C}_A(\mu)|=\binom{n}{\mu}.$$
     \end{corollary}

\subsection{Inversion Codes}\label{CochargeWords}

  We can now generalize the inversion code of a permutation to arbitrary fillings $\rho$ with $\maj(\rho)=0$.
  
   \begin{definition}
     Let $\rho$ be a filling of $\mu^\ast$ having $\maj(\rho)=0$.  Order its entries by size with ties broken in reading order to form a totally ordered alphabet $A=\{a_1,\ldots,a_n\}$. Then its \textit{inversion code}, denoted $\invcode(\rho)$, is the sequence $c_1\cdots c_n$ whose $i$th entry $c_i$ is the number of attacking pairs having $a_i$ as its smaller entry.
   \end{definition}
  
    \begin{example}
     Consider the following tableau.
     \begin{center}
     \includegraphics{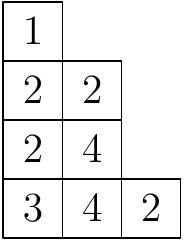}
     \end{center}
     
     There are three attacking pairs in this diagram: the $2$ in the bottom row is attacked by the $3$ and $4$ in its row, and the $3$ is attacked by the $4$ in the second row.  When we order the entries in reading order and record the number of larger numbers that attack it, we get the following table.
     \begin{center}
     $\begin{array}{l|cccccccc}
       \text{Entries} & 1 & 2 & 2 & 2 & 2 & 3 & 4 & 4\\\hline
       \text{Code}    & 0 & 0 & 0 & 0 & 2 & 1 & 0 & 0
     \end{array}$
     \end{center}
     
     Therefore, the inversion code of the filling above is $00002100$.
     
    \end{example}
  
  \begin{theorem}\label{Invcode}  
  The inversion code of any filling $\rho\in \fancy{F}_{\mu^\ast}^{\alpha}$ is $\alpha$-weakly increasing and $\mu$-sub-Yamanouchi.  Moreover, the map $$\invcode:\fancy{F}_{\mu^\ast}^\alpha\mid_{\maj=0}\to C_{\mu,A}$$ is an isomorphism of weighted sets.
  \end{theorem}
  
  The proof of Theorem \ref{Invcode} is somewhat technical, and so we defer it to Section \ref{Proofs}.

   \subsection{Major Index Codes}
   
   To complete the proof of the Hall-Littlewood case, it now suffices to find a weighted set isomorphism $$\majcode:(\fancy{F}_{\mu}^\alpha|_{\inv=0};\maj)\to C_{\mu,A}$$ where $\alpha$ is the content of the alphabet $A$.
   
   Recall the recursion for the $\mu$-sub-Yamanouchi codes of content $A$ from Proposition \ref{recursion}:
   $$\fancy{C}_{A}(\mu)=\bigsqcup_{i=1}^{\mu_1^\ast}x_n^{i-1}\cdot \fancy{C}^{(i-1)}_{A}(\mu^{(i)}).$$
   Using this recursion, one possible strategy for constructing $\majcode$ is by showing combinatorially that $\fancy{F}_{\mu}^\alpha|_{\inv=0}$ satisfies a similar recursion.  
   
   In this section, we present some partial progress towards finding the map $\majcode$.  All of our work is based on the following four-step approach to the problem.
   
    \begin{enumerate}
    \item[\textbf{Step 1.}] Consider the content $(1^n)$ corresponding to fillings with distinct entries, and find an explicit weighted set isomorphism $$\psi:(\fancy{F}_{\mu}^{(1^n)}|_{\inv=0};\maj)\to \bigsqcup_{d}(\fancy{F}_{\mu^{(d+1)}}^{(1^{n-1})}|_{\inv=0};\maj+d).$$  That is, $\psi$ should send an inversion-free filling $T$ of $\mu$ to an inversion-free filling $\psi(T)$ of $\mu^{(d+1)}$ for some $d$, such that $$\maj(\psi(T))=\maj(T)-d.$$
    \item[\textbf{Step 2.}] Define the $\majcode$ of a filling $T$ having content $(1^n)$ to be $d_1d_2\ldots d_n$ where $$d_k=\maj(\psi^{k}(T))-\maj(\psi^{k-1}(T)).$$ 
    \item[\textbf{Step 3.}] Check the base case of a single square, and conclude that because the recursion is satisfied, $\majcode$ is an isomorphism of weighted sets $$(\fancy{F}_{\mu}^{(1^n)}|_{\inv=0},\maj)\to (C_{\mu,[n]},\Sigma),$$ where $C_{\mu,[n]}$ are the generalized Carlitz codes of shape $\mu$ and content $[n]$.
    \item[\textbf{Step 4.}] Show that there is a standardization map $$\Standardize:\fancy{F}_{\mu}^\alpha|_{\inv=0}\to \fancy{F}_{\mu}^{1^n}|_{\inv=0}$$ that respects $\maj$, such that the composition $\majcode\circ\Standardize$ is a bijection to $C_{\mu,A}$ where $A$ is the alphabet with content $\alpha$.  That is, show that after standardizing, we get a major index code which is $A$-weakly increasing, and none of these codes are mapped to twice.
    \end{enumerate}
      
   \subsubsection{Killpatrick's Method for Standard Fillings}

   For Step 1 in our strategy, in which $A=\{1,2,\ldots,n\}$ is an alphabet with no repeated letters, such a map can easily be extracted from the work of Killpatrick \cite{Killpatrick}.  In this paper, the author gives a combinatorial proof of a recursion for a generating function involving \textit{charge}, written $\ch$, and defined in terms of cocharge as $\ch(\mu)=n(\mu)-\cc(\mu)$ where $n(\mu)=\sum_{i}(i-1)\cdot \mu_i$.  Killpatrick defines $W_\mu$ to be the set of words of content $\mu$, and lets $r_{i,\mu}=|\{j>i:\mu_j=\mu_i\}|$.  The recursion is stated as:
    $$\sum_{w\in W_\mu} q^{\ch(w)}=\sum_{i}q^{r_{i,\mu}} \sum_{w\in W_{\mu^{(i)}}}q^{\ch(w)}.$$ 
    
    If we substitute $q\to 1/q$ and multiply both sides by $q^{n(\mu)}$, this becomes
    $$\sum_{w\in W_\mu} q^{\cc(w)}=\sum_{i}q^{i-1}\sum_{w\in W_{\mu^{(i)}}} q^{\cc(w)},$$
    which is equivalent to the recursion we stated in step 1 above.  Killpatrick's map $\psi$ allows us to define a map ${\majcode'}$ that statisfies Steps 1-3 above.  We therefore immediately obtain the following result.
    
    \begin{theorem}
      In the case $\alpha=(1^n)$ of fillings with distinct entries, we have that $\varphi={\majcode'}^{-1}\circ \invcode$ is an isomorphism of weighted sets $$\varphi:\fancy{F}_\mu^{(1^n)}|_{\maj=0}\to \fancy{F}_{\mu^\ast}^{(1^n)}|_{\inv=0}.$$
    \end{theorem}
    
    However, Killpatrick's map $\majcode'$ does not satisfy the requirements of Step 4.  To illustrate this, we consider the case in which $\mu=(1^n)$ is a straight column shape.  In this case, Killpatrick's bijection $\majcode'$ is defined by the following process:
    
    \begin{enumerate}
      \item Given a filling $w$ of a straight column shape such as the one with reading word $1432$ in the diagram below, check to see if the bottommost entry is the largest entry.  If not, cyclically increase each entry by $1$ modulo the number of boxes $n$.  Each such cyclic increase, or \textit{cyclage}, can be shown to decrease the major index by exactly one (see Section \ref{Applications} for details in the language of cocharge).  We perform the minimal number of cyclages to ensure that the bottommost letter is $n$, and let $c_1$ be the number of cyclages used.  (In the figure, $c_1=2$.)
           \begin{center}
            \includegraphics{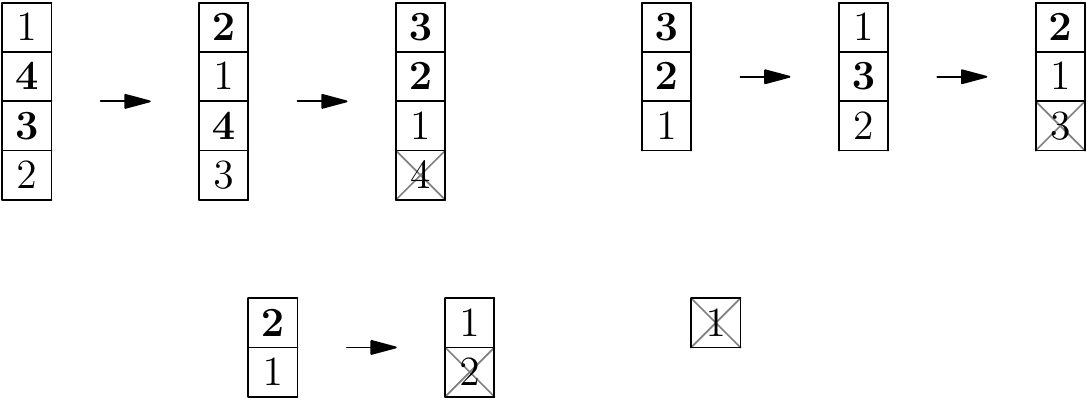}
           \end{center}
           
      \item Once the bottommost entry is $n$, remove the bottom box, and repeat step $1$ on the new tableau.  The resulting number of cyclages used is recorded as $c_2$.  (In the figure, $c_2=2$.)
      \item Continue until there are no boxes left, and set $\majcode'(w)=c_1c_2\cdots c_n$.  (In the figure, $\majcode'(w)=2210$.)
    \end{enumerate}
    
    Now, suppose we had a standardization map $\Standardize$ as in Step $4$.  Consider the one-column tableaux having entries from the alphabet $\{2,2,1,1,1,1\}$ and major index $4$.  There are three such tableaux:
        \begin{center}
            \includegraphics{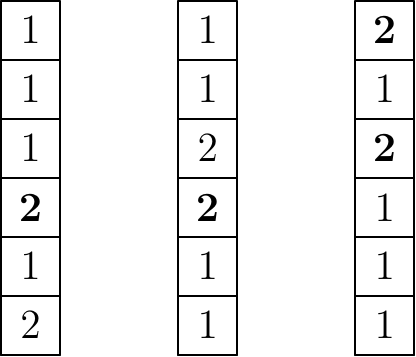}
        \end{center}
                
    There are also three $(1^6)$-sub-Yamanouchi codes that are $\{2,2,1,1,1,1\}$-weakly increasing and sum to $4$, namely:  $$040000$$  $$130000$$  $$220000$$
    
    It follows that $\Standardize$ maps these three tableaux to the three standardized fillings whose codes $\majcode'$ are $040000$, $130000$, and $220000$, respectively.  But these three tableaux are:
      \begin{center}
         \includegraphics{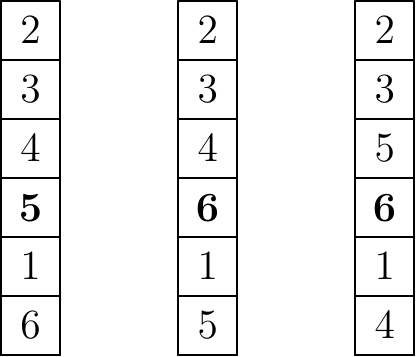}
      \end{center}
    
    Therefore, the map $\Standardize$ cannot preserve the relative ordering of the entries, or even the positions of the descents.  This makes it unlikely that a simple rule for such a standardization map exists.  However, it is possible that there exists a more complicated combinatorial rule for such a map, and we leave this as an open question for future investigation.
    
    \begin{question}
    Is there a natural map $\Standardize$ that satisfies the conditions of Step $4$ for Killpatrick's map $\majcode'$?
    \end{question}
    
    We now return to Carlitz's bijection in the next section, in which we generalize $\majcode$ to arbitrary inversion-free fillings of certain infinite families of shapes.
    
    \subsubsection{Reducing Rectangles to Columns}\label{Rectangles}
    
    The Carlitz bijection on words, defined in section \ref{CarlitzWords}, gives a map $\majcode$ for arbitrary fillings of one-column shapes $\mu$.  We now show present a strategy towards a generalization to all shapes $\mu$, and show that rectangles behave similarly to one-column shapes.
    
    Our primary tool is the following technical result.  This lemma generalizes the fact that if we remove the largest entry $n$ from the bottom of a one-column shape, we get a major index code entry $d=0$.
    
    \begin{proposition}[Main Lemma]\label{ZeroBump}
      Suppose $\sigma:\mu\to \ZZ_+$ is a filling for which $\inv(\sigma)=0$ and the largest entry $n$ appears in the bottom row.  Let $\sigma_\downarrow:\mu^{(1)}\to \ZZ_{+}$ be the filling obtained by:
      \begin{enumerate}
        \item Removing the rightmost $n$ from the bottom row of $\sigma$, which must be in the rightmost column since $\inv(\sigma)=0$,
        \item Shifting each of the remaining entries in the rightmost column down one row,
        \item Rearranging the entries in each row in the unique way so that $\inv(\sigma_\downarrow)=0$.
      \end{enumerate}
     Then the major index does not change: $$\maj(\sigma)=\maj(\sigma_\downarrow).$$
    \end{proposition}
        
    We defer the proof to the Section \ref{Proofs}.
    
    It turns out that the construction $\sigma\to \sigma_\downarrow$ is reversible, and to see this we require the following lemma.
    
    \begin{lemma}\label{PullUp}
      Given two collections of letters $b_1,\ldots,b_{w-1}$ and $a_1,\ldots,a_w$, there is a unique element $a_i$ among $a_1,\ldots, a_w$ such that, in any two-row tableau with $a_1,\ldots,\hat{a_i},\ldots,a_w$ as the entries in the bottom row and $b_1,\ldots,b_{w-1},a_i$ as the entries in the top, with no inversions in the top row, the entry $a_i$ occurs in the rightmost position in the top row.
    \end{lemma}
    
    This lemma allows us to recover $\sigma$ from a tableau $\sigma_\downarrow$ whose second-longest row $\mu_k$ is one square shorter than its longest rows ($\mu_1$ through $\mu_{k-1}$).  We simply raise the appropriate entry $a_i$ from row $\mu_{k-1}$ to row $\mu_k$, then do the same from row $\mu_{k-2}$ to $\mu_{k-1}$, and so on, and finally insert a number $n$ in the bottom row, where $n$ is larger than all of the other entries in $\sigma_\downarrow$.
    
    \begin{example}
     Applying the process in the tableau below, the major indexes of the starting tableau and the ending tableau are both $10$.
      \begin{center}
      \includegraphics{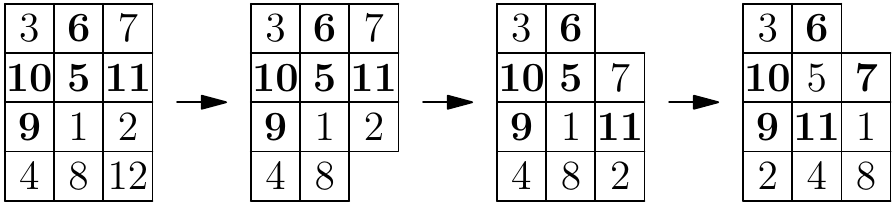}
      \end{center}
    \end{example}
        
    Using Proposition \ref{ZeroBump}, we can provide a new combinatorial proof of the recurrence of Garsia and Procesi for all rectangular shapes $\mu=(a,a,a,\ldots,a)$.  This also will provide the first letter of $\majcode$ for rectangular shapes.
    
    \begin{theorem}\label{RectangleBump}
        Let $A=\{1,2,\ldots,n\}$ be the alphabet with content $\alpha=(1^n)$, and let $\mu=(a,a,a,\ldots,a)$ be a rectangle shape of size $n$.  Then there is a weighted set isomorphism
          $$\psi:(\fancy{F}_{\mu}^{(1^n)}|_{\inv=0};\maj)\to \bigsqcup_{d=0}^{\mu_1^\ast-1}(\fancy{F}_{\mu^{(d+1)}}^{(1^{n-1})}|_{\inv=0};\maj+d)$$
        defined combinatorially by the following process.  
       \begin{enumerate}
        \item Given a filling $\sigma:\mu\to \ZZ_+$ with distinct entries $1,\ldots,n$ and $\inv(\sigma)=0$, let $i$ be the row containing the entry $n$.  Split the filling just beneath row $i$ to get two fillings $\sigma_{top}$ and $\sigma_{bot}$ where $\sigma_{bot}$ consists of rows $1,\ldots,i-1$ of $\sigma$ and $\sigma_{top}$ consists of rows $i$ and above.
        \item Rearrange the entries of the rows of $\sigma_{top}$ in the unique way that forms a filling $\widetilde{\sigma_{top}}$ for which $\inv(\widetilde{\sigma_{top}})=0$.
        \item Apply the procedure of Proposition \ref{ZeroBump} to $\widetilde{\sigma_{top}}$, that is, removing the $n$ from the bottom row and bumping each entry in the last column down one row.  Let the resulting tableau be called $\tau$.
        \item Place $\tau$ on top of $\sigma_{bot}$ and rearrange all rows to form a tableau $\rho$ having $\inv(\rho)=0$.  Then we define $\psi(\sigma)=\rho$.
        \end{enumerate}
         Moreover, if $\maj(\sigma)-\maj(\psi(\sigma))=d$, then $0\le d<\mu_1^\ast$ and we assign $\psi(\sigma)$ to the $d$th set in the disjoint union.
    \end{theorem}

        \begin{remark}
          Theorem \ref{RectangleBump} gives a new combinatorial proof of the recursion $$\sum_{\sigma:\mu\to \ZZ_+ \atop \inv(\sigma)=0}q^{\maj(\sigma)}=\sum_{d} q^{d-1}\sum_{\rho:\mu^{(d)}\to \ZZ_+ \atop \inv(\rho)=0}q^{\maj(\rho)}$$ of Garsia and Procesi for rectangular shapes $\mu$.
        \end{remark}
        
        The map $\psi$ of Theorem \ref{RectangleBump} is illustrated by the example below.
            
            \begin{center}
            \includegraphics{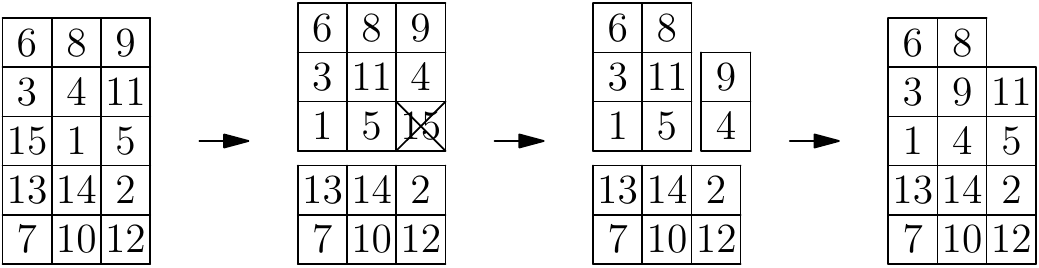}
            \end{center}
            
     The proofs of these results are deferred to sections \ref{Proofs1} and \ref{Proofs2}.  For now, we state some facts pertaining to Theorem \ref{RectangleBump} that will be useful in extending this map to other shapes and alphabets.  Proofs of these facts can also be found in section \ref{Proofs2}.
         
     \begin{proposition}\label{RectBumpCor1} 
     Let $h=\mu_1^\ast$ be the height of the rectangle shape $\mu$.  Since $\mu$ is a rectangle, the shape $\mu^{(d+1)}$ is independent of $d\in \{0,\ldots,h-1\}$, so let $\mu_\downarrow$ be this shape.  Let $\rho\in\fancy{F}_{\mu_\downarrow}^{(1^{n-1})}$, so that there is a copy $\rho_d$ of $\rho$ in $\fancy{F}_{\mu^{(d+1)}}^{(1^{n-1})}|_{\inv=0}$ for all $d=0,\ldots,h-1$.  Let $\sigma_d=\psi^{-1}(\rho_d)$ for each $d$.  Then for each $i=1,\ldots,h$, the largest entry $n$ occurs in the $i$th row in exactly one of $\sigma_0,\ldots,\sigma_{h-1}$.
     \end{proposition}
     
     The next theorem suggests that the standardization map for rectangle shapes can be inherited from the standardization map for single-column shapes described above.
     
     \begin{theorem}[Reducing rectangles to columns]\label{RectBumpCor2}
        For $\sigma\in \fancy{F}_{\mu}^{(1^n)}|_{\inv=0}$ with $\mu$ a rectangle, the value of $d=\maj(\sigma)-\maj(\psi(\sigma))$ can be determined as follows.  Let $\sigma_1$ be the unique element of $\fancy{F}_{\mu}^{(1^n)}|_{\inv=0}$ for which $n$ is in the bottom row and $\psi(\sigma_1)=\psi(\sigma)$, so that ${\sigma_1}_{\downarrow}=\psi(\sigma_1)=\psi(\sigma)$.  Let $a_{h-1},\ldots,a_1,n$ be the entries of the rightmost column of $\sigma_1$ from top to bottom.  Then $d$ is the same as the difference in the major index obtained from inserting $n$ into the $i$th position in the one-column shape with reading word $a_{h-1},\ldots,a_1$.
     \end{theorem}
     
     This theorem is so crucial to the proofs of the results in the next section that it is helpful to give the sequence of $a_i$'s its own name.  We call it the \textit{bumping sequence} of $\sigma$.
     
     \begin{definition}
       Let $\sigma$ be a filling of a rectangle shape $\mu$ having height $h$, with distinct entries $1,2,\ldots,n$.  The \textit{bumping sequence} of $\sigma$ is the collection of entries $a_1,a_2,\ldots,a_{h-1}$ defined as in Theorem \ref{RectBumpCor2} above.  If $n$ is in the $i$th row of $\sigma$, then $a_1,\ldots,a_{i-1}$ are in rows $1$ through $i-1$ respectively, and $a_{i},\ldots,a_{h-1}$ are in rows $i+1$ through $h$.
     \end{definition}
     
     We can also say something about the position of these $a_i$'s given the position of the largest entry.
            
     \begin{proposition}\label{RectBumpCor3}
       Let $\mu$ be a rectangle shape of height $h$, and let $\sigma\in \fancy{F}_{\mu}^{(1^n)}$ with its largest entry $n$ in row $i$.  Then if $a_1,\ldots,a_{h-1}$ is the bumping sequence of $\sigma$, then $a_{i+2},\ldots,a_{h-1}$ all occur in columns weakly to the right of the $n$, and each $a_j$ is weakly to the right of $a_{j-1}$ for $j\ge i+3$.
     \end{proposition}

    \subsubsection{Three Row Shapes}\label{ThreeRowsSection}
       
       We now provide a complete bijection $\majcode$ in the case that $\mu=(\mu_1,\mu_2,\mu_3)$ is a partition with at most three rows.
       
       We start with the definition of $\majcode$ for two-row shapes, which we will use as part of the algorithm for three rows.
      
       \begin{lemma}\label{tworows}
         Let $\mu=(\mu_1,\mu_2)$ be any two-row shape of size $n$.   Then there is a weighted set isomorphism 
           $$\psi:(\fancy{F}_{\mu}^{(1^n)}|_{\inv=0};\maj)\to \bigsqcup_{d=0}^{1}(\fancy{F}_{\mu^{(d+1)}}^{(1^{n-1})}|_{\inv=0};\maj+d)$$
         defined combinatorially by the following process.  Given an element $\sigma$ of $\fancy{F}_{\mu}^{(1^n)}|_{\inv=0}$, that is, a filling of the two-row shape $\mu$ having no inversions, consider its largest entry $n$.
          \begin{enumerate}
            \item If the $n$ is in the bottom row, define $\psi(\sigma)=\sigma_{\downarrow}$ as in Proposition \ref{ZeroBump}.
            \item If the $n$ is in the second row, remove it and re-order the remaining entries in the top row so that there are no inversions.  Let $\psi(\sigma)$ be the resulting filling.
          \end{enumerate}
       \end{lemma}
       
       \begin{proof}
         We first show that $\psi$ is a morphism of weighted sets.  If the $n$ we remove is in the bottom row, then by Proposition \ref{ZeroBump}, the new filling $\sigma_\downarrow=\psi(\sigma)$ is in $\fancy{F}_{\mu^{(1)}}^{(1^{n-1})}|_{\inv=0}$ and has the same major index as $\sigma$.  This means that $\sigma_\downarrow$ is in the $d=0$ component of the disjoint union $$\bigsqcup_{d=0}^{1}(\fancy{F}_{\mu^{(d+1)}}^{(1^{n-1})}|_{\inv=0};\maj+d),$$ and the statistic is preserved in this case.
         
         Otherwise, if the $n$ is in the second (top) row, then $\sigma'=\psi(\sigma)$ is in $\fancy{F}_{\mu^{(2)}}^{(1^{n-1})}|_{\inv=0}$.  We wish to show that the difference in major index, $d=\maj(\sigma)-\maj(\sigma')$, is $1$ in this case.  Indeed, notice that the bottom row remains unchanged after removing the $n$, and so the difference in major index will be the same as if we ignore the extra $\mu_1-\mu_2$ numbers at the end of the bottom row and consider just the rectangle that includes the second row instead.  By Theorem \ref{RectangleBump}, it follows that $d=1$.   Therefore $\psi$ is a morphism of weighted sets.
         
         To show that $\psi$ is bijective, we construct an inverse map $\phi$.  First, let $\sigma'\in \fancy{F}_{\mu^{(1)}}^{(1^{n-1})}|_{\inv=0}$.  Then we can insert $n$ into the bottom row, and if $\mu$ is a rectangle also bump up one of the entries of the bottom row according to Lemma \ref{PullUp}.  This creates a filling $\sigma$ of shape $\mu$ having the same major index as $\sigma'$.  We define $\phi(\sigma')=\sigma$, which defines an inverse map for $\psi$ on the restriction of $\psi$ to $\psi^{-1}\left(\fancy{F}_{\mu^{(1)}}^{(1^{n-1})}|_{\inv=0}\right)$.
         
         Now let $\sigma'$ be a filling of shape $\mu^{(2)}$.  The shape $\mu^{(2)}$ has a longer first row than second row, so we can insert $n$ into the second row and rearrange the row entries to obtain an inversion-free filling $\sigma$ of shape $\mu$ and content $\alpha$.  We define $\phi(\sigma')=\sigma$, and by Theorem \ref{RectangleBump} applied to the two-row rectangle inside $\mu$ of width equal to the top row of $\mu$, the major index increases by $1$ from $\sigma'$ to $\sigma$.  Thus $\phi$ is an inverse to $\psi$ on $\fancy{F}_{\mu^{(1)}}^{(1^{n-1})}|_{\inv=0}$, and $\psi$ is bijective.
       \end{proof}
       
       We now complete the entire bijection for two rows by defining a standardization map for two-row fillings.
       
       \begin{definition}
         For a two-row shape $\mu=(\mu_1,\mu_2)$, we define the map $$\Standardize:\fancy{F}_{\mu}^\alpha|_{\inv=0}\to \fancy{F}_{\mu}^{1^n}|_{\inv=0}$$ as follows.  Given a filling $\sigma\in \fancy{F}_{\mu}^\alpha|_{\inv=0}$, define $\Standardize(\sigma)$ to be the filling of $\mu$ with content $(1^n)$ that respects the ordering of the entries of $\sigma$ by size, with ties broken by reading order.
       \end{definition}
    
       \begin{example}
          The standardization map for two rows is illustrated below.
           \begin{center}\includegraphics{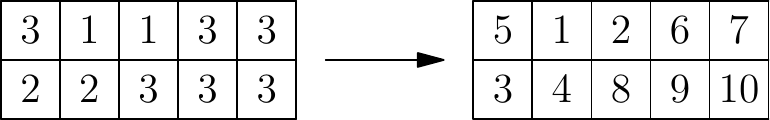} \end{center}
       \end{example}
       We can now define $\majcode$ for two-row shapes.
       
       \begin{definition}
          Let $\mu=(\mu_1,\mu_2)$ be a two-row shape of size $n$.  Given a filling $\sigma$ of $\mu$, let $\overline{\sigma}=\Standardize(\sigma)$.  Then we define $\majcode(\sigma)=d_1d_2\cdots d_n$ where $$d_i=\maj(\psi^{i-1}(\overline{\sigma}))-\maj(\psi^{i}(\overline{\sigma})),$$
          and where $\psi$ is the map defined in Lemma \ref{tworows}.
       \end{definition}
       
       \begin{remark}
          Notice that, given a filling $\sigma$ of $\mu$ having arbitrary content, we have $$\majcode(\sigma)=\majcode(\Standardize(\sigma)).$$ 
       \end{remark}
       
       \begin{theorem}\label{tworowsconclusion}
          The map $\majcode$ defined on two-row shapes $\mu=(\mu_1,\mu_2)$ is an isomorphism of weighted sets $$\fancy{F}_{\mu}^\alpha|_{\inv=0}\to C_{\mu,A}$$ for each alphabet $A$ and corresponding content $\alpha$.  
       \end{theorem}
       
       \begin{proof}
          Putting together the recursions of Lemma \ref{tworows} and Lemma \ref{recursion}, we have that for the content $(1^n)$ corresponding to alphabet $[n]$, the map $\majcode$ is a weighted set isomorphism $$\fancy{F}_{\mu}^{(1^n)}|_{\inv=0}\to C_{\mu,[n]}.$$
          
          Now, let $A$ be any alphabet with content $\alpha$.  Let $\sigma$ be a filling of $\mu$ with content $\alpha$.  Then we know $\majcode(\sigma)=\majcode(\Standardize(\sigma))$, so $\majcode(\sigma)\in C_{\mu,[n]}$.  In other words, $\majcode(\sigma)$ is $\mu$-sub-Yamanouchi.  In addition, since $\Standardize$ is an injective map (there is clearly only one way to un-standardize a standard filling to obtain a filling with a given alphabet), the map $\majcode$, being a composition of $\Standardize$ and the $\majcode$ for standard fillings, is injective as well on fillings with content $\alpha$.
          
          We now wish to show that $\majcode(\sigma)=d_1,\ldots,d_n$ is $A$-weakly increasing, implying that $\majcode$ is an injective morphism of weighted sets to $C_{\mu,A}$.  To check this, let $\widetilde{\sigma}=\Standardize(\sigma)$.  Then any repeated letter from $\sigma$ will become a collection of squares that have consecutive entries and are increasing in reading order in $\widetilde{\sigma}$.  Neither of the two operations of the map $\psi$ affects the reading order of such subcollections, since consecutive integers $a$ and $a+1$ cannot occur in reverse order in a filling with distinct entries and no inversions.  So, it suffices to show that if the largest entry $m$ of $\sigma$ occurs $i$ times, then $d_1\le \cdots \le d_i$.   
          
          In $\widetilde{\sigma}$, the $m$'s of $\sigma$ become the numbers $n-i+1, n-i+2,\ldots, n$, and occur in reading order.  Thus we remove any of these that occur in the bottom row first, and for those we have $d_t=0$.  We continue removing these from the bottom row until there are none left in the bottom row.  Then the remaining $d_t$'s up to $d_i$ will equal $1$.  Therefore, $d_1\le d_2\le \cdots \le d_i$, as required.
          
          Finally, the number of fillings with content $\alpha$ is the same as the number of cocharge-friendly diagrams, which is the same as the number of inversion-friendly dot-diagrams for the reverse alphabet.  This in turn is the same as the cardinality of $C_{\mu,A}$ by the section on inversion codes above.   Thus the injective map $\majcode$ is in fact a bijection.  The result follows.
       \end{proof}

       \begin{corollary}
         For any two-row shape $\mu$ and content $\alpha$, the map $\invcode^{-1}\circ\majcode$ is an isomorphism of weighted sets from $\fancy{F}_{\mu}^{\alpha}|_{\inv=0}\to \fancy{F}_{\mu^\ast}^{r(\alpha)}|_{\maj=0}$.  This gives a combinatorial proof of the identity $$\widetilde{H}_{\mu}(x;0,t)=\widetilde{H}_{\mu^\ast}(x;t,0)$$ for two-row shapes.
       \end{corollary}
       
       \begin{example}
         In Figure \ref{tworowexample}, the map $\majcode$ is applied to a two-row filling $\sigma$.
         \begin{figure}
         \begin{center}
                   \includegraphics{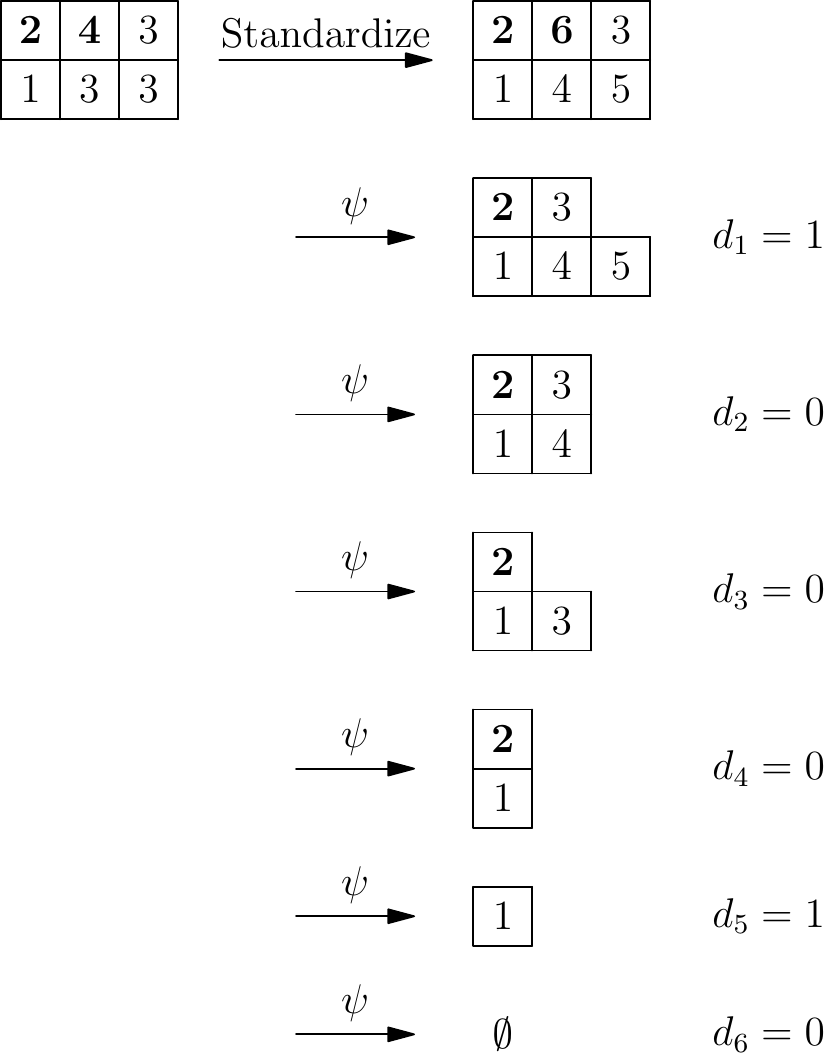}
         \end{center}
         \caption{The map $\majcode$ for two-row shapes.}\label{tworowexample}
         \end{figure}
          The figure shows that $\majcode(\sigma)=100010$.  If we apply $\invcode^{-1}$ to this code using the reversed alphabet, we obtain the the filling $\rho$ below: 
         \begin{center}
          \includegraphics{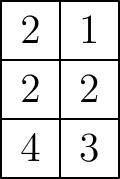}
         \end{center}
            Notice that $\maj(\sigma)=\inv(\rho)=2$.
       \end{example}
       
       We now have the tools to extend our map $\psi$ to three-row shapes.
      
      \begin{definition}
        Let $\sigma$ be any filling of a three-row shape $\mu=(\mu_1,\mu_2,\mu_3)$, and let $\sigma'$ be the $3\times \mu_3$ rectangle contained in $\sigma$.  Let $n$ be the largest entry in $\sigma$.  Choosing one of these $n$'s, say $n_i$, we define $\psi_{n_i}(\sigma)$ by the following process.
         \begin{enumerate}
            \item If $n_i$ is to the right of $\sigma'$, remove the $n$ as in the two-row algorithm to form $\psi_{n_i}(\sigma)$.
            \item If $n_i$ is in the bottom row and in $\sigma'$, then $\sigma$ is a rectangle and we let $\psi_{n_i}(\sigma)=\sigma_\downarrow$.
            \item If $n_i$ is in the second row and in $\sigma'$, let $a_2$ be the top entry of the bumping sequence of $\sigma'$.  Let $b$ be the entry in square $(\mu_2+1,2)$ if it exists, and let $b=n+1$ otherwise.   If $b\ge a_2$, then remove $n_i$ and bump down $a_2$ to the second row, and if $b<a_2$, simply remove $n_i$.  Rearrange the modified rows so that there are no inversions, and let $\psi_{n_i}(\sigma)$ be the resulting filling. 
            \item If $n_i$ is in the top row and in $\sigma'$, let $a_1,a_2$ be the bumping sequence of the $3\times \mu_3$ rectangle in $\sigma$.  If $a_2>a_1$ or $\mu_2=\mu_3$, then remove $n_i$ from $\sigma$.  Otherwise, if $a_2\le a_1$, remove $n$ and bump $a_2$ up to the top row.  Rearrange the modified rows so that there are no inversions, and let $\psi_{n_i}(\sigma)$ to be the resulting filling.
          \end{enumerate}
      \end{definition}

       \begin{lemma}\label{threerows}
         Let $\mu=(\mu_1,\mu_2,\mu_3)$ be any three-row shape of size $n$.  Then the map $\psi=\psi_n$ defined above is a morphism of weighted sets when restricted to fillings having distinct entries.  That is, in the case of distinct entries there is a unique choice of $n$, and $$\psi:(\fancy{F}_{\mu}^{(1^n)}|_{\inv=0};\maj)\to \bigsqcup_{d=0}^{2}(\fancy{F}_{\mu^{(d+1)}}^{(1^{n-1})}|_{\inv=0};\maj+d)$$ is a morphism of weighted sets.
       \end{lemma}
       
       We defer the proof to Section \ref{Proofs}.  In that section, we also show:
       
       \begin{lemma}\label{threerowsbijective}
         The map $\psi$ of Lemma \ref{threerows} is an isomorphism.
       \end{lemma}
       
       We can now complete the three-row case by defining its standardization map for fillings with repeated entries.  This definition is designed to force the $\majcode$ sequences to be $A$-weakly increasing.
       
       \begin{definition}
         Given a filling $\sigma$ of $\mu$, define $\Standardize(\sigma)$ as follows.  First, for any letter $i$ that occurs with multiplicity in $\sigma$, label the $i$'s with subscripts in reading order to distinguish them.  If we bump one of them up or down one row, choose the one to bump from the row in question that preserves their reading order.
         
         Let $n$ be the largest entry that occurs in $\sigma$.  For each such $n_t$ compute $d_t=\maj(\sigma)-\maj(\psi_{n_i}(\sigma))$, and let $d=\min_t(\{d_t\})$.  Let $n_r$ be the last $n$ in reading order for which $d_r=d$.  Form the filling $\psi_{n_r}(\sigma)$, and repeat the process on the new filling.  Once there are no $n$'s left to remove, similarly remove the $n-1$'s, and so on until the empty tableau is reached.
         
         Now, consider the order in which we removed the entries of $\sigma$ and change the corresponding entries to $N,N-1,\ldots,1$ in that order, where $N=|\mu|$.  The resulting tableau is $\Standardize(\sigma)$.
         
       \end{definition}
       
       We can now define $\majcode$ for three-row shapes.
              
      \begin{definition}
          Let $\mu=(\mu_1,\mu_2,\mu_3)$ be a three-row shape of size $n$.  Given a filling $\sigma$ of $\mu$, let $\overline{\sigma}=\Standardize(\sigma)$.  Then we define $\majcode(\sigma)=d_1d_2\cdots d_n$ where $$d_i=\maj(\psi^{i-1}(\overline{\sigma}))-\maj(\psi^{i}(\overline{\sigma})),$$
          and where $\psi$ is the map defined in Lemma \ref{threerows}.
      \end{definition}
              
       \begin{remark}
           Notice that, given a filling $\sigma$ of $\mu$ having arbitrary content, we have $$\majcode(\sigma)=\majcode(\Standardize(\sigma)).$$ 
       \end{remark}
       
      \begin{theorem}\label{threerowsconclusion}
          The map $\majcode$ defined on three-row shapes $\mu=(\mu_1,\mu_2)$ is an isomorphism of weighted sets $$\fancy{F}_{\mu}^\alpha|_{\inv=0}\to C_{\mu,A}$$ for each alphabet $A$ and corresponding content $\alpha$.  
      \end{theorem}
      
      See section \ref{Proofs} for the proof.
      
      \begin{corollary}
         For any three-row shape $\mu$ and content $\alpha$, the map $\invcode^{-1}\circ\majcode$ is an isomorphism of weighted sets from $\fancy{F}_{\mu}^{\alpha}|_{\inv=0}\to \fancy{F}_{\mu^\ast}^{r(\alpha)}|_{\maj=0}$.  This gives a combinatorial proof of the identity $$\widetilde{H}_{\mu}(x;0,t)=\widetilde{H}_{\mu^\ast}(x;t,0)$$ for two-row shapes.
      \end{corollary}
    
    \begin{example}
      We demonstrate all of the above maps on the filling $\sigma$ below, with its repeated entries labeled with subscripts in reading order to distinguish them.
      \begin{center}
         \includegraphics{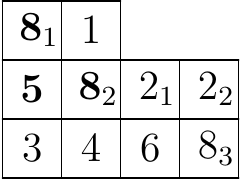}
      \end{center}
      We will standardize and compute $\majcode$ simultaneously.  To decide which of the $8$'s to remove first, we look at which would give the smallest first $\majcode$.  This is clearly the $8_3$ in the bottom row, so we remove it and bump down the $2$.
      \begin{center}
         \includegraphics{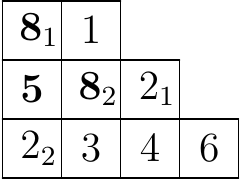}
      \end{center}
      To decide which of the remaining $8$'s to remove next, note that they both would decrease $\maj$ by $2$, and so we remove the one that comes last in reading order, namely $8_2$.  Since $1<2$ we bump down the $1$.
      \begin{center}
         \includegraphics{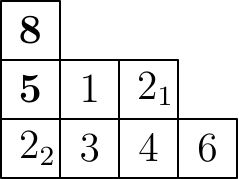}
      \end{center}
      Finally, when we remove the last $8$, the maj decreases by $2$, so we do not have to lift the $5$ up to the third row.
      \begin{center}
         \includegraphics{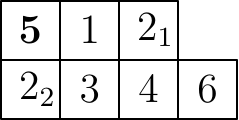}
      \end{center}
      We can now use the two-row algorithm to complete the process, and we find $\majcode(\sigma)=0220100000$.  The corresponding inversion diagram for the reverse alphabet $\{1,1,1,3,4, 5, 6, 7, 7, 8\}$ is shown below.
      \begin{center}
         \includegraphics{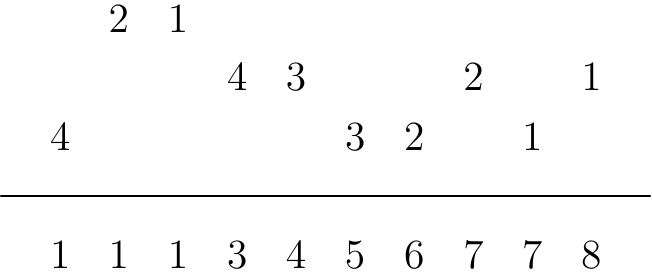}
      \end{center}
      Finally, we can reconstruct from this the filling $\rho=\invcode^{-1}(0220100000)$ below.
      \begin{center}
         \includegraphics{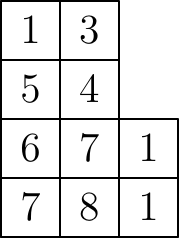}
      \end{center}
      
      Note that $\inv(\rho)=\maj(\sigma)=5$, and $\inv(\sigma)=\maj(\rho)=0$.
    \end{example}
       
    \begin{remark}
      The map above essentially uses the fact that a three-row shape is the union of a rectangle and a two-row shape.  For \textit{any} shape that is the union of a rectangle and two rows, a similar map can be used to remove the first $n$, and so for ``snorkel'' shapes consisting of two rows plus a long column, a similar algorithm also produces a valid $\majcode$ map.
      \begin{center}
      \includegraphics{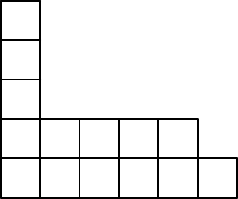}
      \end{center}
      
      In general, however, the resulting shape on removing the first $n$ is no longer the union of a rectangle and a two-row shape, and we cannot use an induction hypothesis.
      
      However, we believe that this method may generalize to all shapes, as follows.
    \end{remark}
    
    \begin{conjecture}
      One can extend the map $\psi$ for three-row shapes to all shapes inductively, as follows.  One would first extend it to shapes which are the union of a three-row shape and a single column, then use this to extend it to shapes which are the union of a three-row shape and a rectangle shape, using Theorem \ref{RectBumpCor2}.  One could then iterate this new map on any four-row shape, so $\majcode$ can then be defined on four-row shapes, and so on.
    \end{conjecture}

    \section{Application to Cocharge}\label{Applications}
    
    Proposition \ref{ZeroBump} reveals an interesting property of the cocharge statistic on words, first defined by Lascoux and Sch\"{u}tzenberger.  To define it, we first recall the definition of Knuth equivalence.  
    
    \begin{definition}
    Given a word $w=w_1\cdots w_n$ of positive integers, a \textit{Knuth move} consists of either:
    \begin{itemize}
    \item A transposition of the form $xyz\to xzy$ where $x,y,z$ are consecutive letters and $y<x\le z$ or $z<x\le y$
    \item A transposition of the form $xyz\to yxz$ where $x,y,z$ are consecutive letters and $x\le z<y$ or $y\le z<x$.
    \end{itemize} 
    Two words $w,\widetilde{w}$ are said to be \textit{Knuth equivalent}, written $w~\widetilde{w}$, if one can be reached from the other via a sequence of Knuth moves.  Knuth equivalence is an equivalence relation on words.
    \end{definition}
    
    Cocharge was originally defined as follows.
    
    \begin{definition}\label{OriginalCocharge}
    Given a word $w=w_1,\cdots, w_n$ with partition content $\mu$, the \textit{cocharge} of $w$, denoted $\cc(w)$ is the unique statistic satisfying the following properties:
    \begin{enumerate}
      \item It is constant on Knuth equivalence classes, that is, if $w$ is Knuth equivalent to $\widetilde{w}$ then $\cc{w}=\cc{\widetilde{w}}$. 
      \item If $w=w_1w_2\cdots w_n$ and $w\neq 1$, let $\cyc(w)=w_2w_3\cdots w_nw_1$ be the word formed by moving the first letter to the end.  Then $$\cc(\cyc(w))=\cc(w)-1.$$
      \item If the letters of $w$ are weakly increasing then $\cc(w)=0$.
    \end{enumerate}
    \end{definition}
    
    There is also an algorithmic way of computing cocharge.
    
    \begin{definition}\label{CochargeAlgorithm}
      Let $w$ be a word with partition content $\mu$, so that it has $\mu_1$ $1$'s, $\mu_2$ $2$'s, and so on.  Let $w^{(1)}$ be the subword formed by scanning $w$ from right to left until finding the first $1$, then continuing to scan until finding a $2$, and so on, wrapping around cyclically if need be.  Let $w^{(2)}$ be the subword formed by removing $w^{(1)}$ from $w$ and performing the same process on the remaining word, and in general define $w^{(i)}$ similarly for $i=1,\ldots,\mu_1$. 
      
    It turns out that $$\cc(w)=\sum_i \cc(w^{(i)}),$$ (see, e.g., \cite{HHL}) and one can compute the cocharge of a word $w^{(i)}$ having distinct entries $1,\ldots,k$ by the following process.
    
    \begin{enumerate}
    \item Set a counter to be $0$, and label the $1$ in the word with this counter, i.e. give it a subscript of $0$. 
    \item If the $2$ in the word is to the left of the $1$, increment the counter by $1$, and otherwise do not change the counter.  Label the $2$ with the new value of the counter.
    \item Continue this process on each successive integer up to $k$, incrementing the counter if it is to the left of the previous letter.
    \item When all entries are labeled, the sum of the subscripts is the cocharge.
    \end{enumerate}
    \end{definition}
    
    The link between the major index of inversion free fillings and the cocharge of words lies in the \textit{cocharge word} construction.
    
    \begin{definition}
      The \textit{cocharge word} of a filling $\sigma:\mu\to \ZZ_+$ is the word $\cw(\sigma)=i_1 i_2\cdots i_n$ consisting of the row indices of the cells $u_k=(i_k,j_k)$, where $u_1,u_2,\ldots,u_n$ is the ordering of the cells of $\mu$ such that $\sigma(u_1)\ge \sigma(u_2)\ge\cdots \ge \sigma(u_n)$, and for each constant segment $\sigma(u_j)=\cdots =\sigma(u_k)$, the cells $u_j,\cdots,u_k$ are in reverse reading order.
    \end{definition}
    
    As mentioned in Section \ref{CochargeWords}, for any filling $\sigma\in \mathcal{F}|_{\inv=0}$ we have $\maj(\sigma)=\cc(\cw(\sigma))$.  (See \cite{HHL} for the proof.)  Therefore, we can translate some of our results regarding such fillings to properties of words and their cocharge.  We first require the following fact.
    
    \begin{lemma}\label{Columns}
    If $\sigma\in \mathcal{F}|_{\inv=0}$ and $w=\cw(\sigma)$, the words $w^{(i)}$ correspond to the columns of $\sigma$, in the sense that the letters in the subword $w^{(i)}$ are in positions corresponding to the entries in column $i$ in $\sigma$.
    \end{lemma}
    
    \begin{proof}
      If $w=\cw(\sigma)$ and $\sigma$ has alphabet $A=a_1\ge \ldots\ge a_n$, the letters $a_i$ for which the corresponding letter $w_i$ equals $r$ are the entries in row $r$.  The smallest - that is rightmost - $a_{i}$, say $a_{i_0}$, for which $w_{i}=1$ is the leftmost entry of the bottom row, i.e. the bottom entry of the first column.  The second entry of the first column is then the first $a_i$ in cyclic order after $a_{i_0}$ for which $w_i=2$.  This corresponds to the $2$ in the subword $w^{(1)}$, and similarly the letters in $w^{(1)}$ correspond to the entries in the first column.  
      
      A similar argument shows that the second column corresponds to $w^{(2)}$, and so on.
    \end{proof} 
    
    In particular, Proposition \ref{ZeroBump} states that if the largest entry of a filling $\sigma \in \mathcal{F}_\mu^{(1^n)}|_{\inv=0}$ is in the bottom row, then we can remove it, bump down any entries in its (rightmost) column, and rearrange the rows to get a filling with no inversions.  By Lemma \ref{Columns}, this translates to the following result in terms of words.
    
    \begin{theorem}
      Let $w=w_1\cdots w_n$ be a word with partition content $\mu$ for which $w_1=1$.  Let $w^{(1)},\cdots,w^{(\mu_1)}$ be its decomposition into subwords as in Definition \ref{CochargeAlgorithm}.  Then $w_1\in w^{(\mu_1)}$, and if $w'$ is the word formed by removing $w_1$ from $w$ and also decreasing each letter that is in $w^{(\mu_1)}$ by one, then $$\cc(w)=\cc(w').$$
    \end{theorem}
    
    This theorem fills a gap in our understanding of cocharge, as it gives a recursive way of dealing with words that start with $1$.  These are the only words that do not satisfy the relation $\cc(\cyc(w))=\cc(w)-1$ of Definition \ref{OriginalCocharge}.
    
    \begin{example}
    Consider the word $15221432313$.  It has three $1$'s, three $2$'s, and three $3$'s, but only one $4$ and $5$, so to find the word $w^{(\mu_{1})}=w^{(3)}$ we can ignore the $4$ and $5$.  The words $w^{(1)}$, $w^{(2)}$, and $w^{(3)}$, ignoring the $4$ and $5$, are the subwords listed below:
    $$\begin{array}{c|ccccccccccc}
    w       & 1 & 5 & 2 & 2 & 1 & 4 & 3 & 2 & 3 & 1 & 3 \\\hline
    w^{(1)} &   &   &   &   &   &   & 3 & 2 &   & 1 &   \\\hline
    w^{(2)} &   &   &   & 2 & 1 &   &   &   &   &   & 3 \\\hline
    w^{(3)} & 1 &   & 2 &   &   &   &   &   & 3 &   &   \\
    \end{array}$$
    and so the word $w'$ is formed by removing the leading $1$ and decreasing the $2$ and $3$ from $w^{(3)}$.  Thus $$w'=5121432213.$$  We also find that $\cc(w)=\cc(w')=12$. 
    \end{example}
    
    \section{Technical proofs}\label{Proofs}
    
    This section contains the proofs of all the results above whose proofs are particularly long or technical.
    
    \subsection{Proof of Proposition \ref{recursion}: The Recursion}

    \begin{GeneralRecursion}[General Recursion]
      For any partition $\mu$ of $n$ and any multiset of positive integers $A=\{a_1,a_2,\ldots,a_n\}$ with $a_1\le a_2\le \cdots \le a_n$, we have $$\fancy{C}_{A}(\mu)=\bigsqcup_{i=1}^{\mu_1^\ast}x_n^{i-1}\cdot \fancy{C}^{(i-1)}_{A}(\mu^{(i)}).$$
    \end{GeneralRecursion}

    \begin{proof}
      The sets forming the union on the right hand side are disjoint because the $i$th set consists only of monomials having $x_n^{i-1}$ as their power of $x_n$.  We now show inclusion both ways.
        
      ($\subseteq$)  Let $x^c\in \fancy{C}_A(\mu)$ where $c=c_1,\ldots,c_n$ is a $\mu$-sub-Yamanouchi word which is $A$-weakly increasing.  Let $i=c_1+1$, so that $c_1=i-1$.  Also let $c'=c_2,\ldots,c_n$.  Notice that if $a_1=a_2$ then $c_2\ge i-1$, and $c'$ is $A\setminus \{a_1\}$-weakly increasing.  Thus, to show $x^c\in x^{i-1}\fancy{C}^{(i-1)}_A(\mu^{(i)})$, we just need to show that $c'$ is $\mu^{(i)}$-sub-Yamanouchi.
        
      Since $c$ is $\mu$-sub-Yamanouchi, there exists a Yamanouchi word $d$ having $\mu_i$ entries equal to $i-1$ for each $i$, for which $x^c | x^d$.  Let $t$ be the highest index such that $\mu_{t+1}=\mu_i$.  Then $\mu^{(i)}=(\mu_1,\mu_2,\ldots,\mu_t-1,\cdots,\mu_k)$.  So, we wish to show that we can form a new $\mu$-Yamanouchi word $b$ from $d$ so that we still have $x^c|x^b$ but $b_1=t$.  This way $c'$ will be $\mu^{(i)}$-sub-Yamanouchi, with respect to $b'=b_2,\ldots,b_n$.
      
      We have $\mu_{t+2}<\mu_{t+1}$ by our assumption defining $t$, so there are strictly more $t$'s than $t+1$'s in $d$.  Notice that this means we can move the \textit{leftmost} $t$ in $d$  any number of spots to the left without changing the fact that the word is Yamanouchi.
        
      Also notice that $d_1\ge c_1=i-1$.  But since there are exactly as many $i-1$'s as $i$'s, $i+1$'s, and so on up to $t$ in $d$, we must in fact have $d_1\ge t$, for otherwise the suffix $d_2,\ldots,d_n$ would not satisfy the Yamanouchi property.  So $d_1\ge t$.
      
      Now, let $d_r$ be the leftmost $t$ in $d$.  We form a subword of $d$ as follows.  Let $d_1$ be the first letter of our subword.  Then let $d_{p_1}$ be the leftmost letter between $d_1$ and $d_r$ with $t\le d_{p_1}\le d_r$, if it exists.  Then let $d_{p_2}$ be the first letter between $d_{p_1}$ and $d_r$ for which $t\le d_{p_2}\le d_{p_1}$, and so on until we reach a point at which no such letter exists.  We now have a subsequence of letters $d_1, d_{p_1}, d_{p_2},\ldots,d_{p_k},d_r=t$ where $d_r$ is the leftmost $t$ in $d$.   We define $b$ to be the word formed from $d$ by cyclically shifting this subsequence, replacing $d_{p_i}$ with $d_{p_{i-1}}$ for all $i>1$, replacing $d_{p_1}$ with $d_1$, and replacing $d_1$ with $d_{p_k}$.
        
      For instance, if $\mu=(4,3,3,2,2)$, $i-1=1$, then $t=2$, and we might have $$c=120412130010100$$ with $$d=\mathbf{4}\mathbf{3}04\mathbf{2}2130021100.$$  Then the subword of $d$ consists of those letters in boldface above, and we cyclically shift the boldface letters to the right in their positions to form $$b=\mathbf{2}\mathbf{4}04\mathbf{3}2130021100,$$ which is still $\mu$-Yamanouchi and still dominates $c$ in the sense that $x^c|x^b$.
        
      To verify that in general $x^c|x^b$, notice that $c_1=i-1\le t$, and since the other letters in the subword decrease to the right, we have $b_i\ge d_i$ for all $i>1$.  Thus each $b_i\ge c_i$ for all $i$, and so $x^c|x^b$.
        
      To show that $b$ is still Yamanouchi, notice that to form $b$ from $d$, we have moved the leftmost $t$ all the way to the left (which, we noted above, preserves the Yamanouchi property) and moved each $d_{p_j}$ to the right without crossing over any element having value $d_{p_j}-1$ (for otherwise our sequence $d_{p_j}$ would have an extra element, a contradiction.)  Thus we have not changed the property of there being at least as many $d_{p_j}-1$'s as $d_{p_j}$'s in each suffix, and we have not changed the property that there are at least as many $d_{p_j}$'s as $d_{p_j}+1$'s in each suffix, because we moved these elements to the right.  The other Yamanouchi conditions remain unchanged, since we are only moving the letters $d_{p_j}$.  Thus $b$ is Yamanouchi as well.
        
      ($\supseteq$) For the other inclusion, let $c=c_1,\ldots,c_n$ be a word such that $x^c\in x^{i-1}\cdot \fancy{C}^{(i)}_A(\mu^{(i)})$.  Then $c'=c_{2},\ldots,c_n$ is $\mu^{(i)}$-sub-Yamanouchi, so there exists a word $d'=d_2,\ldots,d_n$ which is Yamanouchi of content $\mu^{(i)}$ such that $x^{c'}|x^{d'}$.  Let $d_1=t$ where $t$ is the highest index such that $\mu_{t+1}=\mu_i$.  Then $d=d_1,\ldots,d_n$ is Yamanouchi of shape $\mu$ by the definition of $\mu^{(i)}$, and since $c_1=i-1$, we have $c_1\le t=d_1$.  Thus $x^c|x^d$.  Finally, note that if $a_1=a_2$ in $A$, then $c_2\ge i-1$ by the definition of $\fancy{C}^{(i)}$.  Thus $c$ is $A$-weakly increasing. It follows that $x^c\in \fancy{C}_A(\mu)$.
      \end{proof}
    
    \subsection{Proof of Theorem \ref{Invcode}: The Isomorphism $\invcode$}
     
   \begin{InvcodeTheorem}
   The inversion code of any filling $\rho\in \fancy{F}_{\mu^\ast}^{\alpha}$ is $\alpha$-weakly increasing and $\mu$-sub-Yamanouchi.  Moreover, the map $$\invcode:\fancy{F}_{\mu^\ast}^\alpha\mid_{\maj=0}\to C_{\mu,A}$$ is an isomorphism of weighted sets.
   \end{InvcodeTheorem}
    
   To prove Theorem \ref{Invcode}, we first introduce some new notation. In analogy with the cocharge word defined in \cite{HHL}, for fillings $\rho$ having $\maj(\rho)=0$, we can form an associated \textit{inversion word} and describe a statistic on the inversion word that measures $\inv(\rho)$ in the case that $\maj(\rho)=0$.
   
   \begin{definition}
   Let $\rho$ be a filling of shape $\mu$ having $\maj=0$.  We define the \textit{inversion word} of $\rho$ as follows.  Starting with the smallest value that appears in the filling, write the column numbers of the entries with that value as they appear in reading order, and then proceed with the second largest entry and so on.  
   \end{definition}
   
   For instance, the filling:
   \begin{center}
   \includegraphics{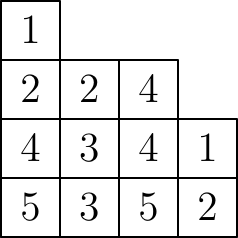}
   \end{center}
   has inversion word $141242231313$.
   
   In order to compute $\inv(\sigma)$ given only its inversion word, we will use a visual representation of the inversion word, which we call a \textit{diagram}.
   
   \begin{definition}
     Fix a linearly ordered finite multiset $A$, with elements $a_1\le a_2\le \cdots \le a_n$.  The \textit{diagram} a function $f:A\to \ZZ_+$ is the plot of the function with respect to the ordering on $A$.  We say that the diagram has \textit{shape} $\mu$ if $|f^{-1}(i)|=\mu_i$ for each $i$. 
   \end{definition}
   
   The diagrams we will be using are essentially the plot of the inversion word, considered as a function on a multiset.
   
   \begin{definition}  
   Let $\rho$ be a filling of $\mu^\ast$ having $\maj(\rho)=0$, and let $w$ be the inversion word of $\rho$.  Let $A$ be the multiset consisting of the entries of $\rho$, ordered from least to greatest and in reading order in the case of a tie.  Let $f:A\to \ZZ_+$ be the function given by $f(a_i)=w_i$.  We define $\InvPlot(\rho)$ to be the diagram of the function $f$, whose plot has $\mu_j$ dots in the $j$th row.  
   \end{definition}
   
   Notice that the $\InvPlot$ of a filling of shape $\mu^\ast$ has shape $\mu$, the conjugate shape.  For instance, the tableau
   \begin{center}
   \includegraphics{pic4-eps-converted-to.pdf}
   \end{center}
   has $\maj=0$, and its inversion word is $11213122$.  Its plot is as follows.
   \begin{center}
   \includegraphics{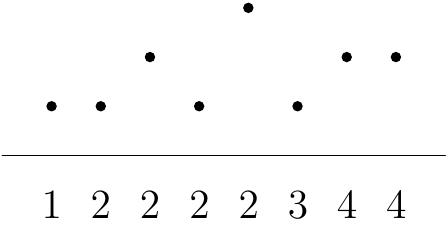}
   \end{center}
   
    To compute the number of inversions, we define the \textit{inversion labeling} of a diagram to be the result of labeling each row of dots $\mu_i$ in the diagram with the numbers $1,2,\ldots,\mu_i$ from right to left:
    
   \begin{center}
   \includegraphics{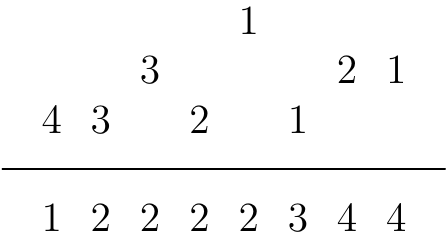}
   \end{center}
   Finally, an \textit{inversion} in the diagram of a function $f:A\to \ZZ_+$, labeled as above, is a pair of entries $a<b$ in the ordered multiset $A$ for which either:
   \begin{itemize}
   \item[I.] The dots above $a$ and $b$ have the same label and $f(a)>f(b)$, or
   \item[II.] The dot in position $a$ is labeled $i$ and the dot in position $b$ is labeled $i+1$, and $f(b)>f(a)$.
   \end{itemize}
   So there are $3$ inversions in the diagram above, two of type I and one of type II:
   
   \begin{center}
   \includegraphics{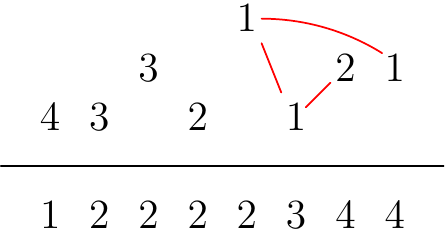}
   \end{center}
   
    For fillings $\sigma$ with $\maj(\sigma)=0$, there are no descents, and so the number of inversions in $\InvPlot(\sigma)$ is equal to $\inv(\sigma)$.  In particular, type I and II inversions correspond to attacking pairs in the same row or on adjacent rows, respectively.
    
    \begin{remark}
     The type I and II inversions also correspond to the two types of inversions used to define the $\dinv$ statistic on parking functions.  Indeed, this was the original motivation for the full definition of the $\inv$ statistic.  \cite{HaglundGenesis}
    \end{remark}
   
   We now classify the types of diagrams that arise as the $\InvPlot$ of a filling.
   
   \begin{definition}
      A consecutive subsequence is in \textit{inversion-friendly order} if, when each row is labeled from right to left as above, all dots of label $i+1$ in the subsequence occur before the dots of label $i$ for all $i$, and the dots of any given label appear in increasing order from bottom to top.
   \end{definition}
   
   An example of an inversion-friendly subsequence is shown below.
   \begin{center}
    \includegraphics{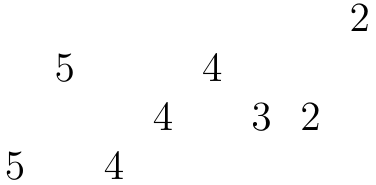}
   \end{center}
   
   It is easy to check that, in the plot of any filling $\rho$ having $\maj\rho=0$, every subsequence above a fixed letter of the alphabet $A$ is in inversion-friendly order.  We claim that the converse is true as well, namely, that every diagram having all such subsequences in inversion-friendly order corresponds to a unique Young diagram filling $\rho$ having $\maj(\rho)=0$.
   
   \begin{definition} 
     A diagram is of \textit{inversion word type} if every subsequence determined by a fixed letter of $A$ is in inversion-friendly order.
     
     We let $\ID_{\mu,A}$ the set of all diagrams of shape $\mu$ of inversion word type over $A$.  We equip $\ID_{\mu,A}$ with its $\inv$ statistic to make it into a weighted set.
   \end{definition}
   
   \begin{proposition}\label{InversionFriendly}
     Let $\mu$ be a partition of $n$, and let $A$ be a multiset of $n$ positive integers with content $\alpha$.  The map $\InvPlot$ is an isomorphism of weighted sets $$\InvPlot:(\fancy{F}_{\mu^\ast}^\alpha|_{\maj=0};\inv)\to (\ID_{\mu,A};\inv).$$ 
   \end{proposition}
   
   \begin{proof} (Sketch.)
     As noted above, this is a map of sets that preserves the $\inv$ statistic since there are no descents.  To show it is bijective, we construct its inverse.
     
     Let $D$ be an arbitrary diagram in $\ID_{\mu,A}$, and let $f:A\to \ZZ_+$ be the corresponding map.  For any $a\in A$ let $\ell(a)$ be the label on the dot at height $f(a)$.  Then let $\rho$ be the filling of shape $\mu^\ast$ in which $a\in A$ is placed in the square in column $f(a)$ from the left, and height $\ell(a)$ from the bottom.  By the definition of $\InvPlot$, we have that $\InvPlot(\rho)=D$, and furthermore if $D=\InvPlot(\sigma)$ then $\rho=\sigma$.  Thus the map sending $D$ to $\rho$ is the inverse of $\InvPlot$.
   \end{proof}
    
    We now show that the inversion-friendly diagrams are in weight-preserving bijection with generalized Carlitz codes.  This will complete the proof of Theorem \ref{Invcode}.

    \begin{definition}
    The \textit{inversion code} of a diagram $w$, denoted $\invcode(w)$ is the sequence $\{c_i\}$ whose $i$th entry $c_i$ is the number of inversion pairs of the form $(w(i),b)$.
    \end{definition}
      
    \begin{example}
     The inversion code of the following diagram is $00002100$.
     \begin{center}
     \includegraphics{pic6-1.pdf}
      \end{center}
    \end{example}
      
     Using Proposition \ref{InversionFriendly}, we can also define the inversion code of a filling $\rho\in \fancy{F}_\mu^\alpha|_{\maj=0}$ to be $$\invcode(\rho):=\invcode(\InvPlot(\rho)).$$  It is easy to see that this matches the definition of inversion code in Section \ref{CochargeWords}.
    
       \begin{theorem}
         The map $\invcode:\ID_{\mu,A}\to C_{\mu,A}$ is an isomorphism of weighted sets.
       \end{theorem}
    
      We break the proof into several lemmas for clarity.
      
      \begin{lemma}\label{welldefinedinv}
      The map $\invcode$ is a well-defined morphism from $\ID_{\mu,A}\to C_{\mu,A}$ for all $\mu$ and $A$.
      \end{lemma}
      
      \begin{proof} 
       Let $w:A\to \ZZ_+$ be a diagram in $\ID_{\mu,A}$, and let $c=\invcode(w)$. 
       
       We first show that $c$ is $\mu$-sub-Yamanouchi.  Let $i>0$ and consider the subset of dots labeled $i$ in the inversion labeling of $w$, say $w(r_1),\ldots,w(r_t)$ from left to right.  We claim that $w(r_{t-j})$ is the left element of at most $j$ inversions for each $j=0,\ldots,t-1$.  Indeed, $w(r_{t-j})$ is to the left of exactly $j$ dots labeled $i$; those dots in a lower row form the Type I inversions with $w(r_{t-j})$.  For Type II, the dots labeled $i+1$ in a higher row must have an $i$ to the right of them, so correspond to one of the dots labeled $i$ in a higher row and to the right of $w(r_{t-j})$.  Thus $w(r_{t-j})$ is the left element of at most $j$ inversions, and so $c_{r_{t-j}}\le j$.
       
       It follows that $c_{r_1},\ldots,c_{r_t}$ is an ordinary Carlitz code.  Therefore, $c$ can be decomposed into several Carlitz codes, one for each label, of lengths $\mu_1^\ast,\mu_2^\ast,\ldots$.  Let $d_i$ be the resulting upper bound on $c_i$ for each $i$.  Then $d$ is a union of the sequences $$\mu_i^\ast,\mu_i^\ast-1,\ldots,2,1,0$$ for each $i$, arranged so that each of these sequences retains its order.  Thus $d$ is a Yamanouchi code, since every entry $d_i$ can be matched with a unique entry having value $d_{i}-1$ to its right, namely the next entry in the corresponding subsequence.  Note also that $d$ is Yamanouchi of shape $\mu$, since there are $\mu_1$ zeroes, $\mu_2$ ones, etc in $d$.  Since $c$ is bounded above component-wise by $d$, we have that $c$ is $\mu$-sub-Yamanouchi.
       
       We now show that $c$ is $A$-weakly increasing.  It suffices to show that for any two consecutive dots $w(t),w(t+1)$ of $w$ that are in inversion-friendly order, we have $c_t\le c_{t+1}$.  Suppose the dot $w(t)$ is labeled $i$ in the inversion labeling, and $w(t+1)$ is labled $j$.  Then by assumption, since they are in inversion-friendly order, we have either $i=j$ with the $j$ in a higher row than $i$, or $j<i$.  The $i$ is the left element of $c_t$ inversions and the $j$ is the left element of $c_{t+1}$ inversions.
       
       First suppose $i=j$ and the $j$ is in a higher row than the $i$, that is, $w(t+1)>w(t)$.  If $b$ is an index to the right of the $i$ such that $(w(t),w(b))$ is an inversion, then there are three possibilities: First, $w(b)$ could be labeled $i$ and be below $w(t)$, in which case $(w(t+1),w(b))$ is also an inversion.  Second, $w(b)$ could be labeled $i+1$ and be above $w(t)$ but below $w(t+1)$, in which case there is a dot labeled $i$ in row $w(b)$ to the right of $b$, forming an inversion with $w(t+1)$.  And third, $w(b)$ could be labeled $i+1$ and be above row $w(t+1)$, in which case $(w(t+1),w(b))$ is also an inversion.  Thus there is at least one inversion with $w(t)$ as its left element for every inversion with $w(t+1)$ as the left element, and so $c_{t}\le c_{t+1}$ in this case.
       
       Similarly, if $j<i$, then any dot labeled $i$ or $i+1$ has a dot labeled $j$ and a dot labeled $j+1$ to its right, and so $c_t\le c_{t+1}$ in this case as well.
       
       It follows that $\invcode$ is a well-defined map.  
      \end{proof}
      
      \begin{lemma}\label{injectiveinv}
        The map $\invcode$ is injective.
      \end{lemma}
      
      \begin{proof}[Proof of Theorem \ref{Invcode}]
      We will show that given a code $c$, we can form an inversion-friendly diagram by placing dots above $c_1,c_2,\ldots,c_n$ from left to right.  We claim that there is a unique height that is compatible with $c$ at each step.
      
      With the empty word as a trivial base case, we proceed inductively.  Suppose we have already placed the first $t-1$ dots from the left.  There may be several possible dot heights available for the $(t)$th dot, depending on the shape $\mu$ and which dot heights have already been chosen.  We claim that each possible height would result in a different value of the code number $c_t$.  To show this, let $h_1<h_2$ be two possible heights of the $(t)$th dot.  Since the first $t-1$ dots have been chosen and we know the shape of the diagram, the labels $i$ and $j$ of a dot at height $h_1$ or $h_2$ respectively are uniquely determined.  We also note that the inversion code number $c_t$ is uniquely determined by the choice of the $(t)$th dot (given the first $t-1$ dots), since any row of length $\mu_r\ge i$ that did not have a dot labeled $i$ among the first $t$ values must necessarily have one afterwards, and so the set of label values in each row to the right of the $(t)$th entry is determined.
      
      So, let $r$ be the inversion code number $c_t$ that would result from the dot at height $h_1$ labeled $i$, and $s$ the code number for $h_2$ labeled $j$.  We wish to show that $s\neq r$, and we consider the cases $j\le i$ and $j> i$ separately.
      
      If $j\le i$, let $k$ be the number of dots labeled $i$ that would be below and to the right of the $w(t)$ if $w(t)=h_1$ (labeled $i$).  Then $r-k$ would be the number of $i+1$'s above and to the right of it.  Each of the $k$ rows having the $i$'s also have $j$'s weakly to the right of them because $j\le i$, and each of the $r-k$ rows with the $i+1$'s have both a $j+1$ and a $j$ to the right.  Thus if $w(t)=h_2$ (labeled $j$) instead, the $j$ would have at least $r$ inversions, and so $s\ge r$.  But if $w(t)=h_2$, then this $j$ also forms an inversion with the $j$ in row $h_1$, giving an extra inversion.  Thus $s>r$, and so $s\neq r$ in this case.
      
      If $j>i$, consider the $s$ dots labeled $j$ or $j+1$ that would form an inversion with $w(t)$ if $w(t)=h_2$.  Then each of these rows would also contain an $i$ or $i+1$ that would form an inversion with the $i$ at height $h_1$, in addition to the row $h_2$ itself, showing that $r>s$.  Thus $s\neq r$, as desired.
      \end{proof}
      
      We have that $|C_{\mu,A}|=|\fancy{C}_A(\mu)|$ by our definition of $\fancy{C}$.  Furthermore, when $A=\{1,2,\ldots,n\}$ we have $|\ID_{\mu,A}|=\binom{n}{\mu}$ because we are simply counting the number of unrestricted diagrams having $\mu_1$ dots in the first row, $\mu_2$ in the second row, and so on.  We can now conclude bijectivity in this case.
       
      \begin{StandardInvCorollary}
        The map $\invcode$ is bijective in the case $A=\{1,2,\ldots,n\}$.
      \end{StandardInvCorollary}
       
       We are now ready to prove Theorem \ref{Invcode}.
       
       \begin{proof}
         We already have shown (Corollary \ref{standardinv}) that $\invcode$ is a bijective map $\ID_{\mu,[n]}\to C_{\mu,[n]}$.  Notice that for any other alphabet $A=\{a_1,\ldots,a_n\}$, we have $\ID_{\mu,A}\subset \ID_{\mu,[n]}$ and $C_{\mu,A}\subset C_{\mu,[n]}$.  We also know that the map $$\invcode:\ID_{\mu,[n]}\to C_{\mu,[n]}$$ restricts to an injective map $\invcode:\ID_{\mu,A}\to C_{\mu,A}$ by Lemmas \ref{welldefinedinv} and \ref{injectiveinv}.  It remains to show that it is surjective onto $C_{\mu,A}$.
         
         Let $c\in C_{\mu,A}\subset C_{\mu,[n]}$.  Then $c$ is $A$-weakly increasing on constant letters of $A$.  Let $d=\invcode^{-1}(c)\in \ID_{\mu,[n]}$.  We wish to show that $d$ is of inversion word type with respect to $A$, so that $d\in \ID_{\mu,A}$, that is, if $r<s$ and $a_r=a_s$ in $A$ then $(d(a_r),d(a_s))$ is not an inversion.  Suppose $(d(a_r),d(a_s))$ is an inversion.  Then either $d(a_r)$ and $d(a_s)$ are both dots labeled $i$ with $d(a_s)<d(a_r)$, or $d(a_r)$ is labeled $i$ and $d(a_s)$ labeled $i+1$ with $d(a_s)>d(a_r)$.
         
         In the first case, if $(d(a_s),d(a_t))$ is another inversion involving $a_s$, then either $d(a_t)$ is lower than $d(a_s)$ (and hence lower than $d(a_r)$) and labeled $i$, or it is above it and labeled $i+1$.  If the former then $(d(a_r),d(a_t))$ is an inversion, and if the latter, either there is an $i$ in the same row forming an inversion with $d(a_r)$, or the $i+1$ is above $d(a_r)$, forming an inversion with it.  Thus $d(a_r)$ is the left element of at least as many inversions as $d(a_s)$, plus one for the inversion $(d(a_r),d(a_s))$.  Thus $c_r>c_s$.
         
         In the second case, if $(d(a_s),d(a_t))$ is another inversion, then $d(a_t)$ is either lower (but possibly above $d(a_r)$) and labeled $i+1$, or higher and labeled $i+2$.  In the former case either $d(a_t)$ itself forms an inversion with $d(a_r)$ or the $i$ in its row does.  In the latter case the $i+1$ in its row forms an inversion with $d(a_r)$.  Since $(d(a_r),d(a_s))$ is an inversion as well, we again have $c_r>c_s$.  But this contradicts the fact that $c$ is $A$-weakly increasing.
         
         Hence $\invcode$ is surjective, and thus bijective, from $\ID_{\mu,A}$ to $C_{\mu,A}$.  Clearly the map preserves the statistics: the sum of all the entries of the inversion code of a diagram is the total number of inversions of the diagram, so $\invcode$ sends $\inv$ to $\Sigma$.  Therefore, $$\invcode:\ID_{\mu,A}\to C_{\mu,A}$$ is an isomorphism of weighted sets.
       \end{proof}

    \subsection*{Proof of Proposition \ref{OneColumnStandardize}: Standardizing Columns}
    
    \begin{OneColumnStandardize}
     For any column filling $\sigma$ with alphabet $A$, let $\rho=\Standardize(\sigma)$.  Then $\rho$ and $\sigma$ have the same major index, and $\majcode(\rho)$ is $A$-weakly increasing.
    \end{OneColumnStandardize}
    
    We first prove the following technical lemma.  Define a \textit{consecutive block} of $n$'s in a filling to be a maximal consecutive run of entries in a column which are all filled with the letter $n$.
         
    \begin{lemmaDuplicates}
      Given a filling of a one-column shape $\mu=(1^r)$ having largest entry $n$, there is a unique way of ordering the $n$'s in the filling, say $n_1,\ldots,n_{\alpha_n}$, such that the following two conditions are satisfied.
      \begin{enumerate}
         \item Any consecutive block of $n$'s in the column appears in the sequence in order from bottom to top, and
         \item If we remove $n_1,\ldots,n_{\alpha_n}$ in that order, and let $d_i$ be the amount that the major index of the column decreases at the $i$th step, then the sequence $d_1,d_2,\ldots,d_{\alpha_n}$ is weakly increasing.
      \end{enumerate}
    \end{lemmaDuplicates}
    
    \begin{proof}
     We first show (1) that there is a \textit{unique} choice of entry labeled $n$ at each step which minimizes $d$ and is at the bottom of a consecutive block, and then that (2) the resulting sequence $d_i$ is weakly increasing.  For any entry $x$, we define $\psi_x(\sigma)$ to be the column formed by removing the entry $x$ from $\sigma$.
         
     To prove (1), consider the bottommost entries of each consecutive block of $n$'s.  We wish to show that no two of these $n$'s have the same value of $d=\maj(\sigma)-\maj(\psi_n(\sigma))$ upon removal.  So, suppose there is an $n$ in the $i$th square from the top and an $n$ in the $j$th square from the top, each at the bottom of their blocks, and call them $n_i$ and $n_j$ to distinguish them.  Assume for contradiction that removing either of the $n$'s results in a decrease by $d$ of the major index.
         
     Suppose an entry $n$ has an entry $a$ above it and $b$ below.  In $\psi_n(\sigma)$, $a$ and $b$ are adjacent, and they can either form a descent or not.  If they do, then $d=\maj(\sigma)-\maj(\psi(\sigma))$ is equal to the number of descents below and including that $n$, and if they do not, then $d$ is equal to the sum of the number of descents \textit{strictly} below the $n$ \textit{plus} the position of the $n$ from the top.  We consider several cases based on the two possibilities for each of $n_i$ and $n_j$.
     
     If either $n_i$ or $n_j$ is at the very bottom of the filling, then removing that entry results in $d=0$, and the other does not, so we may assume neither of $n_i$ or $n_j$ is in the bottom row.  
     
     \textit{Case 1:} Each of $n_i$ and $n_j$ forms a new descent upon removal, in $\psi_{n_i}(\sigma)$ and $\psi_{n_j}(\sigma)$.  Assume without loss of generality that $i<j$, and let $t$ be the number of descents weakly below position $j$ (meaning its position from the top is greater than or equal to $j$) and $s$ the number of descents weakly below position $i$.  Then since the $n_i$ is at the bottom of its block, it is a descent, so $s>t$.  Since $s$ and $t$ are the values of $d$ for the removal of the two $n$'s, we have a contradiction.
    
     \textit{Case 2:} Neither $n_i$ nor $n_j$, upon removal, forms a new descent.  In this case, assume without loss of generality that $i<j$ and let $t$ be the number of descents \textit{strictly} below position $j$.  Let $r$ be the number of descents strictly between rows $i$ and $j$.  Since the $n$'s are at the bottom of their blocks, the two $n$'s are descents as well, so the values of $d$ upon removing the $n$'s are $i+r+t+1$ and $j+t$.  By our assumption, these are equal, and so we have 
     \begin{eqnarray*}
     i+r+1+t &=& j+t \\
     j-i-1 &=& r
     \end{eqnarray*}
     But $j-i-1$ counts the number of squares strictly between positions $i$ and $j$.  Since $r$ is the number of squares in this set which are descents, this means that every square between $i$ and $j$ must be a descent.  But the square in position $j$ has the highest possible label $n$, so the square just before it (above it) cannot be a descent.  Hence we have a contradiction.
     
     \textit{Case 3:} One of the two $n$'s, say the one in position $i$, forms a new descent upon removal, and the other does not.  Then in this case defining $t$ as the number of descents \textit{strictly} below position $j$ and $s$ the number of descents \textit{weakly} below position $i$, the two values of $d$ are $j+t$ and $s$.  So $j+t=s$ by our assumption, and so $j=s-t$, which implies $s-t>0$, or $s>t$.  Thus, necessarily $i<j$. 
         
     Now, $s-t$ is the number of descents between positions $i$ and $j$, inclusive.  Since $i\ge 1$ there are at most $j$ such squares, and the one preceding $j$ cannot be a descent since there is an $n$ in the $j$th position.   Thus this quantity $s-t$ is strictly less than $j$, but we showed before that $j=s-t$, a contradiction.  This completes the proof of claim (1).
         
     For claim (2), consider any two consecutive $d$ values in this process, say $d_1$ and $d_2$ for simplicity, that correspond to the largest value $n$.  Let $n_1$ and $n_2$ be the corresponding copies of $n$.  We wish to show that $d_1\le d_2$.
         
     First, notice that if $n_1$ and $n_2$ were in the same consecutive block before removal, we have $d_1=d_2$ unless $n_2$ is a block of length $1$ in $\psi(\sigma)$, in which case $d_2\ge d_1$.
         
     So we may assume that $n_1$ and $n_2$ were in different consecutive blocks before removal.  In this case the removal of $n_1$ may only change the value of $d$ on removing $n_2$ by at most one, namely by either shifting it back by one position if $n_1$ is above $n_2$ in the column, or by removing one descent from below $n_2$, if $n_1$ is below $n_2$.  Thus $d_2=\maj(\psi_{n_1}(\sigma))-\maj(\psi_{n_2}(\psi_{n_1}(\sigma)))$ is at most one less than $\maj(\sigma)-\maj(\psi_{n_2}(\sigma))$.  Since $n_1$ was chosen so as to minimize $d_1$, and we showed in our proof of (1) that the choice is unique, this implies that $d_2+1>d_1$.  Thus $d_2\ge d_1$, as desired.
         
     This completes the proof of (2).
    \end{proof}
    
    Proposition \ref{OneColumnStandardize} now follows from the proof of the above lemma.

    \subsection{Proof of Main Lemma: Proposition \ref{ZeroBump}}\label{Proofs1}

    \begin{MainLemma}[Main Lemma]
      Suppose $\sigma:\mu\to \ZZ_+$ is a filling for which $\inv(\sigma)=0$ and the largest entry $n$ appears in the bottom row.  Let $\sigma_\downarrow:\mu^{(1)}\to \ZZ_{+}$ be the filling obtained by:
      \begin{enumerate}
        \item Removing the rightmost $n$ from the bottom row of $\sigma$, which must be in the rightmost column since $\inv(\sigma)=0$,
        \item Shifting each of the remaining entries in the rightmost column down one row,
        \item Rearranging the entries in each row in the unique way so that $\inv(\sigma_\downarrow)=0$.
      \end{enumerate}
     Then the major index does not change: $$\maj(\sigma)=\maj(\sigma_\downarrow).$$
    \end{MainLemma}

    To prove Proposition \ref{ZeroBump}, we require a new definition and several technical lemmata.  We write $(i,j)$ to denote the square in row $i$ and column $j$ of a Young diagram.
        
    \begin{definition}
    The \textit{cocharge contribution} $\cc_{(i,j)}(\sigma)$ of an entry $\sigma(i,j)$ of a filling $\sigma$ is the number of descents that occur weakly below the entry $(i,j)$ in its column, $j$.
    \end{definition}
    
    It is easy to see that the cocharge contributions add up to the major index.
        
    \begin{proposition}\label{CochargeContributions}
     Let $\sigma:\mu\to \ZZ_+$ be any filling.  Then $\maj(\sigma)$ is equal to the sum of the cocharge contributions of the entries of $\sigma$, i.e. $$\maj(\sigma)=\sum_{(i,j)\in \mu}\cc_{(i,j)}(\sigma).$$
     \end{proposition}
         
\begin{figure} 
     \begin{center}
     \includegraphics{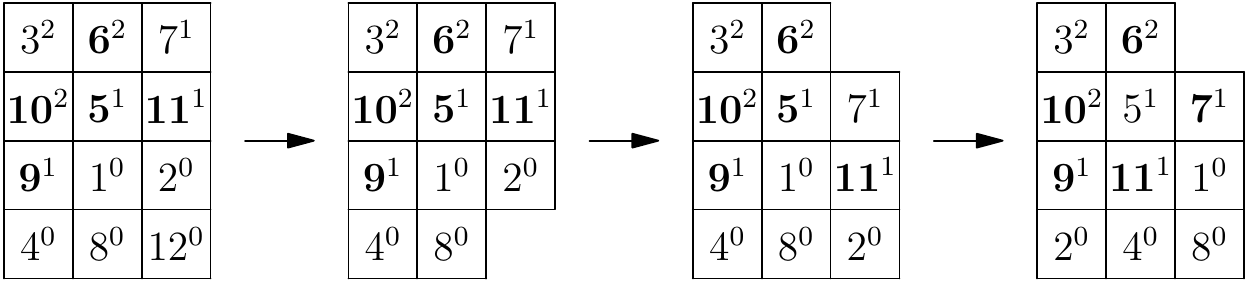}
     \end{center}
     \caption{\label{ContributionFigure} The cocharge contribution of the entries in each tableau is shown as a superscript.  Notice that the sum of the cocharge contributions of a tableau is equal to its major index.  In addition, the three-step process of Proposition \ref{ZeroBump} does not change the major index.}
     \end{figure}

      We omit the proof, and refer the reader to the example in Figure \ref{ContributionFigure}.
         
     \begin{definition} 
     Let $w$ be any sequence consisting of $k$ $0$'s and $k$ $1$'s, and let $a_1,a_2,\ldots,a_k$ be any ordering of the $0$'s.  We define the \textit{crossing number} of $w$ with respect to this ordering as follows.  Starting with $a_1$, let $b_1$ be the first $1$ to the right of $a_1$ in the sequence, possibly wrapping around cyclically if there are no $1$'s to the right of $a_1$.  Then let $b_2$ be the first $1$ cyclically to the right of $a_2$ other than $b_1$, and so on.  Then the \textit{crossing number} is the number of indices $i$ for which $b_i$ is to the left of $a_i$.
        \end{definition}
        
    \begin{example}
      If we order the $0$'s from left to right, the word $10110010$ has crossing number $2$.
    \end{example}
        
    \begin{lemma}\label{aabb}
      Let $w$ be any sequence consisting of $k$ $0$'s and $k$ $1$'s.  Then its crossing number is independent of the choice of ordering of the $0$'s.
    \end{lemma}
      
    \begin{proof}
       Say that a word is \textit{$0$-dominated} if every prefix has at least as many $0$'s as $1$'s.  First, we note that there exists a cyclic shift of $w$ which is $0$-dominated.  Indeed, consider the partial sums of the $(-1)^{w_i}$'s in the sequence, so that any $0$ contributes $+1$ and any $1$ contributes $-1$.  The total sum is $0$, and we can shift to start at the index of the minimal partial sum; the partial sums will now all be positive.
          
       Now, we show by induction that any $0$-dominated sequence has crossing number $m=0$.  It is clearly true for $k=1$, since the only $0$-dominated sequence is $01$ in this case.  
       
       Suppose the claim holds for any $0$-dominated sequence of $k-1$ $0$'s and $k-1$ $1$'s and let $s$ be an $0$-dominated sequence with $k$ $0$'s.  Choose an arbitrary $0$ to be $a_1$, and denote it $\hat{0}$.  Then since $s$ is $0$-dominated, the last term in $s$ is a $1$ and so $\hat{0}$ will be paired with a $1$, denoted $\hat{1}$, to the right of it.  Remove both $\hat{0}$ and $\hat{1}$ from $s$ to form a sequence $s'$ having $k-1$ $0$'s and $k-1$ $1$'s.
          
       We claim that $s'$ is $0$-dominated.  Note that all prefixes of $s'$ that end to the left of $\hat{0}$ are unchanged, and hence still have at least as many $0$'s as $1$'s.   Any prefix $P'$ that ends between $\hat{0}$ and $\hat{1}$ is the result of removing $\hat{0}$ from a corresponding prefix $P$ of $s$, which had at least as many $0$'s as $1$'s.  If there were an equal number of $0$'s as $1$'s in $P$, then its last term is a $1$.  This means that $\hat{1}$ was not the first $1$ to the right of the $0$, a contradiction.  So $P$ has strictly more $0$'s than $1$'s, and so $P'=P\setminus\{\hat{0}\}$ has at least as many $0$'s as $1$'s.  Finally, any prefix which ends to the right of $\hat{1}$ has one less $0$ and one less $1$ than the corresponding initial subsequence of $s$, and so it also has at least as many $0$'s as $1$'s.  It follows that $s'$ is $0$-dominated.
          
       By the inductive hypothesis, no matter how we order the remaining $0$'s, there are no crossing pairs.  Since the choice of $a_1$ was arbitrary, the crossing number is $0$ for any ordering of the $0$'s.
          
       Returning to the main proof, let $w=w_1w_2\cdots w_{2k}$ and let $i$ be such that the cyclic shift $w'=w_iw_{i+1}\cdots w_{2k}w_1w_2\cdots w_{i-1}$ is $0$-dominated.  Then every pairing in $w'$ has the $0$ to the left of the $1$, and so the crossing number of $w$ is the number of pairings in which the $0$ is among $w_i\cdots w_{2k}$ and the $1$ is among $w_1\cdots w_{i-1}$.   Hence, the crossing number is equal to the difference between the number of $1$'s and $0$'s among $w_1w_2\cdots w_{i-1}$.  This is independent of the choice of order of the $0$'s, and the proof is complete.
        \end{proof}
        
      In the rest of the paper, if a row $r$ is above a row $s$ in a filling, we say that we \textit{rearrange $r$ with respect to $s$} if we place the entries of $r$ in the unique ordering for which there are no inversions in row $r$, given that $s$ is below it.
        
     \begin{lemma}\label{twist}
        Let $\sigma$ be a filling of the two-row shape $(k,k)$ with $\inv(\sigma)=0$.  Let $\sigma'_\pi$ be formed by rearranging the bottom row via the permutation $\pi$, nad rearranging the top row with respect to the new bottom row.  Then $\maj(\sigma)=\maj(\sigma'_\pi)$.
        \end{lemma}
        
      \begin{proof}
       Let $w$ be cocharge word of the diagram.  No matter what the permutation of rows, the cocharge word will remain unchanged, a sequence of $k$ $1$'s and $k$ $2$'s.  But the permutation of the bottom row determines a permutation of the $1$'s, and the subsequent ordering of the top row is determined by the process of selecting the first remaining $2$ cyclically to the right of the $1$ at each step.  It forms a descent if and only if that $2$ is to the left of the $1$, i.e. if it contributes to the crossing number.  So the number of descents is equal to the crossing number of the cocharge word (thinking of the $1$'s as $0$'s and the $2$'s as $1$'s), and by Lemma \ref{aabb} the proof is complete.
      \end{proof}
         
       We now have the tools to prove the next technical lemma.
         
     \begin{lemma}\label{dangle}
       Let $a_1,\ldots,a_{w-1}$ be any positive integers, and suppose $b_1,\ldots,b_w$ are positive integers such that in the partial tableau
        $$\begin{array}{ccccc}
             b_1 & b_2 & \cdots & b_{w-1} & b_w \\
             a_1 & a_2 & \cdots & a_{w-1} &
        \end{array}$$
       has no inversions among the $b_i$'s.  Then if we rearrange $a_1,\ldots,a_{w-1}$ in any way and then rearrange the $b$'s in the unique way that guarantees no inversions among the $b$'s, then the entry $b_w$ is still in the last position.  Furthermore, the total number of descents among $b_1,\ldots,b_w$ is unchanged after this operation.
     \end{lemma}
         
         \begin{proof}
           Consider the cyclic ordering of $a_1,\ldots,a_{w-1},b_1,\ldots,b_{w}$.  Since there are no inversions among the $b$'s, we have that $a_i,b_i,b_w$ are in cyclic order for each $i$, possibly with $b_i=b_w$ or $a_i=b_w$.
           
           Let $b_w,t_1,\ldots,t_{2w-2}$ be the ordering of these letters that is in cyclic order, with ties broken in such a way that $b_w,a_i,b_i$ occur in that order in the sequence for each $i$.  Then if we replace the $a_i$'s with $0$'s and the $b_i$'s with $1$'s, the suffix $t_1,\ldots,t_{2w-2}$ has crossing number $0$ since each $a_i$ is paired with $b_i$ to its right.
           
          It follows from Lemma \ref{aabb} that, if we rearrange the $a_i$'s, the crossing number is still $0$ and so $b_w$ still corresponds to the $1$ at the beginning of the sequence.  It follows that $b_w$ is still in the last position in the new filling.  Finally, by considering only the first $w-1$ columns, we can apply Lemma \ref{twist} to see that the total number of descents among $b_1,\ldots,b_{w-1}$ remains unchanged.
         \end{proof}
         
         We require one more technical lemma regarding two-row fillings.  First, notice that in a two-row shape with the bottom row ordered least to greatest and no inversions in the second row, the descents must be ``left-justified": they must occur in columns $1,\ldots,k$ for some $k$.  For, if $b_r>a_r$ is a descent and $b_{r-1} \le a_{r-1}$ is not, then $b_r>a_{r-1}$ by transitivity and we have $b_{r-1}\le a_{r-1}<b_{r}$, forming an inversion.  Moreover, after the descents the $b_i$'s are weakly increasing: $b_i\le b_j$ for $k<i<j$ - this follows directly from the fact that none of these $b_i$'s are descents.  The descents $b_1,\ldots,b_k$ are also weakly increasing; otherwise we would have an inversion.
         
         We will use these facts repeatedly throughout.
         
         \begin{lemma}\label{TwoRowDrop}
           Let $a_1\le\cdots\le a_{w-1}$ and let $b_1,b_2,\ldots,b_{w}$ be numbers such that the partial tableau
           $$\begin{array}{ccccc}
             b_1 & b_2 & \cdots & b_{w-1} & b_w \\
             a_1 & a_2 & \cdots & a_{w-1} &
           \end{array}$$
           has no inversions in the second row.  Then if we bump $b_w$ down one row so that $$a_1\le a_2\le \cdots\le a_t\le b_w< a_{t+1}\le\cdots\le a_{w-1}$$ is the bottom row, and leave $b_1,\ldots,b_{w-1}$ unchanged, then the new tableau still has no inversions, and the descents in the second row remain the same (and left-justified).
         \end{lemma}
         
         \begin{proof}
         
           Let $k$ be the number of descents among the $b$'s.  If $k=0$, there are no descents, and we must have $b_w\le a_1$ so as not to have inversions.  In this case, $b_w$ drops down into the first position in the bottom row, and there are still no descents and no inversions since $b_1\le b_2\le \ldots\le b_{w-1}$ in this case.
           
           If $k\ge 1$, then $b_k>a_k$ is the last descent.  Since $b_k$ and $b_w$ do not form an inversion in the original tableau, we must either have $a_k<b_k\le b_w$ or $b_w\le a_k<b_k$.   We consider these cases separately.
           
          \textit{Case 1:} Suppose $a_k<b_k\le b_w$.  Then $t>k$, i.e. $b_w$ drops to a position to the right of the last descent, after which point we have $b_i\le a_i$ for all such $i$.  Thus, for instance, $b_{t+1}<a_{t+1}$, and since $b_w$ and $b_{t+1}$ did not originally form a descent, we must have $b_{t+1}\le b_w\le a_{t+1}$.  This means that $b_{t+1}\le b_w$, so $b_{t+1}$ still does not form a descent in the new tableau.  Then, similarly we have $b_{t+2}\le b_w$, and so $b_{t+2}\le a_{t+1}$, and so on.  Thus the descents have stayed the same in the new tableau.  
          
          Furthermore, since $b_i<b_w$ for all $i\ge t+1$ in this case, we have $b_{i}<b_w<a_i$ for all $i\ge t+1$, and since the $b_i$'s after position $k$ are weakly increasing, none of these form inversions.  Since $b_1,\ldots,b_t$ are above the same letters $a_1,\ldots,a_t$ as before and are in the same positions relative to the other $b_i$'s, they cannot be the left elements of inversions either.
           
          \textit{Case 2:} Suppose now that $b_w\le a_k<b_k$.  If $b_w=a_k$ then in fact it drops to the right of $a_k$ and it is the same as the previous case.  So we can assume that $b_w<a_k<b_k$. 
          
          Then $t\le k$, i.e. $b_w$ drops to a position underneath a descent of the original tableau shape.  Since $b_w\le a_{t+1}$ and $a_{t+1}<b_{t+1}$ is a descent, we have $b_{w}<b_{t+1}$ and so $b_{t+1}$ is still a descent in the new tableau.  Similarly $b_i$ is still a descent for all $i\le k$.  To check that $b_{k+1}$ is still \textit{not} a descent, assume it is: that $a_k<b_{k+1}$.  Then $b_w\le a_k<b_{k+1}$, and so $b_w\le a_{k+1}\le b_{k+1}$ since the original filling had no inversions.  If $a_{k+1}<b_{k+1}$, we get a contradiction, so $a_{k+1}=b_{k+1}$.  But then $b_w=a_{k+1}$, contradicting the fact that $b_w<a_{k+1}$.  Thus there is not an inversion in the $(k+1)$st position.  Hence the descents stay the same in this case as well.
           
           Furthermore, consider $b_i$ and $b_j$ with $i<j<w$: if $i$ is among $1,\ldots,t$ then $b_i$ and $b_j$ do not form an inversion since $b_i$ is still above $a_i$.  If $i$ and $j$ are both among $t+1,\ldots,k$, then they do not form an inversion, since $b_i$ and $b_j$ are both descents and $b_i<b_j$.  If $i$ is among $t+1,\ldots,k$ and $j>k$, note that $b_j<b_w$ since it is in the run of non-descents of the $b$'s.  Hence $b_j<a_{i}$ by transitivity, and so $b_j<a_i<b_i$ since $b_i$ is a descent.  This implies that $b_i$ and $b_j$ do not form an inversion.  Finally, if $i>k$ and $j>i$, we are once again in the run of non-descents at the end, which is weakly increasing, and hence there are no inversions since none are descents.  We conclude that the $b_i$'s have no inversions among them in this case either.
         \end{proof}
         
        \begin{lemma}\label{OnTop}
         Let $a_1,\ldots,a_{w-1}$, $b_1,\ldots,b_{w}$, and $c_w$ be numbers such that the partial filling 
         $$\begin{array}{ccccc}
                    &     &        &         & c_w \\
                b_1 & b_2 & \cdots & b_{w-1} & b_w \\
                a_1 & a_2 & \cdots & a_{w-1} &
            \end{array}$$
          has no inversions in the second row.  Then there exists an ordering $t_1,\ldots,t_w$ of $a_1,\ldots,a_{w-1},b_w$ such that if $s_1,\ldots,s_w$ is the unique ordering of $b_1,\ldots,b_{w-1},c_w$ for which the partial filling 
                $$\begin{array}{cccc}
                        s_1 & s_2 & \cdots & s_{w}  \\
                       t_1 & t_2 & \cdots & t_{w} 
                  \end{array}$$
           has no inversions in the second row, then the entry $c_w$ is directly above $b_w$ in the new filling.
         \end{lemma}
                 
         \begin{proof}
          Let $T'$ be the two-row filling consisting of the $s$'s and $t$'s as in the statement of the lemma.  Let $x$ be the cocharge word of $T'$, with the bottom row indexed by $0$ and the top by $1$.  Then $x$ consists of $0$'s and $1$'s, and as in Lemma \ref{aabb}, the number of descents in $T'$ is the crossing number of this word.  So $b_w$ is one of the $0$'s in this word, and $c_w$ is one of the $1$'s, and we wish to show that there is some ordering of the $0$'s in which $b_w$ is paired with $c_w$.
          
          Assume to the contrary that $b_w$ cannot be paired with $c_w$ no matter how we order the $0$'s.  Choose a cyclic shift $\widetilde{x}$ of $x$ whose crossing number is $0$, as we did in Lemma \ref{aabb}.  If $b_w$ is to the left of $c_w$ in $\widetilde{x}$, then since it can't be paired with $c_w$ there must be an index $k$ between that of $b_w$ and $c_w$ at which the prefix of the first $k$ letters is $0$-dominated.  For, if there were more $0$'s than $1$'s at every step up to $c_w$ then we can pair off the other $0$'s starting from the left until $c_w$ is the first $1$ to the right of $b_w$.  This means we can choose a different cyclic ordering, starting at the $k+1$st letter, for which the crossing number is also $0$.  In this cyclic shift, $c_w$ is to the left of $b_w$.  So we have reduced to the case that $c_w$ is to the left of $b_w$.
          
          In this case, $c_w$ is one of the $1$'s, and $b_w$ is one of the $0$'s, e.g. in the $0$-dominated sequence $001011$, we might have $c_w$ be the third entry and $b_w$ the fourth.  Before we dropped down the $b_w$ and $c_w$, we had a tableau whose cocharge word looked like this word except with the $0$ of $b_w$ replaced by a $1$, and the $1$ of $c_w$ replaced by a $2$ (in the example, this would give us the word $002111$.)  Remove the $2$ from this word.  In the resulting word of $0$'s and $1$'s, since we have bumped up a $0$ to a $1$ but removed one of the $1$'s before it, every prefix is $0$-dominated except the entire word, which has one more $1$ than it has $0$'s.  Thus the very last $1$ is the only entry which is not paired.  But $b_w$ is, by assumption, the entry which is unpaired in the original ordering.  This is a contradiction, since $b_w$ was a $0$ in the bumped-down word and hence could not have been in the last position.
          
          It follows that there must exist an ordering of the $0$'s in which $b_w$ is paired with $c_w$.  This completes the proof.
         \end{proof}
         
         In the next two lemmas, we let $\sigma:\mu\to \ZZ_+$ be a filling with $\inv(\sigma)=0$ whose largest entry appears in the bottom row, and let $\sigma_\downarrow:\mu^{(1)}\to \ZZ_+$ be constructed from $\sigma$ as in the statement of Proposition \ref{ZeroBump}.
         
        \begin{lemma}\label{CochargeConservation1}
          Suppose $\inv(\sigma)=0$.  Let $i\ge 1$ be an index such that $\mu_{i+1}=\mu_1$, i.e. the $(i+1)$st row of $\mu$ is as long as the bottom row.  Then we have 
          $$\cc_{(i+1,\mu_1)}(\sigma)+\sum_{1\le j\le \mu_1-1} \cc_{(i,j)}(\sigma)=\sum_{1\le j\le \mu_1} \cc_{(i,j)}(\sigma_\downarrow).$$
        \end{lemma}
        
        \begin{proof}
          We induct on $i$.  For the base case, $i=1$, the left hand side is the total cocharge contribution of the entries $(1,1),(1,2),\ldots,(1,\mu_1-1)$ and the entry $(2,\mu_1)$.  The square $(1,\mu_1)$ is filled with the largest number $n$, by our assumption that $n$ appears in the bottom row and the fact that $\inv(\sigma)=0$.  Thus the entry in $(2,\mu_1)$ cannot be a descent, and so the cocharge contribution of all of these entries are $0$.  Thus the left hand side is $0$. The right hand side is also $0$, since it is the sum of the cocharge contributions from the bottom row of $\sigma_\downarrow$.  
          
          For the induction, let $i>1$ and suppose the claim is true for $i-1$.  Then the induction hypothesis states that 
          $$s:=\cc_{(i,\mu_1)}(\sigma)+\sum_{1\le j\le \mu_1-1} \cc_{(i-1,j)}(\sigma)=\sum_{1\le j\le \mu_1} \cc_{(i-1,j)}(\sigma_\downarrow).$$  Then if there are $k$ descents among the entries $(i+1,\mu_1)$ and $(i,1),\ldots,(i,\mu_{1}-1)$ of $\sigma$, then their total cocharge contribution is equal to $s+k$, since they are the entries strictly above those that contribute to the left hand side of the equation above.
          
          So, to show that $$\cc_{(i+1,\mu_1)}(\sigma)+\sum_{1\le j\le \mu_1-1} \cc_{(i,j)}(\sigma)=\sum_{1\le j\le \mu_1} \cc_{(i,j)}(\sigma_\downarrow),$$ it suffices to show that the total cocharge contribution of the $i$th row of $\sigma_\downarrow$ is also $s+k$.  By the induction hypothesis it is equivalent to show that there are $k$ descents among the entries in the $i$th row of $\sigma_\downarrow$.
          
          Now, let $w=\mu_1$ be the width of the tableau, and let $a_1,\ldots,a_{w-1}$ be the first $w-1$ entries in row $i-1$ of $\sigma$.  Let $b_1,\ldots,b_w$ be the elements of row $i$, and let $c_w$ be the entry in square $(i+1,w)$, above $b_w$. 
                  $$\begin{array}{ccccc}
                           &     &        &         & c_w \\
                       b_1 & b_2 & \cdots & b_{w-1} & b_w \\
                       a_1 & a_2 & \cdots & a_{w-1} &
                  \end{array}$$
          
          Consider the $2\times w$ tableau $T'$ with bottom row elements $a_1,\ldots,a_{w-1},b_w$ and top row elements $b_1,\ldots,b_{w-1},c_w$.   By Lemma \ref{OnTop}, there is a way of rearranging the bottom row of $T'$ such that if we rearrange the top row respectively, then $c_w$ lies above $b_w$.  This suffices, for now the remaining columns will form a tableau with no inversions in the second row, with $a_1,\ldots,a_{w-1}$ and $b_1,\ldots,b_{w-1}$ as the entries of the rows.  By Lemma \ref{twist} this has the same number of descents independent of the ordering of the $a_i$'s, and $c_w$ will be a descent or not depending on whether it was a descent before.  Thus there are still $k$ descents in the $i$th row.
        \end{proof}
        
        Lemma \ref{CochargeConservation1} shows that the cocharge contribution is conserved for rows $i$ for which $\mu_{i+1}=\mu_{1}$.  The next lemma will show that the cocharge contribution is unchanged for higher rows as well.  Again, here $\sigma$ is a filling having its largest entry $n$ occurring in the bottom row.
        
        \begin{lemma}\label{CochargeConservation2}
          Suppose $\inv(\sigma)=0$, and the rightmost ($w$th) column of $\mu$ has height $\mu_w^\ast=h$.  Then in $\sigma_\downarrow$, row $h$ consists of the first $w-1$ letters of row $h$ of $\sigma$ in the same order, and their cocharge contributions are the same as they were in $\sigma$.
        \end{lemma}
        
        It follows from this lemma that all higher rows are unchanged as well, and combining this with Lemma \ref{CochargeConservation1}, it will follow that $\maj(\sigma)=\maj(\sigma_\downarrow)$.
        
        \begin{proof}
          We induct on $h$, the height of the rightmost column.  For $h=1$ and $h=2$, we are done by previous lemmata (see Lemma \ref{TwoRowDrop}).  So, suppose $h\ge 3$ and the claim holds for all smaller $h$.
          
          Performing the operation of Proposition \ref{ZeroBump}, suppose we have bumped down all but the topmost entry (in row $h$) of the rightmost column and rearranged each row with respect to the previous.  Let rows $h-2$, $h-1$, and $h$ have contents:
         $$\begin{array}{ccccc}
            d_1 & d_2 & \cdots & d_{w-1} & d_w \\
              c_1 & c_2 & \cdots & c_{w-1} &   \\
              x_1 & x_2 & \cdots & x_{w-1} & x_w 
           \end{array}$$   
         Notice that, by the induction hypothesis, the entries $c_1,\ldots,c_{w-1}$ are the same as they were in $\sigma$ before bumping down $c_w$ and have the same cocharge contributions as they did before.  Thus the row of $d$'s as shown is currently the same as row $h$ of $\sigma$.  So, we wish to show that upon bumping $d_w$ down and rearranging all rows so that the filling has no inversions, the entries in row $h$ are still $d_1,d_2,\ldots,d_{w-1}$ in that order, and that these entries have the same cocharge contributions as they did before.
         
         We first show that the entries $d_1,\ldots,d_{w-1}$ do not change their positions upon bumping $d_w$ down to row $h-1$ (and rearranging so that there are still no inversions.)  We proceed by strong induction on the width $w$.  For the base case, $w=2$, we have that $d_1$ is the only entry left in the top row, and therefore cannot change its position.
         
         Now, assume that the claim is true for all widths less than $w$.  If $d_w$ bumps down and inserts in a row $t$ above $x_t$, then the numbers $c_1,\ldots,c_{t-1}$ are still above $x_1,\ldots,x_{t-1}$ respectively since they are still first in cyclic order after each.  Likewise the entries $d_1,\ldots,d_{t-1}$ remain the same in this case.  Thus we may delete the first $t-1$ columns and reduce to a smaller case, in which the claim holds by the induction hypothesis.  This allows us to assume that when $d_w$ bumps down, it is in the first column, above $x_1$, and so the tableau looks like:
              $$\begin{array}{ccccc}
                 d_\ast & d_\ast & \cdots & d_{\ast} & \\
                   d_w & c_\ast & \cdots & c_{\ast} & c_{\ast}   \\
                   x_1 & x_2 & \cdots & x_{w-1} & x_w 
                \end{array}$$  
                where the $\ast$'s are an appropriate permutation of the indices for $d_1,\ldots,d_{w-1}$ and $c_1,\ldots,c_{w-1}$.
         
         We now show that $d_1,\ldots,d_r$ remain in their respective positions for all $r\ge 1$, by induction on $r$.  (So, we are doing a triple induction on the height, the width of the tableau, and the index of the $d$'s).  For the base case, we wish to show that $d_1$ is the entry above $d_w$ in the new tableau.  We have, from the fact that $\inv(\sigma)=\inv(\sigma_\downarrow)=0$, that the following triples are in cyclic order for any $k$ such that $2<k<w$:
         \begin{enumerate}
           \item $(x_1,d_w,c_1)$, with possible equalities $x_1=c_1$, $d_w=c_1$
           \item $(x_1,c_1,c_k)$, with possible equalities $x_1=c_k$, $c_1=c_k$
           \item $(c_k,d_k,d_w)$, with possible equalities $c_k=d_w$, $d_k=d_w$
           \item $(c_1,d_1,d_k)$, with possible equalities $c_1=d_k$, $d_1=d_k$
         \end{enumerate}
         Combining (1) and (2) above, we have that $(x_1,d_w,c_1,c_k)$ are in cyclic order, and so in particular $(d_w,c_1,c_k)$ is in cyclic order.  Combining this with (3) above, we have $(d_k,d_w,c_1,c_k)$ are in cyclic order, and in particular so are $(d_k,d_w,c_1)$.  Using this and (4), we have $(d_k,d_w,c_1,d_1)$ are in cyclic order, and in particular either $c_1\neq d_1$ or $c_1=d_1=d_k$, and so $(d_w,d_1,d_k)$ are in cyclic order with either $d_w\neq d_1$ or $d_w=d_k=d_1$; this implies that $d_1$ and $d_k$ will not form an inversion if $d_1$ is placed above the $d_w$.  Thus $d_1$ does indeed stay in the leftmost column.
         
         For the induction step, suppose $d_1,\ldots, d_{r-1}$ are in columns $1,\ldots,r-1$ respectively in $\sigma_\downarrow$.  We wish to show that $d_r$ must be in the $r$th position.  To do so, first notice that since $c_i$ is first in cyclic order after $x_i$ among $c_i,c_{i+1},\ldots,c_{w-1}$ for each $i$, we have that for each $i$, the element that appears above $x_i$ after bumping $d_w$ down is among $c_1,\ldots,c_i$.  
         
         Suppose $c_k$ is above $x_k$ in the new tableau for some $k\le r-1$.  Then $d_k$ is in this column as well by the induction hypothesis, and so removing this entire column will not affect the relative ordering of the remaining entries.  But now $d_r$ is the $(r-1)$st of the $d$'s in question, and therefore must be in the $(r-1)$st position by the induction hypothesis, and so must be in the $r$th position in the full tableau (prior to removing the $k$th column).
         
         Otherwise, if $c_k$ is never above $x_k$ for any $k\le r-1$, we have that $c_1$ must appear above $x_2$, since it can only be $c_1$ or $c_2$ but is not $c_2$ by assumption.  Then, $c_2$ must appear above $x_3$, and so on, up to $c_{r-2}$ appearing above $x_{r-1}$.  If $c_r$ appears above $x_r$, then $d_r$ must be above that since we knew from the previous tableau that it is first in cyclic order after $c_r$ among $d_r,\ldots,d_{w-1}$.  So the only case that remains is where $c_{r-1}$ appears in column $r$, above $x_r$.  The diagram is as follows:
          $$\begin{array}{ccccccc}
                d_1 & d_2 & d_3 & \cdots & d_{r-1} & d_{\ast}   & \cdots \\
                d_w & c_1 & c_2 & \cdots & c_{r-2} & c_{r-1} & \cdots   \\
                x_1 & x_2 & x_3 & \cdots & x_{r-1} & x_r     & \cdots
            \end{array}$$  
          We wish to show that the entry $d_\ast$ above is $d_r$.  First, we claim that $(d_w,c_1,c_2,\ldots,c_{r})$ are in cyclic order.  For, we have $(x_1,c_1,c_2)$ and $(x_1,d_w,c_1)$ are in cyclic order, so $(x_1,d_w,c_1,c_2)$ are.  Since $(x_2,c_1,c_2)$ and $(x_2,c_2,c_3)$ are in cyclic order, we have that $(x_2,c_1,c_2,c_3)$ are in cyclic order.  Since $(x_1,c_1,c_3)$ are in cyclic order as well, we can combine this with the last two observations to deduce that $$(x_1,d_w,c_1,c_2,c_3)$$ are in cyclic order.  Now, we can use the triples $(x_3,c_3,c_4)$, $(x_3,c_2,c_3)$, and $(x_2,c_2,c_4)$ to deduce that $(x_2,c_1,c_2,c_3,c_4)$ are in cyclic order as well.  But since $(x_1,c_1,c_4)$ are in cyclic order, this means that $$(x_1,d_w,c_1,c_2,c_3,c_4)$$ are in cyclic order as well, and so on.  At each step, to add $c_k$ to the list we only need consider rows up to that of $x_{k-1}$.  Hence, the process continues up to $k=r$.
          
          Finally, notice that since we are only concerned with relative cyclic order of the entries to determine their positions, we may cyclically increase all the entries modulo the highest entry in such a way that $d_w\le c_1\le c_2\le \cdots \le c_{r}$ in actual size.  Furthermore, since we are currently only concerned with the position of $d_r$, which is determined by its relative ordering with $d_i$ for $i>r$ and with $c_{r-1}$, we may assume that $c_r\le c_{r+1}\le c_{r+2}\le \ldots\le c_{n}$ are increasing as well; it will make no difference as to the value of $d_\ast$.  But then the top two rows behave exactly as in the two-row case of Lemma \ref{TwoRowDrop}.  We know that $d_r$ occurs in the $r$th column from this lemma, and the induction is complete.
          
          We have shown that $d_1,\ldots,d_{w-1}$ retain their ordering, and it remains to show that they retain their cocharge contributions.  If any $c_k$ lies above $x_k$, and hence $d_k$ above it, the column has not changed and so $d_k$ does indeed retain its cocharge contribution.  So, as before, we may remove such columns and reduce to the case in which the entries are:
           $$\begin{array}{cccccc}
                d_1 & d_2 & d_3 & \cdots & d_{w-1} &     \\
                d_w & c_1 & c_2 & \cdots & c_{w-2} & c_{w-1}   \\
                x_1 & x_2 & x_3 & \cdots & x_{w-1} & x_w    
            \end{array}$$  
          For the first column, we have that $(x_1,d_w,c_1)$ are in cyclic order since $d_w$ and $c_1$ do not form an inversion.  Moreover, either $x_1\neq d_w$ or $x_1=d_w=c_1$, in which case we may assume that $d_w$ is in fact located in the second column instead, and reduce to a smaller case.  So we may assume $x_1\neq d_w$.  In addition, $(c_1,d_1,d_w)$ are in cyclic order, with $c_1\neq d_1$ unless $c_1=d_1=d_w$, and if $d_1=d_w$ then we must have $d_1=d_2=\cdots=d_w$ so that $d_1$ does not form an inversion with any element in the new tableau.  We now consider three cases based on the actual ordering of $x_1,d_w,c_1$ (which are in cyclic order):
          
          \textit{Case 1:}  Suppose $x_1<d_w\le c_1$.  Then since $(c_1,d_1,d_w)$ are in cyclic order, either $d_1$ is greater than both $c_1$ and $d_w$ or less than or equal to both.  Since both $c_1$ and $d_w$ are descents when over $x_1$, the cocharge contribution of $d_1$ is unchanged in this case.
          
          \textit{Case 2:}  Suppose $d_w\le c_1\le x_1$.  Then in this case neither $c_1$ nor $d_w$ is a descent when in the first column, and the same analysis as in Case 1 shows that $d_1$ has the same cocharge contribution in either case.
                
          \textit{Case 3:}  Suppose $c_1\le x_1<d_w$.  Then $c_1\le d_1\le d_w$.  If $d_1$ is strictly greater than $c_1$, it forms a descent with $c_1$ and not with $d_w$.  But note that $d_w$ is a descent when in the first column, and $c_1$ is not, so the total number of descents weakly beneath $d_1$ balances out and is equal in either case.  If $d_1=c_1$, then $d_1=d_2=\cdots=d_w$, which is impossible since then $c_1=d_w$.  So the cocharge contribution of $d_1$ is the same in this case as well.
          
          This completes the proof that $d_1$ retains the same cocharge contribution.  We now show the same holds for an arbitrary column $i$.
          
          In the $i$th column, we have $d_i$ above $c_{i-1}$ above $x_i$.  Note that $(c_i,d_i,d_w)$ and $(d_w,c_{i-1},c_i)$ are in cyclic order (the latter by the above argument which showed that $d_w,c_1,c_2,\ldots,c_{w-1}$ are in cyclic order given that the $c_i$'s are arranged as above), so $(c_i,d_i,d_w,c_{i-1})$ are in cyclic order.  In particular $(c_i,d_i,c_{i-1})$ are in cyclic order.   Moreover, if $c_i=d_i$ then $d_i=d_w=c_i$.  Since $d_w,c_{i-1},c_i$ are in cyclic order we must have $c_i=d_i=c_{i-1}$ in this situation.
          
          We also have that $(x_i,c_{i-1},c_i)$ are in cyclic order, and by a similar argument as above we can assume $x_i\neq c_{i-1}$.  So either $x_i< c_{i-1}\le c_{i}$, $c_i\le x_i<c_{i-1}$, or $c_{i-1}\le c_i\le x_i$.  The exact same casework as above for these three possibilities then shows that $d_i$ retains its cocharge contribution.
        \end{proof}
        
        Proposition \ref{ZeroBump} now follows immediately from Lemmas \ref{CochargeConservation1} and \ref{CochargeConservation2} and Proposition \ref{CochargeContributions}.
        
   \subsection{Proof of Theorem \ref{RectangleBump}: Reducing Rectangles}\label{Proofs2}
    
    We first recall the statement of the theorem.
    
    \begin{theorem}\label{RectangleBump}
            Let $A=\{1,2,\ldots,n\}$ be the alphabet with content $\alpha=(1^n)$, and let $\mu=(a,a,a,\ldots,a)$ be a rectangle shape of size $n$.  Then there is a weighted set isomorphism 
              $$\psi:(\fancy{F}_{\mu}^{(1^n)}|_{\inv=0};\maj)\to \bigsqcup_{d=0}^{\mu_1^\ast-1}(\fancy{F}_{\mu^{(d+1)}}^{(1^{n-1})}|_{\inv=0};\maj+d)$$
            defined combinatorially by the following process.  
           \begin{enumerate}
            \item Given a filling $\sigma:\mu\to \ZZ_+$ with distinct entries $1,\ldots,n$ and $\inv(\sigma)=0$, let $i$ be the row containing the entry $n$.  Split the filling just beneath row $i$ to get two fillings $\sigma_{top}$ and $\sigma_{bot}$ where $\sigma_{bot}$ consists of rows $1,\ldots,i-1$ of $\sigma$ and $\sigma_{top}$ consists of rows $i$ and above.
            \item Rearrange the entries of the rows of $\sigma_{top}$ in the unique way that forms a filling $\widetilde{\sigma_{top}}$ for which $\inv(\widetilde{\sigma_{top}})=0$.
            \item Apply the procedure of Proposition \ref{ZeroBump} to $\widetilde{\sigma_{top}}$, that is, removing the $n$ from the bottom row and bumping each entry in the last column down one row.  Let the resulting tableau be called $\tau$.
            \item Place $\tau$ on top of $\sigma_{bot}$ and rearrange all rows to form a tableau $\rho$ having $\inv(\rho)=0$.  Then we define $\psi(\sigma)=\rho$.
            \end{enumerate}
             Moreover, if $\maj(\sigma)-\maj(\psi(\sigma))=d$, then $0\le d<\mu_1^\ast$ and we assign $\psi(\sigma)$ to the $d$th set in the disjoint union.
        \end{theorem}

        We first need to prove Lemma \ref{PullUp}, which is a sort of inverse to Lemma \ref{dangle}.
        
        \begin{PullUpLemma}
          Given two collections of letters $b_1,\ldots,b_{w-1}$ and $a_1,\ldots,a_w$, there is a unique element $a_i$ among $a_1,\ldots, a_w$ such that, in any two-row tableau with $a_1,\ldots,\hat{a_i},\ldots,a_w$ as the entries in the bottom row and $b_1,\ldots,b_{w-1},a_i$ as the entries in the top, with no inversions in the top row, the entry $a_i$ occurs in the rightmost position in the top row.
        \end{PullUpLemma}
        
        \begin{proof}
          As usual, let us think of the $a_i$'s as $0$'s and the $b_i$'s as $1$'s in a cocharge word, arranged according to the magnitudes of the $a_i$'s and $b_i$'s.  Then we have a sequence of $w$ $0$'s and $w-1$ $1$'s, and we wish to show that there is a unique $0$ that, when we change it to a $1$, is not paired with any $0$ when computing the crossing number.  By Lemma \ref{dangle}, there is a unique such $1$ in any word of $w-1$ $0$'s and $w$ $1$'s.
          
          So, by Lemma \ref{dangle}, it suffices to find a $0$ in the original tableau such that upon removal, the remaining sequence starting with the entry to its right is $0$-dominated.  For instance, in the sequence $001110100$, which has $5$ zeros and $4$ ones, if we remove the second-to-last zero and cyclically shift the letters so that the new sequence starts with the $0$ to its right, we get the sequence $00011101$, which is $0$-dominated.
          
          To show that there is a unique such $0$, consider the up-down walk starting at $0$ in which we move up one step for each $0$ in the sequence and down one step for each $1$.  Then we end at height $1$, since there is one more $0$ than $1$ in the sequence.  For instance, the sequence $001110100$ corresponds to the up-down walk:
          
          \begin{center}
          \includegraphics{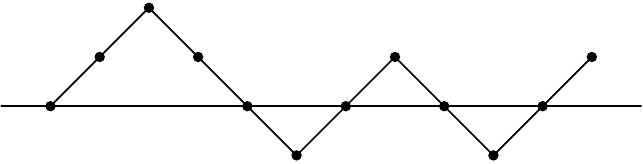}
          \end{center}
          
          Consider the last visit to the minimum height of this walk. If the minimum height is $0$ then we simply remove the last $0$ in the sequence and we are done.  If the minimum height is less than $0$, then there are at least two up-steps ($0$'s) following it since it is the last visit to the min.  The first of these up-steps corresponds to a $0$ which we claim is our desired entry.  Indeed, if we remove this $0$, the walk starting at the next step and cycling around the end of the word is a positive walk, corresponding to a $0$-dominated sequence.  
          
          It is easy to see that if we do the same with any of the other $0$-steps, the resulting walk will not be positive and so the corresponding sequence will not be $0$-dominated.  This completes the proof.
        \end{proof}
        
        We now have the tools to prove Theorem \ref{RectangleBump}.
        
        \begin{proof}[Proof of Theorem \ref{RectangleBump}]
           It is clear that $\psi$ is a morphism of weighted sets, preserving the statistics, so we only need to show that $\psi$ is a bijection.  To do so, we construct an inverse map $\phi=\psi^{-1}$ that takes a pair $(\rho,d)$ and returns an appropriate filling $\sigma:\mu\to\ZZ_+$, where $\rho:\mu^{(d-1)}\to \ZZ_+$ is a filling with no inversions using the letters $1,\ldots,n-1$, and $d$ is a number with $0\le d\le \mu_1^{\ast}-1$.  For simplicity let $h=\mu_1^\ast$ be the height of $\mu$.
           
           Let $(\rho,d)$ be such a pair.  Consider the fillings $\sigma_1,\sigma_2,\ldots,\sigma_h$ formed as follows.  Let $\sigma_h$ be the tableau obtained by inserting the number $n$ into the top row of $\rho$ and rearranging the entries of the top row so that $\inv(\sigma_h)=0$.  Let $\sigma_{h-1}$ be the tableau formed from $\rho$ by first moving the unique element of the $(h-1)$st row given by Lemma \ref{PullUp} to the top row, and then inserting $n$ into the $(h-1)$st row and rearranging all rows so that there are no inversions again.  Then, let $\sigma_{h-1}$ be formed from $\rho$ by first moving the same element, call it $a_{h-1}$, up to the top row, then using Lemma \ref{PullUp} again to move an element $a_{h-2}$ from row $h-2$ to row $h-1$, and finally inserting $n$ in row $h-2$ and rearranging the rows again so that there are no inversions.  Continuing in this manner, we define each of $\sigma_1,\ldots,\sigma_h$ likewise, and it is easy to see that $\psi(\sigma_i)=\rho$ for all $i$, by using Lemma \ref{dangle} repeatedly.
           
           Now, we wish to show that the numbers $d_i=\maj(\sigma_i)-\maj(\rho)$ for $i=1,\ldots,h$ form a permutation of $0,\ldots,h-1$.  Let $a_1,\ldots,a_{h-1}$ be the elements of rows $1,\ldots,h-1$ that were moved up by $1$ in each of the steps as described above.  By Proposition \ref{ZeroBump}, the filling $\sigma_1$, whose rightmost column has entries $a_{h-1},a_{h-2},\ldots,a_1,n$ from top to bottom, has the same major index as $\rho$.  So $d_1=0$, and $\maj(\sigma_1)=\maj(\rho)$.  We will now compare all other $\sigma_i$'s to $\sigma_1$ rather than to $\rho$.
           
           We claim that the difference in the major index from $\sigma_1$ to $\sigma_i$ is the same as the difference obtained when moving $n$ up to row $i$ (and shifting all lower entries down by one) in the one-column filling having reading word $a_{h-1},a_{h-2},\ldots,a_1,n$.  Then, by Carlitz's original bijection, we will be done, since each possible height gives a distinct difference value $d$ between $0$ and $h-1$.  
           
           To proceed, consider the total number of descents in each row.  In $\sigma_i$, the entry $n$ is in the $i$th row.  Let $\tau$ consist of the top $h-i$ rows of this filling, arranged so that $\inv(\tau)=0$.  Then the top $h-i-1$ rows (row $2$ to $h-i$ of $\tau$) are the same as in $\sigma_1$, with the same descents.  Thus if we rearrange \textit{every} row with respect to the one beneath, including rows $i-1$ and below to form $\sigma_i$, each row also has the same number of descents as it does in $\sigma_1$ by Lemma \ref{twist}.
           
           We now show the same is true for row $i+1$.  In $\tau$, we have $a_{i}$ above $n$, and the remaining entries in that row are above the same set of entries they were in $\sigma_1$.  So the number of descents in row $i+1$ goes down by $1$ from $\sigma_1$ to $\sigma_i$ if $a_i>a_{i-1}$, and otherwise it remains the same.
           
           For rows $i$ and below, we use Lemma \ref{OnTop}.  For any row $t$ from $2$ to $i$, the entries of row $t-1$ can be rearranged so that if row $t$ is arranged on top of it with no inversions, the entry $a_t$ lies in the space above $a_{t-1}$ (or $n$ lies above $a_{i-1}$ in the case $t=i$.)  The remaining entries in the top row of this two-row arrangement are then above the same set of entries they were in $\sigma_1$, with no inversions between them, and by Lemma \ref{twist} they have the same number of descents among them.  So, the descents have only changed by what the comparison of each $a_{t}$ with $a_{t-1}$ (or $n$ with $a_{i-1}$) contributes.  
           
           Therefore, the number of descents in a given row of $\sigma_i$, relative to $\sigma_1$, can either increase by $1$, stay the same, or decrease by $1$, according to whether it does in the one-column shape filled by $a_h,\ldots,a_1,n$ when we move $n$ up to height $i$. 
           
           Now, for rectangular shapes, if $p_t$ is the total number of descents in row $t$, it is easy to see that the total cocharge contribution (major index) of the filling is the sum of the partial sums $$p_1+(p_1+p_2)+(p_1+p_2+p_3)+\cdots+(p_1+\cdots+p_h).$$
           Since the values of $p_t$ in $\sigma_i$ differ by $0$ or $\pm 1$ from the corresponding values of $\sigma_1$, it follows that the difference $d_i$ is the sum of the partial sums of these differences.  But this is the same as the difference in the one-column case we are comparing to.  This completes the proof.
        \end{proof}
        
      Note that Proposition \ref{RectBumpCor1} and Theorem \ref{RectBumpCor2} also follow immediately from the above proof.
      
      Finally, we prove Proposition \ref{RectBumpCor3} of section \ref{Rectangles}.
      
      \begin{RectangleBumpCorollary3}
         Let $\mu$ be a rectangle shape of height $h$, and let $\sigma\in \fancy{F}_{\mu}^{(1^n)}$ with its largest entry $n$ in row $i$.  Then if $a_1,\ldots,a_{h-1}$ is the bumping sequence of $\sigma$, then $a_{i+2},\ldots,a_{h-1}$ all occur in columns weakly to the right of the $n$, and each $a_j$ is weakly to the right of $a_{j-1}$ for $j\ge i+3$.
      \end{RectangleBumpCorollary3}

      \begin{proof}
        Let $c_1,\ldots,c_r,n,c_{r+1},\ldots,c_{m-1}$ be the entries in row $i$ from left to right.  Consider the reordering of row $i$ given by $c_1,\ldots,c_{m-1},n$ and order row $i+1$ with respect to this ordering.  Let the numbers in the new ordering in row $i+1$ be $b_1,\ldots,b_{m-1},a_{i}$.  Then $a_{i}$ is the same as the value of $a_{i}$ from Theorem \ref{RectBumpCor2} by Lemma \ref{dangle}; that is, $a_i$ would lie above $n$ if we ordered $c_1,\ldots,c_{m-1}$ by size as well.  
        
        Now, since $c_1,\ldots,c_r$ are the first $r$ entries in both orderings of row $i$, it follows that $b_1,\ldots,b_r$ must be the first $r$ entries in both corresponding orderings of row $i+1$.  Thus $a_i$, not being equal to any of $b_1,\ldots,b_r$, must be weakly to the right of the column that $n$ is in.
        
        The same argument can be used to show that $a_{i+1}$ is weakly to the right of $a_i$ as well, and so on.  This completes the proof.
      \end{proof}
      
    \subsection{Proofs for Three Row Shapes}
    
      We first prove Lemma \ref{threerows}, restating the definition of $\psi$ as part of the statement.
       \begin{ThreeRowsLemma}
         Let $\mu=(\mu_1,\mu_2,\mu_3)$ be any three-row shape of size $n$.  Then there is a  morphism of weighted sets $$\psi:(\fancy{F}_{\mu}^{(1^n)}|_{\inv=0};\maj)\to \bigsqcup_{d=0}^{2}(\fancy{F}_{\mu^{(d+1)}}^{(1^{n-1})}|_{\inv=0};\maj+d)$$ defined combinatorially by the following process.  Given an element $\sigma$ of $\fancy{F}_{\mu}^{(1^n)}|_{\inv=0}$, consider its largest entry $n$.  Let $\sigma'$ be the $3\times\mu_3$ rectangle contained in $\sigma$.
         \begin{enumerate}
            \item If the $n$ is to the right of $\sigma'$, remove the $n$ according to the process in Lemma \ref{tworows} to form $\psi(\sigma)$.
            \item If the $n$ is in the bottom row and in $\sigma'$, then the shape is a rectangle and we remove it according to Theorem \ref{RectangleBump} to obtain $\psi(\sigma)$.
            \item If the $n$ is in the second row and in $\sigma'$, let $a_2$ be the top entry of the bumping sequence of $\sigma'$.  Let $b$ be the entry in square $(\mu_2+1,2)$ if it exists, and set $b=n+1$ if it does not.   If $b\ge a_2$, then remove the $n$ and bump down $a_2$ to the second row.  If $b<a_2$, simply remove the $n$.  Rearrange all rows of the resulting tableau so that there are no inversions, and let $\psi(\sigma)$ be this filling. 
            \item If the $n$ is in the top row and in $\sigma'$, let $a_1,a_2$ be the bumping sequence in $\sigma'$.  If $a_2>a_1$ or $\mu_2=\mu_3$, then simply remove the $n$, and otherwise, if $a_2\le a_1$ and $\mu_2\neq \mu_3$, then remove $n$ and bump $a_2$ up to the top row.  Rearrange all rows so that there are no inversions, and setet $\psi(\sigma)$ to be the resulting filling.
          \end{enumerate}
       \end{ThreeRowsLemma}
       
       \begin{proof}
         We wish to show that $\psi$ is a morphism of weighted sets, i.e. that it preserves the statistics on the objects.  If the $n$ is in the bottom row, then $\psi(n)$ is in the $d=0$ component of the disjoint union and the maj is preserved, by Proposition \ref{ZeroBump}.  If $n$ is in the second row and to the right of column $\mu_3$, then by Lemma \ref{tworows} the difference in maj upon removing it is $1$ and we obtain a filling in the $d=1$ component of the disjoint union.
         
         This leaves us with two possibilities: $n$ is in the second row and weakly to the left of column $\mu_3$, or $n$ is in the top (third) row.  In either case, if $\mu_2=\mu_3$ then the mapping is the same as that in Theorem \ref{RectangleBump}, and we get a map to either the $d=1$ or $d=2$ component of the disjoint union.  So we may assume $\mu_2\neq \mu_3$.  
         
         \textit{Case 1.}  Suppose that $n$ is in the second row.  We have two subcases to consider: $b< a_2$ and $b\ge a_2$.  
         
         If $b<a_2$, $\psi(\sigma)$ is formed by removing the $n$ and rearranging so that there are no inversions.  Note that any entry $i$ to the right of $n$ in row $2$ is less than the entry directly south of $n$.  Furthermore, such entries $i$ are not descents and are increasing from left to right.  Thus these entries simply slide to the left one space each to form $\psi(\sigma)$ after removing the $n$.  So $b$ is the only new entry to be weakly to the left of column $\mu_3$ in $\psi(\sigma)$.  Since $b$ is not a descent, the effect on the major index is the same as if we simply replaced $n$ by $b$ in $\sigma'$.  Consider any arrangement of the second row of $\sigma'$ in which $n$ is at the end, and arrange the top row relative to this ordering.  Then $a_2$ is at the end of this top row by its definition, and so replacing $n$ by $b$ will make $a_2$ a descent and thereby increase the total cocharge contribution by $1$.  By Lemma \ref{twist} this is the same as the increase in the the cocharge contribution from $\sigma$ to $\psi(\sigma)$.  Hence $\psi(\sigma)$ lies in the $d=1$ component of the disjoint union.
         
         If $b\ge a_2$, we claim that if $a_1$ is the entry in the bottom row of the bumping sequence, then $b<a_1$.  If $a_1$ is to the right of the column that $n$ is in then the claim clearly holds.  Otherwise, let $a_1,d_1,d_2,\ldots,d_i$ be the consecutive entries in the bottom row starting from $a_1$ and ending at the entry $d_i$ beneath the $n$, and let $c_1,\ldots,c_i$ be the entries in the second row from the entry above $a_1$ to the entry just before the $n$.  The $c_j$'s are all descents, and the $c_j$'s and $d_j$'s are both increasing sequences.  Since there are no inversions in the second row, we have $b<d_i$.  Since removing the $n$ and bumping up $a_1$ results in the $a_1$ at the end of the second row by definition, upon doing this the $d_i$'s all slide to the left one space, and the $c_i$'s must also remain in position and remain descents by Proposition \ref{ZeroBump}.  In particular, this means that $d_i<c_i$, and so $b<c_i$ as well.  But then since there are no inversions it follows that $b<d_{i-1}$, which is less than $c_{i-1}$, and so on.  Continuing, we find that $b<a_1$ as claimed.
         
         Since $b \ge a_2$ by assumption, it follows that $a_2<a_1$ and so removing the $n$ and bumping down $a_2$ in the rectangle results in a difference in major index of $2$ by Theorem \ref{RectBumpCor2}.  Note also that if we perform this bumping in the entire filling $\sigma$, the entry $a_2$ ends up to the left of column $\mu_3+1$ since $a_2\le b$ and hence it is to the left of $b$ in the second row.  Thus the entries to the right of the rectangle are preserved, and $\maj(\psi(\sigma))=\maj(\sigma)-2$.  It follows that $\psi(\sigma)$ lies in the $d=2$ component of the disjoint union.
         
         \textit{Case 2.} Suppose $n$ is in the top row.  If $a_2>a_1$, then removing $n$ results in the major index decreasing by $d=2$, and so $\psi(\sigma)$ is in the $d=2$ component of the disjoint union.  Otherwise, $a_2\le a_1$. Since $\mu_2\neq \mu_3$, we remove the $n$ and bump $a_2$ up to the top row.
         
         Since $a_2\le a_1$, by Theorem \ref{RectBumpCor2} we find that simply removing the $n$ results in a decrease by $1$ in the major index.  Since the top row has had a descent removed (by the proof of Theorem \ref{RectangleBump}), it follows that the empty space created in the top row was not above a descent, for otherwise the major index would decrease by $2$.  Thus in particular $b$ is not a descent.
         
         It follows that if $\widetilde{\sigma}$ is formed by bumping up $a_2$ and inserting $n$ in the second row, then the $n$, being the last descent in the second row, will appear among the first $\mu_3$ columns of $\widetilde{\sigma}$.  In addition, since $a_2\le a_1$ this results in an \textit{increase} in major index by $1$ from $\sigma$ to $\widetilde{\sigma}$, by Proposition \ref{RectBumpCor1}.  
         
         We now wish to show that $b\ge a_2$; if so, we claim removing $n$ from $\widetilde{\sigma}$ will result in a decrease by $2$ in the major index, and will also result in the tableau $\psi(\sigma)$, thereby showing that $\maj(\psi(\sigma))=\maj(\sigma)-1$ and so $\psi(\sigma)$ is in the $d=1$ component.  To see that the major index decreases by $2$ on removing $n$, note that by Proposition \ref{RectBumpCor2}, the effect of removing the $n$ is the same as replacing $n$ by $b$ in the one-column shape with entries $a_2,n,a_1$.  If $b\ge a_2$ then we have that $b<a_1$ by the same argument as in Case 1 above, and so the major index decreases by $2$.  Thus it suffices to show $b\ge a_2$.  
         
         If $a_2$ is not a descent of $\sigma$, this is clear, so suppose $a_2$ is a descent of $\sigma$ in the second row.  Let $c$ be the entry directly below $a_2$, and assume for contradiction that $b<a_2$.  Then $b<c$, and furthermore the first non-descent in row $2$, say $e$, is less than $c$.  Note that by our above argument we know that $e$ lies within the rectangle $\sigma'$.
         
         Now, we restrict our attention to $\sigma'$ and bump $a_2$ and $a_1$ up one row each, and consider the ordering of the bottom row in which we place $c$ in the column one to the left of the column that $e$ was contained in and shift the remaining entries to the left to fill the row.  Rearranging the new second row with respect to the first, we consider the position of $a_1$ relative to $c$.  If $a_1$ is to the left of the $c$ we have a contradiction since $a_1$ must land in column $\mu_3$ by Lemma \ref{dangle} and the definition of bumping sequence.  Therefore the entries in the second row to the left of $c$ are unchanged.  Since $a_1\ge a_2\ge c$, and all remaining entries in the second row are either $a_1$ or are less than $c$, we have that $a_1$ must be on top of the $c$ in the second row.  This is again a contradiction, since this implies that $a_1$ does not land in column $\mu_3$.  It follows that $b\ge a_2$, as desired.
         
         This completes the proof that $\psi$ is a well-defined morphism of weighted sets.
       \end{proof}
       
       \begin{ThreeRowsBijectiveLemma}
         The map $\psi$ of Lemma \ref{threerows} is an isomorphism.
       \end{ThreeRowsBijectiveLemma}
       
       \begin{proof}
         We know from the lemma above that $\psi$ is a morphism; it suffices to show that it is bijective.  First notice that the cardinality of $\left.\left(\fancy{F}_{\mu}^{(1^n)}\right)\right|_{\inv=0}$ is $$\binom{n}{\mu},$$ and the cardinality of $\left.\left(\bigsqcup_{d=0}^2 \fancy{F}_{\mu^{(d+1)}}^{(1^{n-1})}\right)\right|_{\inv=0}$ is $$\binom{n-1}{\mu_1-1,\mu_2,\mu_3}+\binom{n-1}{\mu_1,\mu_2-1,\mu_3}+\binom{n-1}{\mu_1,\mu_2,\mu_3-1}.$$  Thus the cardinalities of the two sets are equal, and so it suffices to show that $\psi$ is surjective.
         
         To do so, choose an element $\zeta$ of the codomain.  Then $\zeta$ can lie in any one of the three components of the disjoint union $\left(\bigsqcup_{d=0}^2 \fancy{F}_{\mu^{(d+1)}}^{(1^{n-1})}\right)|_{\inv=0}$, and we consider these three cases separately.
         
         \textit{Case 1:} Suppose $\zeta$ lies in the $d=0$ component.   Then we can insert $n$ in the bottom row so as to reverse the map of Proposition \ref{ZeroBump}, and we obtain an element $\sigma$ of $\left(\fancy{F}_{\mu}^{(1^n)}\right)|_{\inv=0}$ which maps to $\zeta$ under $\psi$.
         
         \textit{Case 2:} Now, suppose $\zeta$ lies in the $d=1$ component.  If $\mu_2=\mu_3$ then $\mu^{(1)}=(\mu_1,\mu_2,\mu_3-1)$ and so we can find a filling $\sigma$ of $\mu$ that maps to $\zeta$ by Proposition $\ref{RectangleBump}$.  Otherwise, the shape of $\zeta$ is $(\mu_1,\mu_2-1,\mu_3)$ and we wish to find a filling $\sigma$ of shape $\mu$ for which $\psi(\sigma)=\zeta$.  Let $\rho$ be the filling of $\mu$ formed by inserting $n$ into the second row and rearranging entries so that there are no inversions.  Notice that if the $n$ lies to the right of column $\mu_3$ then $\psi(\rho)=\zeta$ and we are done.
         
         So, suppose $n$ lies in the $3\times \mu_3$ rectangle in $\rho$.  Let $a_1$ and $a_2$ be the bumping sequence of this rectangle.  Since $n$ is the rightmost descent in the second row of $\rho$, inserting it did not change the cocharge contribution of the portion to the right of column $\mu_3$; there were no descents there in $\sigma$ and there are none in $\rho$.  Let $b$ be the entry in column $\mu_3+1$, row $2$ of $\rho$.   If $b<a_2$, then $\psi(\rho)=\zeta$ and we are done.  
         
         Otherwise, if $b\ge a_2$, then by the argument in Lemma \ref{threerows} we know that $\maj(\rho)-\maj(\zeta)=2$.  We have that $\tau:=\psi(\rho)$ is the filling formed by removing the $n$ and bumping $a_2$ down to the second row, and that $\maj(\rho)-\maj(\tau)=2$.  Hence $\maj(\tau)=\maj(\zeta)$.   Since $b\ge a_2$, $a_2$ lies to the left of $b$ in $\tau$ and hence is weakly to the left of column $\mu_3$.  So, let $\sigma$ be the tableau formed by inserting $n$ in the top row of $\tau$.  Now $\sigma$ has shape $\mu$, and can be formed directly from $\rho$ by shifting the position of $n$ among $a_1$ and $a_2$ as in Theorem \ref{RectBumpCor2}.  
         
         It follows that $\maj(\sigma)-\maj(\rho)=\pm 1$, and so $\maj(\sigma)-\maj(\tau)$ is equal to $1$ or $3$.  It is not $3$ because $\tau$ is formed from $\sigma$ by removing an $n$ from the top row, which changes the major index by at most $2$ by Theorem \ref{RectangleBump}.  It follows that $\maj(\sigma)-\maj(\tau)=1$, and therefore $a_2\le a_1$ by Theorem \ref{RectBumpCor2}.  Thus $\psi(\sigma)=\zeta$ by the definition of $\psi$.
           
          \textit{Case 3:} Suppose $\zeta$ is in the $d=2$ component.  If $\mu_2=\mu_3$ then we simply insert $n$ into $\zeta$ in either row $2$ or $3$ according to Theorem \ref{RectBumpCor2} to obtain a tableau $\sigma$ with $\psi(\sigma)=\zeta$.
          
          Otherwise, if $\mu_2\neq \mu_3$, $\zeta$ has shape $(\mu_1,\mu_2,\mu_3-1)$.  Let $\rho$ be the tableau of shape $\mu$ formed by inserting $n$ in the top row of $\zeta$.  Let $a_1$ and $a_2$ be the entries in row $1$ and $2$ corresponding to this $n$ in the $3\times \mu_3$ rectangle contained in $\rho$.  Then if $a_2>a_1$, $\psi(\rho)=\zeta$ and we're done.
          
          If instead $a_2\le a_1$, then removing $n$ from $\rho$ decreases its major index by $1$.  Since the number of descents in the top row goes down by exactly $1$ by Lemma \ref{tworows}, we can conclude that the entry in row $2$, column $\mu_3$ is a non-descent; otherwise removing $n$ from $\rho$ would decrease the major index by $2$.  So, let $\sigma$ be the filling formed by removing $n$ from $\rho$, bumping $a_2$ to the top row, and inserting $n$ in the second row.  Since there are non-descents in the rectangle we have that $n$ lies in the rectangle in $\sigma$ as well.
          
          Finally, again by the argument used for Lemma \ref{threerows} we have that $a_2\le b$ where $b$ is the entry in row $2$, column $\mu_3+1$ in $\sigma$.  Thus $\psi(\sigma)=\zeta$ as desired.
       \end{proof}
       
       We can now complete the three-row case by using the standardization map $\Standardize$ defined in Section \ref{ThreeRowsSection} for fillings with repeated entries.  We first state a structure lemma about three-row shapes with no inversions.
       
       \begin{lemma}\label{patterns}
         If the consecutive entries $b_1,\ldots,b_n$ in some row of a filling with no inversions are directly above a weakly increasing block of squares $c_1\le \cdots \le c_n$ in the row below, then there exists a $k$ for which $b_1,\ldots,b_k$ are descents and $b_{k+1},\ldots,b_{n}$ are not descents.  Moreover $b_1\le \cdots \le b_k$ and $b_{k-1}\le \cdots \le b_n$ are both increasing blocks of squares.
       \end{lemma}
       
       \begin{proof}
          This is clear by the definition of inversions.  
       \end{proof}
              
        In particular, the second row has one (possibly empty) block of descents and one (possibly empty) block of non-descents.  The third row has up to two blocks of descents, one for each of the blocks in the second row, and so on.
       
       We also need to show that the cardinalities of the sets are equal in the case of repeated entries.
       
       \begin{lemma}\label{Cardinalities}
        We have $$\left|\fancy{F}_\mu^\alpha|_{\inv=0}\right|=|C_{\mu,A}|$$ for any alphabet $A$ with content $\alpha$ and any shape $\mu$.
       \end{lemma}
       
       \begin{proof}
        Given an alphabet $A$, the cocharge word of any filling using the letters in $A$ has the property that it is weakly increasing on any run of a repeated letter, where we list the elements of $A$ from largest to smallest.  Furthermore, the cocharge word has content $\mu$.  It is not hard to see that a word is the cocharge word of a filling in $\fancy{F}_\mu^\alpha|_{\inv=0}$ if and only if it has content $\mu$ and is weakly increasing over repeated letters of $A$, listed from greatest to least.
        
        Recall that the fillings in $\fancy{F}_{\mu^\ast}^{r(\alpha)}|_{\maj=0}$ can be represented by their \textit{inversion word}, and a word is an inversion word for such a filling if and only if it has content $\mu$ and every subsequence corresponding to a repeated letter of the reversed alphabet is in inversion-friendly order.  By swapping the inversion-friendly order for weakly increasing order above each repeated letter, we have a bijection between inversion words and cocharge words, and hence a bijection (of sets, not of weighted sets) from $\fancy{F}_\mu^\alpha|_{\inv=0}$ to $\fancy{F}_{\mu^\ast}^{r(\alpha)}|_{\maj=0}$.  By Theorem \ref{Invcode}, we have that $$\left|\fancy{F}_{\mu^\ast}^{r(\alpha)}|_{\maj=0}\right|=|C_{\mu,A}|,$$ and so the cardinality of $\fancy{F}_\mu^\alpha|_{\inv=0}$ is equal to $|C_{\mu,A}|$ as well.
       \end{proof}
       
      \begin{ThreeRowsConclusionTheorem}
        The map $\majcode$ defined on three-row shapes $\mu=(\mu_1,\mu_2,\mu_3)$ is an isomorphism of weighted sets $$\fancy{F}_{\mu}^\alpha|_{\inv=0}\to C_{\mu,A}$$ for each alphabet $A$ and corresponding content $\alpha$.  
      \end{ThreeRowsConclusionTheorem}
      
      \begin{proof}
       By Lemmas \ref{threerows}, \ref{threerowsbijective}, and \ref{recursion}, we have that for the content $(1^n)$ corresponding to alphabet $[n]$, the map $\majcode$ is a weighted set isomorphism $$\fancy{F}_{\mu}^{(1^n)}|_{\inv=0}\to C_{\mu,[n]}.$$
                  
       Now, let $A$ be any alphabet with content $\alpha$.  Let $\sigma$ be a filling of $\mu$ with content $\alpha$.  Then we know $\majcode(\sigma)=\majcode(\Standardize(\sigma))$, so $\majcode(\sigma)\in C_{\mu,[n]}$.  In other words, $\majcode(\sigma)$ is $\mu$-sub-Yamanouchi.  In addition, since $\Standardize$ is an injective map (there is clearly a unique way to un-standardize a standard filling to obtain a filling with a given alphabet), the map $\majcode$, being a composition of $\Standardize$ and the $\majcode$ for standard fillings, is injective as well on fillings with content $\alpha$.
                  
       We now wish to show that $\majcode(\sigma)=d_1,\ldots,d_n$ is $A$-weakly increasing, implying that $\majcode$ is an injective morphism of weighted sets to $C_{\mu,A}$.  By Lemma \ref{Cardinalities} this will imply that it is an isomorphism of weighted sets.  It suffices to show this for the largest letter $m$ of $A$ by the definition of standardization.  Suppose $m$ occurs $i$ times.   We wish to show that $d_j\le d_{j+1}$ for all $j\le i-1$.  So choose $j\le i-1$.
       
       Suppose $d_j=0$.  Then by the definition of $\Standardize$, we have that the $m$ we removed from $\psi^{j-1}(\sigma)$ was in the bottom row.  If there are still $m$'s in the bottom row of $\psi^{j}(\sigma)$ then $d_{j+1}=0$ as well.  Otherwise $d_{j+1}>0$, so $d_j\le d_{j+1}$ in this case.
       
       Suppose $d_j=1$.  Then the $m$ we removed from $\psi^{j-1}(\sigma)$ was in either the first or second row and there were no $m$'s in the bottom row.  By the definition of $\psi$, there are therefore no $m$'s in the bottom row of $\psi(\sigma)$ either, and so $d_{j+1}\ge 1=d_j$.
       
       Finally, suppose $d_j=2$.  Let $m_j$ be the $m$ we remove from $\psi^{j-1}(\sigma)$ to obtain $d_j=2$.  As in the previous case we have $d_{j+1}\ge 1$, and we wish to show $d_{j+1}\neq 1$.  Let $m_{j+1}$ be the corresponding $m$.  Since $d_j$ is minimal for $\psi^{j-1}(\sigma)$, there are no $m$'s in $\psi^{j-1}(\sigma)$ which we can treat as the largest entry and remove according to $\psi$ to form $d_j=1$.  Therefore if we removed $m_{j+1}$ before $m_j$ we would also have a difference of $2$ in the major index.
        
       We consider three subcases separately for the locations of $m_j$ and $m_{j+1}$: they can either both be in the second row, $m_j$ can be in the second row with $m_{j+1}$ in the third (top) row, or they can both be in the top row.  No other possibilities exist because they must occur in reverse reading order, and cannot be in the bottom row since $d_j=2$.
       
       \textit{Subcase 1:}  Suppose both $m_j$ and $m_{j+1}$ are in the second row.  Then $m_{j+1}$ and $m_{j}$ are at the end of the block of descents in that order, and weakly to the left of column $\mu_3$.  Let $b$ be the entry in row $2$, column $\mu_3+1$.  Let $a_2$ be the entry in the third row in the bumping sequence of $m_j$, and let $a_2'$ be the entry in the bumping sequence of $m_{j+1}$ in $\psi^j(\sigma)$.  Since $d_j=2$, we have $a_2\le b$ and $b< m$, and so $a_2\neq m$.  Therefore no new $m$'s are dropped down.  In other words, $m_{j+1}$ is indeed the $m$ that will be removed upon applying $\psi$ the second time.
       
       We now need to check that $m_{j+1}$ remains to the left of column $\mu_3$ after applying $\psi$.  Indeed, by Proposition \ref{RectBumpCor2}, we have that the number of descents in row $2$ goes down by one, and the number of descents in the top row remains the same, upon applying $\psi$ to $\psi^{j-1}(\sigma)$.  Since there are no $m$'s in the bottom row, $m_{j+1}$ is the rightmost descent in the second row of $\psi^{j}(\sigma)$, and the descent we lost was $m_j$, so $m_{j+1}$ remains in its column.
       
       We now just need to show that $a_2'\le b'$, where $b'$ is the entry in row $2$, column $\mu_3$ after applying $\phi$.  Either $b'=b$, $b'=a_2$, or $b'$ is the entry $b_0$ that is bumped out from the first $\mu_3-2$ columns when we drop down $a_2$.  
       
       Consider any ordering of the first $\mu_3$ entries of the second row of $\psi^{j-1}(\sigma)$ such that the two $m$'s ($m_{j+1}$ and $m_j$) are at the end in that order, and also place $b_0$ in the third-to-last position.  Now, rearrange the entries above these so that there are no inversions.  We know that $a_2$ is at the end of the top row, above $m_j$, by its definition.  Let $a$ be the entry above $m_{j+1}$ and let $c$ be the entry to the left of that (if such a column exists.)
       
       If $b'=a_2$, then the first $\mu_3-2$ entries of the second row are unchanged on applying $\psi$.  In our new ordering above, this means $a_2'=a$, and since $a$ and $a_2$ occur above the two $m$'s in our new ordering, we have $a\le a_2$.  It follows that $a_2'\le b'$.
       
       If $b'\neq a_2$, then $b'$ is either $b$ or $b_0$.  To find $a_2'$ in our new ordering, bumping down the $a_2$ can be thought of as replacing the $b_0$ with $a_2$ and rearranging the top row again so that there are no inversions.  The first $\mu_3-3$ entries remain in the same positions, and either $c$ or $a$ lies above $m_{j+1}$ based on which comes later in cyclic order after $a_2$.  So either $a_2'=a$ or $a_2'=c$.  
       
       We now have to show that whether $a_2'$ is $a$ or $c$, it is less than both $b$ and $b_0$.  Notice that $a_2\le b_0$: Since $m_{j+1}$ stays in its place, either $a_2$ replaces a larger entry among the descents in the second row, which in turn bumps out a larger entry $b_0$ among the non-descents, or it replaces a non-descent itself and displaces a larger non-descent $b_0$ to its right.  So if $a_2'=a$, then since $a\le a_2$ we have $a_2'\le a_2\le b_0$ and also $a_2'\le a_2\le b$ since $a_2\le b$.
       
       Finally, if $a_2'=c$, then $a_2,a, c$ are in cyclic order.  If $c\le a_2$ we are done by the above argument.  Otherwise $a_2< a\le c$ or $a_2=a=c$, in which case $a_2'=a$ and we are done by the previous case.  So $a_2<a\le c$, but we already know $a\le a_2$, so we have a contradiction.  It follows that $a_2'\le b'$ as desired.
       
       \textit{Subcase 2:}  Suppose $m_j$ and $m_{j+1}$ are in the top row.  Then by Lemma \ref{patterns} and since there are no $m$'s in the second row by the definition of $\Standardize$, the $m$'s are either in the first or second block of descents in the third row.  If either of them is in the second block, it is clear that removing $m_j$ results in $d_j=1$, not $2$, a contradiction.  So they are both in the first block, themselves above descents in the second row, with $m_{j+1}$ and $m_j$ adjacent and in that order.
       
       Now, removing $m_j$ will cause the block of non-descents to its right to slide to the left one space (since they are necessarily less than the entry beneath $m_j$).  If the \textit{second} block of non-descents in the third row is nonempty, one of these will replace the last entry above the descents in the second row, since all of these are still less than the entry below $m_j$ and the least among the entries to the right will replace it.  In that case the number of descents to the right of $m_j$ is unchanged, and so $d_j=1$, a contradiction.  Thus there are no non-descents in the second block, i.e. above the non-descents in row $2$.
       
       Because of this, removing $m_j$ simply causes all the entries to its right to slide to the left one space, and the first descent to its right becomes a non-descent.  The same then happens when we remove $m_{j+1}$ by the same argument.  It follows that $d_{j+1}=2$ in this case.
       
       \textit{Subcase 3:}  Suppose $m_j$ is in the second row and $m_{j+1}$ in the top.  Then the $m_{j+1}$ is to the left of $m_j$, in the first block of descents in the third row, since otherwise we would have a difference of $1$ on removing $m_{j+1}$.  Moreover, as in the previous case, the top row has no non-descents above the non-descents in row $2$.
       
       So, let $a_1,\ldots, a_r$ be the entries in row $3$ that lie weakly to the right of $m_j$'s column.  Then $a_1$ is not a descent and each of $a_2,\ldots,a_r$ are descents.  Let $m_j,b_2,\ldots,b_r$ be the entries below them.  If we rearrange these in the second row in the increasing order $b_2,\ldots,b_r,m_j$, and then rearrange the $a_i$'s above them as $a_{\sigma(1)},\ldots,a_{\sigma(r)}$ so that there are no inversions, there are still $r-1$ descents among the $a_{\sigma(i)}$'s by Lemma \ref{twist}.  These descents must be $a_{\sigma(1)},\ldots,a_{\sigma(r-1)}$ by Lemma \ref{patterns}, and the last entry $a_{\sigma(r)}$ above the $m_j$ is the entry in $m_j$'s bumping sequence.
       
       Now, to form $\psi^j(\sigma)$, we remove $m_j$ and drop down $a_{\sigma(r)}$.  Notice that the entries in the top row to the left of where $m_j$ was are unchanged: consider the $3\times \mu_3$ rectangle and bump the $m_j$ down to the bottom row according to Theorem \ref{RectangleBump}.  Then bump it out according to Proposition \ref{ZeroBump}, which leaves us with the same top row as that of $\psi^j(\sigma)$.  The entire top row save for the last entry is unchanged upon applying Proposition \ref{ZeroBump}, and so having the $m_j$ inserted into the second row instead can only change the entries to the right of it in the row above.  Thus the entries to its left in the top row are unchanged, and have the same cocharge contribution as well.
       
       Finally, in the columns weakly to the right of the column that $m_j$ was in, the entries in the top row are $a_{\sigma(1)},\ldots, a_{\sigma(r-1)}$ in some order.  We claim that the entries in the second row are formed by replacing at most one of $b_2,\ldots,b_r$ by a smaller entry, which is either $a_{\sigma(r)}$ or something bumped to the right by $a_{\sigma(r)}$ if $a_{\sigma(r)}$ lands in a column to the left of the $b_i$'s.  Indeed, the only way it would be a larger entry replacing them is if a descent replaced $m_j$, but in this case we would have $d_j=1$ since the number of descents in the second row would be the same, and the number of descents in the top row would decrease by only $1$. 
       
       Therefore, the entries $a_{\sigma(1)},\ldots,a_{\sigma(r-1)}$ are all descents in the top row, and so removing $m_{j+1}$ still results in a difference $d_{j+1}=2$.  In particular, the descents formed by $m_{j+1}$ and one of the $a_i$'s are removed, since the $a$'s all slide one position to the left, and did not form new descents upon removing the $m_{j+1}$ before the $m_j$.
       
       This completes the proof.   
      \end{proof}
       
      \begin{corollary}
         For any three-row shape $\mu$ and content $\alpha$, the map $\invcode^{-1}\circ\majcode$ is a weighted set isomorphism from $\fancy{F}_{\mu}^{\alpha}|_{\inv=0}\to \fancy{F}_{\mu^\ast}^{r(\alpha)}|_{\maj=0}$.  This gives a combinatorial proof of the identity $$\widetilde{H}_{\mu}(x;0,t)=\widetilde{H}_{\mu^\ast}(x;t,0)$$ for three-row shapes.
      \end{corollary}
    
    \section{Acknowledgments}
    
    This work was partially funded by the Hertz foundation and the NSF Graduate Research Fellowship Program.
    
    The author thanks Mark Haiman and Angela Hicks for their guidance and support.

\end{document}